\documentclass[10pt,a4paper]{article}

\usepackage[a4paper,margin=2.5cm]{geometry}

\usepackage{amsmath,amssymb,amsthm,mathtools}
\usepackage{mathabx,stmaryrd,bm}
\input{Ustmry.fd}
\DeclareFontShape{U}{stmry}{b}{n}{<->ssub*stmry/m/n}{}

\usepackage{array,booktabs,enumitem}
\usepackage{graphicx,calc,xcolor,tikz}

\usepackage{xurl}
\usepackage[breaklinks,bookmarks=false]{hyperref}
\hypersetup{colorlinks=true,linkcolor=blue,citecolor=blue,urlcolor=blue,
  plainpages=false,pdfwindowui=false,pdfstartview={FitH}}

\newtheorem{lemma}{Lemma}[section]
\newtheorem{theorem}[lemma]{Theorem}
\newtheorem{proposition}[lemma]{Proposition}
\newtheorem{corollary}[lemma]{Corollary}
\newtheorem{assumption}[lemma]{Assumption}
\newtheorem{remark}[lemma]{Remark}
\newtheorem{example}[lemma]{Example}

\newcommand{\email}[1]{\href{mailto:#1}{#1}}
\newcommand{\HHOabstractstart}{\noindent}
\newcommand{\HHOmetadata}{
  \bigskip\noindent\textbf{Keywords: }\KW\par
  \noindent\textbf{MSC2020:} 35J60, 65N12, 65N15, 65N30.}

\numberwithin{equation}{section}

\def\bae#1\eae{\begin{equation}\begin{aligned}#1\end{aligned}\end{equation}}
\def\ba#1\ea{\begin{align}#1\end{align}}
\def\ban#1\ean{\begin{align*}#1\end{align*}}
\newcommand{\bse}{\begin{subequations}}
\newcommand{\ese}{\end{subequations}}
\newcommand{\bt}{\begin{theorem}}
\newcommand{\et}{\end{theorem}}
\newcommand{\bpr}{\begin{proposition}}
\newcommand{\epr}{\end{proposition}}
\newcommand{\bl}{\begin{lemma}}
\newcommand{\el}{\end{lemma}}
\newcommand{\bc}{\begin{corollary}}
\newcommand{\ec}{\end{corollary}}
\newcommand{\bp}{\begin{proof}}
\newcommand{\ep}{\end{proof}}
\newcommand{\br}{\begin{remark}}
\newcommand{\er}{\end{remark}}
\newcommand{\bas}{\begin{assumption}}
\newcommand{\eas}{\end{assumption}}

\newcommand\Om{\Omega}

\newcommand\eq{\coloneqq}

\newcommand\tL{\bm{L}}
\newcommand\tW{\bm{W}}

\newcommand{\bxi}{{\bm\xi}}
\newcommand\btau{\bm{\tau}}

\newcommand\F{\mathcal{F}}
\newcommand\I{\mathcal{I}}
\newcommand\T{\mathcal{T}}

\newcommand\RR{{\mathbb R}}
\newcommand\NN{{\mathbb N}}

\newcommand\PP{{\mathbb P}}

  \newcommand{\restri}[1]{\ \!\!_{|_{#1}}}

\newcommand\refover[2]{
\setlength{\dimen0}{0pt-(\widthof{$\overset{#1}{#2}$}-\widthof{$#2$})/2}\hspace{\dimen0}\overset{#1}{#2}
}
\newcommand\refovern[2]{
\overset{#1}{#2}
}

\DeclareMathOperator*{\essinf}{ess\,inf}
\DeclareMathOperator*{\esssup}{ess\,sup}

\newcommand*\CV[2]{\,\smash{\mathop{\xrightarrow{\hspace*{0.8cm}}}\limits_{#1 \to #2}}\ \vspace{0.07cm}}

\def\b{\boldsymbol}
\newcommand*\cst[1]{\mathrm{#1}}

\newcommand\dU[2]{{\u U}_{#1}^{#2}}

\renewcommand\u{\underline}
\newcommand{\Fb}{\F_h^\cst{b}}

\newcommand{\res}[1]{\ \!\!_{|_{#1}}}
\newcommand{\dgrad}[2]{\b{\cst{G}}^{#1}_{#2}}
\newcommand{\drec}[2]{\cst{R}^{#1}_{#2}}
\newcommand{\dfbres}[2]{\Delta^{#1}_{#2}}
\newcommand{\GRAD}{\b\nabla}
\newcommand{\DIV}{\b\nabla{\cdot}}
\newcommand\Interp[2]{\u{I}_{#1}^{#2}}
\newcommand\proj[2]{\pi_{#1}^{#2}}
\newcommand\PROJ[2]{\b{\pi}_{#1}^{#2}}
\newcommand\stress{\b\sigma}

\newcommand\stabfl{\sigma}

\newcommand\sob{p}
\newcommand\sobmax{\widehat{\sob}}
\newcommand\sobmin{\widecheck{\sob}}
\newcommand\petp{{\scalebox{0.5}{$+$}}}
\newcommand\petm{{\scalebox{0.5}{$-$}}}

\newcommand\Plog{\mathcal P^\cst{log}}
\newcommand\Plogm{\mathcal P^\cst{log}_-}
\newcommand\Plogp{\mathcal P^\cst{log}_+}
\newcommand\Plogpm{\mathcal P^\cst{log}_{\pm}}

\title{A Hybrid High-Order method for variable-exponent Leray--Lions problems with convergence under minimal regularity}
\author{Andr\'e Harnist\footnotemark[1]$~^{1}$\and Romain Mottier\footnotemark[2]$~^{2}$}
\def\KW{variable-exponent Leray--Lions equations, Hybrid High-Order methods, discrete compactness, convergence under minimal regularity, polytopal meshes}

\begin{document}
\maketitle

\footnotetext[1]{\email{andre.harnist@utc.fr}, \url{https://andreharnist.fr}, corresponding author}
\footnotetext[2]{\email{romain.mottier@outlook.com}, \url{https://romainmottier.github.io}}
\renewcommand{\thefootnote}{\fnsymbol{footnote}}
\footnotetext[1]{Université de technologie de Compiègne, Alliance Sorbonne Université, LMAC, Compiègne, France}
\footnotetext[2]{University of Basel, Department of Mathematics and Computer Science, Basel, Switzerland}
\renewcommand{\thefootnote}{\arabic{footnote}}

\begin{abstract}
\HHOabstractstart
{
We consider a homogeneous Dirichlet problem for nonlinear elliptic equations in divergence form with Leray--Lions fluxes of variable growth. The class includes the variable-exponent $p(\cdot)$-Laplacian and regularized Carreau--Yasuda-type laws. We propose and analyze a Hybrid High-Order (HHO) discretization on general polytopal meshes, a hybrid nonconforming method based on cell and face polynomial unknowns of arbitrary degree. We develop a variable-exponent polynomial and discrete functional toolkit centered on a mean-exponent norm equivalence, together with inverse, trace, Poincar\'e--Wirtinger, Sobolev--Poincar\'e, discrete Lebesgue embedding, interpolation, and compactness estimates adapted to HHO unknowns. We then establish well-posedness of the nonlinear discrete problem and convergence without assuming additional regularity of the exact solution beyond its natural energy space. The discrete potential converges strongly in every uniformly subcritical variable-exponent Lebesgue space, the reconstructed gradient converges strongly in the natural variable-exponent space, and the stabilization energy vanishes. Numerical experiments in two, three, and four space dimensions confirm convergence for singular solutions and for exponent fields ranging from close to $1$ up to $6$. Quantitative a priori error estimates for the same method are derived in the companion paper~\cite{harnist-romain-a-priori}. 
}
\end{abstract}

\HHOmetadata

\section{Introduction}
\label{intro}

We consider the homogeneous Dirichlet problem
\bse\label{eq:intro.problem}
\begin{alignat}{2}
-\operatorname{div}\bigl(\stress(\cdot,\nabla u)\bigr)&=f &&\quad\text{in }\Omega,\\
u&=0 &&\quad\text{on }\partial\Omega,
\end{alignat}
\ese
where $\Omega\subset\mathbb R^d$ is a bounded polytopal Lipschitz domain and $\stress:\Omega\times\mathbb R^d\to\mathbb R^d$ is a Leray--Lions flux whose growth exponent $p(\cdot)$ varies in space. The $p(\cdot)$-Laplacian is the basic prototype (see Example~\ref{ex:laplace}), and the class considered below also contains regularized Carreau--Yasuda-type constitutive laws with spatially dependent parameters. Such variable-growth models arise in the description of electrorheological and other non-Newtonian fluids and in variational models for image restoration; see, for example, \cite{Ruzicka2000,ChenLevineRao2006,diening2011lebesgue}.

The analytical framework for these problems is provided by Lebesgue and Sobolev spaces with variable exponents, whose basic properties (duality, reflexivity, density, and Sobolev embeddings) were established in \cite{KovacikRakosnik1991,FanZhao2001,FanShenZhao2001} and are presented systematically in \cite{diening2011lebesgue}; a broad account of differential equations with nonstandard growth is given in \cite{HarjulehtoHastoLeNuortio2010}. For elliptic equations and systems with $p(\cdot)$-growth, existence \cite{FanZhang2003,Fan2011}, regularity and gradient estimates \cite{AcerbiMingione2001,AcerbiMingione2005,Fan2007,DelPezzoMartinezH2}, and structural stability for variable-exponent Leray--Lions problems \cite{AndreianovBendahmaneOuaro2010} are now well understood. A feature of this theory that is central here is the role of the regularity of the exponent: under log-H\"older continuity of $p(\cdot)$, smooth functions are dense in the variable-exponent Sobolev space and the usual embeddings hold, whereas for less regular exponents Lavrentiev-type phenomena may occur \cite{Zhikov1987,Zhikov1995,Zhikov2011}.

The numerical analysis of problems with variable growth is more recent. For the stationary $p(\cdot)$-Laplacian, an interior-penalty discontinuous Galerkin (dG) method with convergence of discrete minimizers was developed in \cite{DelPezzoLombardiMartinez2012}, which also provides broken variable-exponent Sobolev and compactness tools used below. Conforming finite element error estimates for the Dirichlet problem were obtained in \cite{BreitDieningSchwarzacher2015,DelPezzoMartinez2015}, and medius and a priori estimates for a Crouzeix--Raviart discretization in \cite{BalciKaltenbach2025}. Finite element estimates for a variable-exponent Stokes model of electrorheological fluids were established in \cite{BerselliBreitDiening2016}, nonconforming discretizations of non-autonomous variational problems including variable-exponent energies were studied in \cite{BalciOrtnerStorn2022}, and control-volume finite element and discrete-duality finite-volume schemes for variable-exponent nonlinear diffusion were analyzed in \cite{AndreianovQuenjel2023}. Evolution problems have been addressed by space--time finite element approximations of variable-exponent wave equations and of generalized Newtonian flows with nonstandard growth \cite{HaehnleProhl2010,CarelliHaehnleProhl2010}. On the algorithmic side, the resulting nonlinear systems have been solved by decomposition--coordination \cite{DelPezzoMartinez2014} and quasi-Newton minimization strategies \cite{CaliariZuccher2017,AragonBonderRubio2023}. These analyses mostly address the $p(\cdot)$-Laplacian or problems with a variational structure.

Hybrid nonconforming methods attach polynomial unknowns both to the mesh cells and to the mesh faces; this family includes, among others, the hybridizable discontinuous Galerkin (HDG)~\cite{CockburnGopalakrishnanLazarov2009}, weak Galerkin (WG)~\cite{WangYe2013}, and Hybrid High-Order (HHO)~\cite{DiPietroErnLemaire2014,DiPietroErn2015} methods, which are closely related~\cite{CockburnDiPietroErn2016}. The two design principles of HHO methods are: (i) a reconstruction of the gradient in each mesh cell from the local cell and face unknowns; (ii) a stabilization that weakly enforces the matching of the trace of the cell unknowns with the face unknowns. We refer to the monographs~\cite{Di-Pietro.Droniou:20,CicuttinErnPignet2021} for a detailed presentation. HHO methods support arbitrary polynomial degrees and general element shapes, are locally conservative, and their cell unknowns can be eliminated by static condensation, so that the global system only couples face unknowns. Beyond linear problems, HHO methods have been applied, e.g., to nonlinear elasticity~\cite{BottiDiPietroSochala2017}, hyperelasticity~\cite{AbbasErnPignet2018}, non-Newtonian flows~\cite{Botti.Castanon-Quiroz.ea:20,harnist3}, and nonlinear acoustic waves~\cite{SteinsErnJamondDrui2023}.

For Leray--Lions equations with a constant growth exponent, convergence under minimal regularity and the associated discrete functional analysis were established in \cite{Di-Pietro.Droniou:17}, approximation properties and a priori estimates in \cite{Di-Pietro.Droniou:17*1}, and refined estimates for degenerate regimes in \cite{harnist2}; HHO formulations for constant-exponent convex minimization and $p$-Laplace problems can also be found in \cite{carstensen:20,RenLiu2024}.

Virtual Element Methods (VEM)~\cite{Beirao-da-Veiga.Brezzi.ea:13} form another family of polytopal discretizations. VEM have been analyzed for non-Newtonian Stokes systems with Leray--Lions-type constitutive laws, with convergence analysis and a priori error estimates in the shear-thinning and shear-thickening regimes~\cite{AntoniettiBeiraoBottiVaccaVerani2024,harnist5}. More recently, a mixed VEM for the scalar $p$-Laplace equation with constant exponent $1<p<\infty$ was introduced in~\cite{HaileVacca2026}, with well-posedness and a priori error estimates for the solution and the flux.

To the best of our knowledge, however, no hybrid nonconforming method, such as HDG, WG, or HHO, has been analyzed for variable-exponent Leray--Lions problems. In particular, no convergence result under minimal regularity is available for a hybrid discretization of arbitrary order on polytopal meshes that covers non-variational Leray--Lions fluxes with variable growth. This is the purpose of the present paper.

Extending the constant-exponent analysis requires accommodating the spatial variation of the growth exponent at the mesh scale. Local polynomial inverse and scaling estimates are naturally expressed with a fixed exponent on each cell, whereas the continuous energy is measured with the pointwise exponent, so that their uniform comparison involves the oscillation of $p(\cdot)$ on each cell. The log-H\"older continuity of the exponent provides the scale control needed for this comparison. Moreover, since no regularity of the exact solution beyond its energy space is assumed, the convergence argument relies on modular bounds and compactness rather than on consistency estimates involving higher derivatives of the solution.

The contributions of this work are as follows. First, we develop a polynomial toolkit that transfers fixed-exponent arguments to variable exponents uniformly at the mesh scale. Its central result is the mean-exponent equivalence of Theorem~\ref{thm:mean_exp_equiv}, and in particular for every $U\subset\Omega$ satisfying the geometric assumptions of that theorem,
\bae
\begin{gathered}
 \|v\|_{L^{p(\cdot)}(U)}\eqsim \|v\|_{L^{p_U}(U)}, \qquad \forall\, v\in\PP^k(U),\\[5pt]
 \text{with }p_U\coloneq \frac{1}{|U|}\int_U p(\b x)\,\cst{d}\b x,
\end{gathered}
\eae
cf.~\eqref{eq:mean.exp.norm}. This norm equivalence is an immediate consequence of the finite-dimensionality of $\PP^k(U)$. The essential point of Theorem~\ref{thm:mean_exp_equiv} is the uniform control of the hidden constants in $\eqsim$. Indeed, they depend only on the fixed parameters listed in the theorem, and not on the mesh size or on any element of the mesh. This equivalence is the basic device of the analysis: it yields the polynomial scaling rules needed for variable-exponent inverse estimates, the uniform boundedness of local orthogonal projectors, and, in its two-exponent form, the discrete Lebesgue embedding of Corollary~\ref{cor:discrete-lebesgue-emb}. Related reconstruction and embedding estimates are a central ingredient of the interior-penalty dG analysis of \cite{DelPezzoLombardiMartinez2012}; the mean-exponent formulation used here is tailored to polynomial HHO unknowns, in that mesh-scaling factors are expressed through local mean exponents rather than local endpoint exponents.

Second, we establish the discrete functional analysis required by the HHO setting: the equivalence of the natural HHO seminorms, variable-exponent inverse and trace inequalities, and discrete Poincar\'e--Wirtinger and Sobolev--Poincar\'e inequalities covering the full variable-exponent range of the method, the latter relying on a variable-exponent bound for the node-averaging operator. We also state the variable-exponent approximation properties of the local polynomial projectors and of the HHO interpolator in Section~\ref{sec:properties}; the polynomial-projector proofs are given in the companion a priori analysis~\cite{harnist-romain-a-priori}, where these estimates are also used to derive quantitative error bounds, while the interpolator boundedness proof is supplied here. This organization reflects the complementary objectives of the two papers: the present work focuses on qualitative convergence, while the companion paper develops the approximation theory underlying the quantitative error estimates.

Third, we construct and analyze an HHO discretization for general variable-exponent Leray--Lions fluxes satisfying continuity and strong-monotonicity assumptions, without requiring a variational structure. The method is posed on general polytopal meshes, allows every polynomial degree $k\ge0$, covers both equal-order and mixed-order cell--face spaces, and uses a nonlinear face stabilization that may carry its own degeneracy parameter. We establish the commutation and polynomial-consistency properties of the reconstructions, the strong monotonicity of the continuous and discrete diffusion operators, and the well-posedness of the discrete problem together with uniform a priori bounds.

Finally, we prove a discrete compactness theorem and convergence under the sole energy regularity of the exact solution. Bounded sequences in the natural HHO norm admit subsequences whose discrete potentials converge strongly in every uniformly subcritical variable-exponent Lebesgue space and whose reconstructed gradients converge weakly in the natural space. For discrete solutions, monotonicity and a Minty argument upgrade the gradient convergence to strong convergence, and the stabilization energy vanishes. The convergence analysis thus requires no higher-order regularity of the exact solution and remains compatible with the general Leray--Lions structure and with independent bulk and face degeneracy parameters.

Numerical experiments in dimensions two, three, and four illustrate these results. For singular solutions, convergence is observed in every configuration tested, while the attainable rates are limited by the regularity of the solution. The error also decreases under refinement for exponent fields approaching the lower admissibility bound, with $p(\cdot)$ ranging over $[1.05,2.95]$, and for large exponents ranging over $[2,6]$, both for the degenerate flux and for its regularization.

The paper is organized as follows. Section~\ref{sec:continuous.setting} presents the continuous problem, the assumptions on the flux, and the weak formulation. Section~\ref{sec:scheme} introduces the mesh, the HHO unknowns, the reconstruction operators, the stabilization, and the discrete problem. Section~\ref{sec:discrete.norms} develops the variable-exponent discrete inequalities, and Section~\ref{sec:properties} establishes the properties of the discrete operators. Section~\ref{sec:well-posedness} addresses discrete well-posedness, and Section~\ref{sec:convergence} proves compactness and convergence. Section~\ref{sec:numerics-convergence} presents numerical experiments in dimensions two, three, and four. The appendices collect the variable-exponent functional framework and polynomial tools.

\section{Continuous problem}\label{sec:continuous.setting}

Let $\Omega \subset \RR^d$, $d\in\NN^*$, denote a bounded, connected, polytopal open set with Lipschitz boundary $\partial\Omega$. We consider problem~\eqref{eq:intro.problem}. The function $f : \Omega \to \RR$ is the volumetric source term, while $\stress : \Omega \times \RR^d \to \RR^d$ is the flux function; it may vary with the spatial position and depends on the potential only through its gradient.

{
\br[Nonhomogeneous Dirichlet data and manufactured solutions]\label{rem:nonhomogeneous}
We treat homogeneous Dirichlet data in the analysis and distinguish this setting from the manufactured tests with prescribed traces. Some manufactured solutions used in the numerical section have a nonzero boundary trace $g$. Provided $g$ admits a lifting $G\in W^{1,p(\cdot)}(\Omega)$ with ${\rm tr}\,G=g$, the associated continuous weak problem is posed in the affine space $G+W^{1,p(\cdot)}_0(\Omega)$ and tested against $W^{1,p(\cdot)}_0(\Omega)$. One can equivalently write $u=G+v$ and obtain a nonlinear problem for $v$ with zero trace; the flux then depends on $\nabla G+\nabla v$. At the discrete level, boundary face unknowns must be prescribed using an appropriate approximation of $g$. The same structural monotonicity arguments are expected to apply to differences of admissible solutions, but a full extension of the convergence and, in particular, of the quantitative error estimates requires control of the chosen discrete lifting and its additional consistency terms. Such an extension is not claimed in the theorems below.
\er
}

\subsection{Flux function}\label{sec:strain.rate.shear.stress.law}

The following assumptions characterize $\stress$, where we refer to Appendix~\ref{app:ve-background} for more details on variable-exponent Lebesgue and Sobolev spaces.

\bas[Flux function]\label{ass:stress}
  The flux function $\stress : \Omega \times \RR^d \to \RR^d$ is measurable, and there exist $\sob \in \Plogpm(\Omega)$, $\delta \in L^{p(\cdot)}(\Om,[0,\infty))$ and $\sigma_\cst{c},\sigma_\cst{m} \in (0,\infty)$ such that, for all $\b\tau,\b\eta \in \RR^d$ with $\b \tau \ne \b \eta$, and a.e. $\b x \in \Omega$, the following H\"older continuity, strong monotonicity, and vanishing-at-rest properties hold, respectively:
 \bse\label{eq:ass:stress}
 \ba
       |\stress(\b x,\b\tau)-\stress(\b x,\b\eta)| &\le \sigma_\cst{c} \left(\delta(\b x)^{\sob(\b x)}+|\b\tau|^{\sob(\b x)}+|\b\eta|^{\sob(\b x)}\right)^\frac{{\sob(\b x)}-2}{{\sob(\b x)}}| \b\tau-\b\eta |,
    \label{eq:ass:stress:cont} 
\\
      \left(\stress(\b x,\b\tau)-\stress(\b x,\b\eta)\right)\cdot\left(\b\tau-\b\eta\right) &\ge \sigma_\cst{m}\left(\delta(\b x)^{\sob(\b x)}+|\b\tau|^{\sob(\b x)}+|\b\eta|^{\sob(\b x)}\right)^\frac{{\sob(\b x)}-2}{{\sob(\b x)}}|\b\tau-\b\eta|^{2},\label{eq:ass:stress:mono}\\
      \stress(\b x,\b 0) &= \b 0.\label{eq:ass:stress:zero}
\ea
\ese
\eas
{In the above inequalities and throughout the paper, a non-negative quantity raised to the power zero is interpreted as $1$, including at a zero base; expressions involving coincident arguments in the structural estimates are understood by the corresponding limiting value.}

{
In the remainder of this section, to avoid a proliferation of generic constants, we write $a\lesssim b$ (respectively, $a\gtrsim b$) for non-negative quantities $a,b$ whenever $a\le Cb$ (respectively, $a\ge Cb$), where $C\in(0,\infty)$ depends only on $p_\petm$, $p_\petp$, $p_\cst{log}$, and $\Omega$ (including the fixed spatial dimension $d$). We also write $a\eqsim b$ if $C^{-1}b\le a\le Cb$, equivalently $a\lesssim b\lesssim a$. In particular, these hidden constants do not depend on $\sigma_\cst{m}$ or $\sigma_\cst{c}$, which are kept explicit when relevant. Additional dependencies are specified when needed.
}

Some remarks are in order.

\br[Residual shear stress]
  Assumption \eqref{eq:ass:stress:zero} can be relaxed by taking $\stress(\cdot,\b 0)  \in L^{{\sob}'(\cdot)}(\Omega,\RR^d)$. This modification requires only minor changes in the analysis, not detailed for the sake of conciseness.
\er

\br[Singular exponents]
  Inequalities \eqref{eq:ass:stress:cont}--\eqref{eq:ass:stress:mono} yield the following bounds, which are convenient for the H\"older arguments used later.
  For all $\b\tau,\b\eta \in \RR^d$ with $\b\tau \neq \b\eta$ and almost every $\b x \in \Omega$,
  \bse\label{eq:ass:stress:sing}
  \ba
      \left|
        \stress(\b x,\b\tau)-\stress(\b x,\b\eta)
        \right| &\lesssim\sigma_\cst{c} \left(\delta(\b x)^{\sob(\b x)}+|\b\tau|^{\sob(\b x)}+|\b\eta|^{\sob(\b x)}\right)^{\frac{{\sob(\b x)}-{\sobmin(\b x)}}{{\sob(\b x)}}}| \b\tau-\b\eta |^{{\sobmin(\b x)}-1},\label{eq:ass:stress:sing:cont}
      \\
        \left(\stress(\b x,\b\tau)-\stress(\b x,\b\eta)\right)\cdot(\b\tau-\b\eta)  &\gtrsim\sigma_\cst{m}\left(\delta(\b x)^{\sob(\b x)}+|\b\tau|^{\sob(\b x)}+|\b\eta|^{\sob(\b x)}\right)^{\frac{{\sob(\b x)-\sobmax(\b x)}}{{\sob(\b x)}}}|\b\tau-\b\eta|^{{\sobmax(\b x)}}.\label{eq:ass:stress:sing:mono}
\ea
\ese
where $\sobmin$ and $\sobmax$ are defined in \eqref{eq:exp:conj.lower.uper}.
{These bounds follow from \eqref{eq:ass:stress:cont}--\eqref{eq:ass:stress:mono} and the triangle inequality $|\b\tau-\b\eta|^p\le 2^{p-1}(|\b\tau|^p+|\b\eta|^p)$, applied separately for $p<2$ and $p\ge2$. The hidden constants depend only on $p_\petm$ and $p_\petp$.}

\er

\br[Order between the continuity and monotonicity constants]
  One can see that $\sigma_\cst{m} \le \sigma_\cst{c}$ by taking $\b\tau,\b\eta \in \RR^d$ such that $|\b\tau| = 1$ and $\b \eta = \b 0$ in \eqref{eq:ass:stress} and writing a.e. in $\Om$,
  \bae
 \label{eq:power-framed:constants.bound}
    \sigma_\cst{m} \le  \frac{\stress(\cdot,\b\tau)\cdot\b\tau}{\left(\delta^{\sob}+|\b\tau|^{\sob}\right)^\frac{{\sob}-2}{{\sob}}|\b\tau|^{2}}\le \frac{|\stress(\cdot,\b\tau)||\b\tau|}{\left(\delta^{\sob}+|\b\tau|^{\sob}\right)^\frac{{\sob}-2}{{\sob}}|\b\tau|^{2}} \le \sigma_\cst{c}.
\eae
\er

We next give examples of flux functions $\stress$ satisfying Assumption~\ref{ass:stress}.

\begin{example}[$p(\cdot)$-Laplace flux function]\label{ex:laplace}
The $p(\cdot)$-Laplace equation $-\DIV(|\GRAD u|^{p(\cdot)-2}\GRAD u)=f$ corresponds to $\stress(\b x,\b\tau)=|\b\tau|^{p(\b x)-2}\b\tau$ and $\delta(\b x)=0$. The standard vectorial $p$-power inequalities, applied pointwise with $p=p(\b x)$ and using $1<p_-\le p(\b x)\le p_+<\infty$, show that Assumption~\ref{ass:stress} holds with positive constants depending only on $p_-$ and $p_+$; see, e.g., \cite[Appendix~A]{Botti.Castanon-Quiroz.ea:20} for the corresponding power-framed inequalities.
\end{example}

\begin{example}[Carreau--Yasuda flux function]\label{ex:Carreau--Yasuda}
  Another example of function $\stress$ which satisfies Assumption \ref{ass:stress}, inspired by the rheology of Carreau--Yasuda fluids, is obtained by setting, for almost every $\b x\in\Omega$ and all $\b\tau \in \RR^d$,
  \bae\label{eq:Carreau--Yasuda}
    \stress(\b x,\b\tau) = \mu(\b x)\left(\delta(\b x)^{\alpha(\b x)}+|\b\tau|^{\alpha(\b x)}\right)^\frac{\sob(\b x)-2}{\alpha(\b x)}\b\tau,
  \eae
  where $\mu : \Omega \to [\mu_\petm,\mu_\petp] \subset (0,+\infty)$ is a measurable function corresponding to the \emph{local flow consistency index}, $\delta \in L^{p(\cdot)}(\Omega,[0,+\infty))$ is the \emph{degeneracy parameter}, $\alpha : \Omega \to [\alpha_\petm,\alpha_\petp] \subset (0,+\infty)$ is a measurable function expressing the \emph{local transition flow behavior index}, and $\sob \in \Plogpm(\Omega)$ is the \emph{flow behavior index}.
  The power-framed estimates in \cite[Appendix~A]{Botti.Castanon-Quiroz.ea:20}, applied pointwise in $\b x$ and uniformly over the prescribed ranges of $p(\b x)$, $\alpha(\b x)$, and $\mu(\b x)$, show that $\stress$ matches Assumption~\ref{ass:stress} with positive continuity and monotonicity constants depending only on these ranges.
Finally, when $\delta \equiv 0$ and $\mu \equiv 1$, we get the $p(\cdot)$-Laplace example \ref{ex:laplace}.
\end{example}

\subsection{Weak formulation}\label{sec:weak.formulation}

We define the following space for the potential embedding the homogeneous boundary condition: 
\bae
U \eq W_0^{1,\sob(\cdot)}(\Omega) \eq \{v \in W^{1,\sob(\cdot)}(\Omega) : v\restri{\partial\Omega} = 0\}.
\eae
Assuming $f \in L^{\sob'(\cdot)}(\Omega)$, the weak formulation of problem \eqref{eq:intro.problem} reads:
Find $u \in U$ such that
\bae\label{eq:lr.weak}
     \cst{a}(u,v) = (f, v)_\Om \qquad \forall v \in U,
\eae
where the \emph{diffusion function} $\cst{a} : U \times U \to \RR$ is defined such that, for all $v,w \in U$,
\bae\label{eq:a.b}
  \cst{a}(w,v) \eq ( \stress(\cdot,\GRAD w) ,\GRAD v)_\Om.
\eae

We recall the following inequality between sums of powers (see \cite[Eq. (15)]{Botti.Castanon-Quiroz.ea:20}):
let an integer $N\ge 1$ and a real number $m \in (0,+\infty)$ be given. Then, for all $x_1,\ldots,x_N \in [0,+\infty) $, we have
\bae\label{eq:sum-power}
  N^{\min(0,m-1)}\sum_{i=1}^N x_i^m \le \left(\sum_{i=1}^N x_i\right)^m  \le N^{\max(0,m-1)}\sum_{i=1}^N x_i^m.
\eae

\bl[Strong monotonicity of $\cst{a}$]\label{lem:a:holder.continuity.strong.monotonicity}
Under Assumption \ref{ass:stress}, for all $v,w\in U$, setting $e\eq v-w$, it holds
\bae\label{eq:a:strong.monotonicity}
\cst{a}(v,e)-\cst{a}(w,e)
\gtrsim\sigma_\cst{m}
\Big(\rho_{L^{\sob(\cdot)}(\Om)}(\delta)
+\rho_{\tL^{\sob(\cdot)}(\Om)}(\GRAD v)
+\rho_{\tL^{\sob(\cdot)}(\Om)}(\GRAD w)\Big)^{1-\frac{\sobmax_\petm}{\sob_\petm}}
\rho_{\tL^{\sob(\cdot)}(\Om)}(\GRAD e)^{\frac{\sobmax_\petm}{\sob_\petm}}.
\eae
{The hidden constant depends only on $p_\petm$ and $p_\petp$.}
\el

\br[Modular and norm expressions]
It is possible to have a version of \eqref{eq:a:strong.monotonicity} written with norms instead of modulars. However, to avoid the proliferation of terms to the powers $1/p_-$ and $1/p_+$ appearing by using \eqref{eq:modular-norm-equiv}, we keep the modular form.
\er

\bp

If $\rho_{\tL^{\sob(\cdot)}(\Om)}(\GRAD e)=0$, the result is immediate. Otherwise, using \eqref{eq:ass:stress:sing:mono} and the variable-exponent Young inequality \eqref{eq:young:gener}, we obtain for $\varepsilon\in(0,1]$
{
\bae\label{eq:a:sm:1}
\rho_{\tL^{\sob(\cdot)}(\Om)}(\GRAD e)
&\lesssim \varepsilon\Big(\rho_{L^{\sob(\cdot)}(\Om)}(\delta)+\rho_{\tL^{\sob(\cdot)}(\Om)}(\GRAD v)+\rho_{\tL^{\sob(\cdot)}(\Om)}(\GRAD w)\Big)\\
&\quad+\varepsilon^{1-\frac{\sobmax_\petm}{\sob_\petm}}\sigma_\cst{m}^{-1}
\big(\stress(\cdot,\GRAD v)-\stress(\cdot,\GRAD w),\GRAD e\big)_\Om.
\eae
}
Here we used the fact that $t\mapsto 1-\max(2,t)/t$ is nondecreasing on $(1,\infty)$ and, for $\varepsilon\in(0,1]$, 
{
\bae\label{eq:eps:pm}
\sup_{\b x\in\Omega}\varepsilon^{1-\sobmax(\b x)/\sob(\b x)}
=\varepsilon^{1-\sobmax_\petm/\sob_\petm}.
\eae
}
The triangle inequality and \eqref{eq:sum-power} also give
{
\bae
\rho_{\tL^{\sob(\cdot)}(\Om)}(\GRAD e)
\lesssim \rho_{L^{\sob(\cdot)}(\Om)}(\delta)+\rho_{\tL^{\sob(\cdot)}(\Om)}(\GRAD v)+\rho_{\tL^{\sob(\cdot)}(\Om)}(\GRAD w).
\eae
}
Consequently, we can choose
{
\bae\label{eq:eps:choice}
\varepsilon=\theta\frac{\rho_{\tL^{\sob(\cdot)}(\Om)}(\GRAD e)}{\rho_{L^{\sob(\cdot)}(\Om)}(\delta)+\rho_{\tL^{\sob(\cdot)}(\Om)}(\GRAD v)+\rho_{\tL^{\sob(\cdot)}(\Om)}(\GRAD w)},
\eae
}
where $\theta>0$ depends only on $p_\petm,p_\petp$ and is sufficiently small both to ensure $\varepsilon\le1$ and to absorb the first right-hand term in \eqref{eq:a:sm:1}. Rearranging then yields \eqref{eq:a:strong.monotonicity}, with a hidden constant depending only on $p_\petm,p_\petp$.
\ep

{We will use the variable-exponent version of the Poincaré inequality, cf. \cite[Theorem 8.2.4]{diening2011lebesgue}: for all $v\in W_0^{1,p(\cdot)}(\Omega)$,}
\bae\label{eq:Poincare}
 \|v\|_{L^{\sob(\cdot)}(\Om)}\lesssim\|\GRAD v\|_{\tL^{\sob(\cdot)}(\Om)}.
\eae
{The hidden constant in this Poincar\'e inequality depends only on $p_\petm$, $p_\petp$, $p_\cst{log}$, and $\Omega$.}

The well-posedness of problem \eqref{eq:lr.weak} is shown in \cite{diening2011lebesgue,Fan2011}. Nevertheless, we give a proof since the same argument will be adapted at the discrete level.

\bt[{ Continuous well-posedness and variable-exponent a priori estimate}]\label{prop:a-priori}
Under Assumption \ref{ass:stress}, the continuous problem \eqref{eq:lr.weak} admits a unique solution $u\in U$ satisfying
\bae\label{eq:continuous.solution:bounds:uh}
\rho_{\tL^{\sob(\cdot)}(\Om)}(\GRAD u)
\lesssim\max(1,\sigma_\cst{m}^{-p'_\petm},\sigma_\cst{m}^{-p'_\petp})\max\!\left(
\rho_{L^{\sob'(\cdot)}(\Om)}(f)^{p'_\petp/p'_\petm},
\rho_{L^{\sob'(\cdot)}(\Om)}(f)^{p'_\petm/p'_\petp},
\rho_{L^{\sob(\cdot)}(\Om)}(\delta)
\right).
\eae
{The hidden constant depends only on $p_\petm$, $p_\petp$, $p_\cst{log}$, and $\Omega$.}
\et

	\bp

(i) \emph{Existence.}
Equip $U$ with the norm $\|v\|_U\eq\|\GRAD v\|_{\tL^{\sob(\cdot)}(\Om)}$, and let $A:U\to U^*$ be defined by $\langle A(w),v\rangle_{U^*,U}\eq\cst{a}(w,v)$. By \eqref{eq:ass:stress:cont} and \eqref{eq:ass:stress:zero}, $A$ is bounded and hemicontinuous, while \eqref{eq:a:strong.monotonicity} yields strict monotonicity. If $\rho_{\tL^{\sob(\cdot)}(\Om)}(\GRAD v)\ge\max(1,\rho_{L^{\sob(\cdot)}(\Om)}(\delta))$, then
{
\bae\label{eq:phi:coercive}
\langle A(v),v\rangle_{U^*,U}=\cst{a}(v,v)
&\refover{\eqref{eq:a:strong.monotonicity}}\gtrsim\sigma_\cst{m}\Big(\rho_{L^{\sob(\cdot)}(\Om)}(\delta)+\rho_{\tL^{\sob(\cdot)}(\Om)}(\GRAD v)\Big)^{1-\frac{\sobmax_\petm}{\sob_\petm}}\rho_{\tL^{\sob(\cdot)}(\Om)}(\GRAD v)^{\frac{\sobmax_\petm}{\sob_\petm}}\\
&\gtrsim\sigma_\cst{m}\rho_{\tL^{\sob(\cdot)}(\Om)}(\GRAD v)
\refovern{\eqref{eq:modular-norm-equiv}}\ge\sigma_\cst{m}\|v\|_U^{p_\petm}.
\eae
}
Hence $A$ is coercive. Since $U$ is reflexive, the Browder--Minty theorem applied to $A$ and the functional $v\mapsto(f,v)_\Om\in U^*$ yields a solution of \eqref{eq:lr.weak}.

(ii) \emph{A priori estimate.}
{Suppose that $\rho_{\tL^{\sob(\cdot)}(\Om)}(\GRAD u)\ge\rho_{L^{\sob(\cdot)}(\Om)}(\delta)$ and set}
{
\bae
 \#\eq\operatorname{sign}\!\left(
 1-\|\GRAD u\|_{\tL^{\sob(\cdot)}(\Om)}\right),
 \qquad \operatorname{sign}(0)=+1.
\eae
}
{By} the strong-monotonicity estimate~\eqref{eq:a:strong.monotonicity} with $v=u$ and $w=0$, using $\rho_{\tL^{p(\cdot)}(\Omega)}(\GRAD u)\ge\rho_{L^{p(\cdot)}(\Omega)}(\delta)$, \eqref{eq:lr.weak}, the variable-exponent H\"older inequality \eqref{eq:holder}, and Poincar\'e's inequality \eqref{eq:Poincare},
{
\bae\label{eq:proof:apriori}
\sigma_\cst{m}\rho_{\tL^{\sob(\cdot)}(\Om)}(\GRAD u)
&\lesssim \cst{a}(u,u)=(f,u)_\Om\\
&\lesssim\|f\|_{L^{\sob'(\cdot)}(\Om)}\|\GRAD u\|_{\tL^{\sob(\cdot)}(\Om)}
\refovern{\eqref{eq:modular-norm-equiv}}\le
\|f\|_{L^{\sob'(\cdot)}(\Om)}\rho_{\tL^{\sob(\cdot)}(\Om)}(\GRAD u)^{1/p_\#}.
\eae
}
If the gradient modular is nonzero, division by its $1/p_\#$ power yields
{
\bae\label{eq:proof:apriori:power}
\rho_{\tL^{\sob(\cdot)}(\Om)}(\GRAD u)
\lesssim\sigma_\cst{m}^{-p'_\#}\|f\|_{L^{\sob'(\cdot)}(\Om)}^{p'_\#}.
\eae
}
Using \eqref{eq:modular-norm-equiv} and the two possible values of $\#$ gives \eqref{eq:continuous.solution:bounds:uh}. If the gradient modular is zero, or is smaller than $\rho_{L^{\sob(\cdot)}(\Om)}(\delta)$, the claimed bound follows immediately.

(iii) \emph{Uniqueness.}
For two solutions $u,u'\in U$, put $e=u-u'$. By \eqref{eq:lr.weak} and \eqref{eq:a:strong.monotonicity},
{
\bae
0=\cst{a}(u,e)-\cst{a}(u',e)
\gtrsim\sigma_\cst{m}\Big(\rho_{L^{\sob(\cdot)}(\Om)}(\delta)+\rho_{\tL^{\sob(\cdot)}(\Om)}(\GRAD u)+\rho_{\tL^{\sob(\cdot)}(\Om)}(\GRAD u')\Big)^{1-\frac{\sobmax_\petm}{\sob_\petm}}\rho_{\tL^{\sob(\cdot)}(\Om)}(\GRAD e)^{\frac{\sobmax_\petm}{\sob_\petm}}.
\eae
}
Therefore $\rho_{\tL^{\sob(\cdot)}(\Om)}(\GRAD e)=0$, so $u=u'$.
\ep

\section{HHO discretization}\label{sec:scheme}

This section presents the discrete spaces, the local operators, and the HHO discretization of problem~\eqref{eq:lr.weak}.

\subsection{Discrete setting}\label{sec:discrete.setting}

\paragraph{Mesh.}

For any set $U\subset\RR^d$, denote by $h_U$ its diameter.
A polytopal mesh is defined as a couple $\mathcal M_h\coloneq(\T_h,\F_h)$, where $\T_h$ is a finite collection of polytopal elements, whose largest diameter $h\coloneq\max_{T\in\T_h}h_T$ is the mesh size, while $\F_h$ is a finite collection of hyperplanar faces.
It is assumed henceforth that the mesh $\mathcal M_h$ matches the geometrical requirements detailed in \cite[Definition 1.7]{Di-Pietro.Droniou:20}, and we denote by $\kappa$ the mesh regularity parameter.
Boundary faces lying on $\partial\Omega$ and internal faces contained in
$\Omega$ are collected in the sets $\F_h^{\rm b}$ and $\F_h^{\rm i}$, respectively.
For every mesh element $T\in\T_h$, we denote by $\F_T$ the subset of $\F_h$ collecting the faces that lie on the boundary $\partial T$ of $T$. For every face $F \in \F_h$, we denote by $\T_F$ the subset of $\T_h$ containing the one (if $F\in\F_h^{\rm b}$) or two (if $F\in\F_h^{\rm i}$) elements on whose boundary $F$ lies.
For each mesh element $T\in\T_h$ and face $F\in\F_T$, $\b n_{TF}$ denotes the (constant) unit vector normal to $F$ pointing out of $T$, and we denote $(\b n_{\partial T})\res{F} \eq \b n_{TF}$.

Our focus is on the $h$-convergence analysis, so we consider a sequence of refined meshes that is regular in the sense of \cite[Definition 1.9]{Di-Pietro.Droniou:20}, with regularity parameter uniformly bounded away from zero. We write $\mathfrak H\subset(0,\infty)$ for the set of mesh sizes of this sequence, whose only accumulation point is $0$.
Consequently, cell diameters and the characteristic face scales are uniformly comparable, and the number of faces per cell remains bounded independently of mesh refinement.

\br[One-dimensional meshes]\label{rm:one-dimensional-mesh}
When $d=1$, every cell $T$ is an interval and each face $F$ is an endpoint. Face integrals are understood with respect to the counting measure $\mathcal H^0$; in particular, $\PP^k(F)=\RR$ and an integral over $\partial T$ is a sum over its two endpoints. Since the geometric diameter of $F$ vanishes, any face length $h_F$ used in the analysis is interpreted as the characteristic scale
\bae
 h_F\coloneq\frac{1}{|\T_F|}\sum_{T\in\T_F}h_T.
\eae
For regular one-dimensional meshes, $h_F\eqsim h_T$ whenever $F\in\F_T$; cf.~\cite[Remark~1.13]{Di-Pietro.Droniou:20}. The HHO unknowns, local reconstructions, and stabilization use cell diameters $h_T$ and therefore require no modification in this dimension.
\er

\paragraph{Projectors and broken spaces.}

Given $U \in \T_h \cup \F_h$ and $\ell \in \NN$, we denote by $\PP^\ell(U)$ the space spanned by the restriction to $U$ of scalar-valued, $d$-variate polynomials of total degree $\le \ell$.
The local $L^2$-orthogonal projector $\proj{U}{\ell} : L^{1}(U) \to \PP^\ell(U)$ is defined such that, for all $v \in L^{1}(U)$,
\bae\label{eq:proj}
  (\proj{U}{\ell}v, w)_U = (v,w)_U \qquad \forall w \in  \PP^{\ell}(U).
\eae
When applied to vector-valued functions in $L^1(U)^d$, the $L^2$-orthogonal projector mapping on $\PP^\ell(U)^d$ acts component-wise and is denoted in boldface font as $\PROJ{U}{\ell}$. We will also use the notations $\proj{h}{\ell}$ and $\proj{\partial T}{\ell}$ such that for all $T \in \T_h$ and all $F \in \F_T$, $(\proj{h}{\ell})\res{T} = \proj{T}{\ell}$ and $(\proj{\partial T}{\ell})\res{F} = \proj{F}{\ell}$.
The elliptic projector $\proj{T}{1,k+1}:W^{1,1}(T)\to\PP^{k+1}(T)$ is characterized by $(\GRAD(\proj{T}{1,k+1}v-v),\GRAD w)_T=0$ for all $w\in\PP^{k+1}(T)$, together with $(\proj{T}{1,k+1}v-v,1)_T=0$.

For the approximation properties used below, we employ the broken Sobolev spaces
\bae
W^{n,\sob(\cdot)}(\T_h)\coloneq\left\{ v\in L^{\sob(\cdot)}(\Omega)\ : \ v\res{T}\in W^{n,\sob(\cdot)}(T)\quad\forall T\in\T_h\right\}.
\eae
The corresponding seminorm is defined by  \eqref{eq:luxembourg} with the modular defined such that 
\ba 
\rho_{\dot W^{n,\sob(\cdot)}(\T_h)}(v) &\eq \sum_{T\in\T_h}\rho_{\dot W^{n,\sob(\cdot)}(T)}(v\res{T}).
\ea

For a piecewise Sobolev function $v$, the broken gradient $\nabla_h v$ is defined cellwise by $(\nabla_h v)|_T=\GRAD(v|_T)$ for each $T\in\T_h$.

\subsection{Hybrid space and norms}

Let an integer $k\ge 0$ be fixed and let $\widetilde k \in \{k,k+1\}$. The HHO unknowns are polynomials of degree $\widetilde k$ on the mesh elements and of degree $k$ on the mesh faces. 
The choice of $\widetilde k$ gives the equal-order formulation when $\widetilde k = k$, and the mixed-order formulation when $\widetilde k = k+1$. This notation allows us to merge both cases in the construction and in the following proofs.  
We define the HHO space by
\bae
\dU{h}{k} \eq \left\{
\u v_h = ((v_T)_{T \in \T_h},(v_F)_{F\in \F_h}) \ : \ v_T \in \PP^{\widetilde k}(T)\ \ \forall T \in \T_h\ \mbox{ and }\ v_F \in \PP^k(F)\ \ \forall F \in \F_h \right\}.
\eae
The discrete potential is sought in the subspace of $\dU{h}{k}$ embedding the homogeneous boundary condition:
\bae
\dU{h,0}{k} \eq \left\{ \u v_h = ((v_T)_{T \in \T_h},(v_F)_{F\in \F_h}) \in \dU{h}{k} \ : \ v_F = 0 \quad \forall F \in \Fb \right\}.
\eae
The interpolation operator $\Interp{h}{k} : W^{1,1}(\Omega) \to  \dU{h}{k} $ maps a function $v \in W^{1,1}(\Omega)$ to the vector of discrete unknowns $\Interp{h}{k}v$ defined as follows:
\bae\label{eq:I:def}
  \Interp{h}{k} v \eq ((\proj{T}{\widetilde k} v\res{T})_{T \in \T_h},(\proj{F}{k} v\res{F})_{F \in \F_h}).
\eae
It is used to compare continuous functions with discrete unknowns.
For all $T \in \T_h$, we denote by $\dU{T}{k}$ and $\Interp{T}{k}$ the restrictions of $\dU{h}{k}$ and $\Interp{h}{k}$ to $T$, respectively, and, for all $\u v_h \in \dU{h}{k}$, we let $\u v_T \eq (v_T,(v_F)_{F\in \F_T}) \in \dU{T}{k}$ denote the vector collecting the discrete unknowns attached to $T$ and its faces.
Furthermore, for all $\u v_h \in \dU{h}{k}$, we define the broken polynomial field $v_h\in\PP^{\widetilde k}(\T_h)$ obtained by patching element unknowns, that is,
\bae
(v_h)\res{T} \eq v_T\qquad\forall T \in \T_h.
\eae
Finally, for all $\u v_T \in \dU{T}{k}$, we define the broken polynomial field $v_{\partial T} \in \PP^{k}(\F_T)$ obtained by patching face unknowns, that is, 
\bae
(v_{\partial T})\res{F} \eq v_F\qquad\forall F \in \F_T.
\eae

\subsection{Operator reconstructions}

The HHO discretization relies on two local ingredients: (i) reconstructions of a gradient and a potential from the cell and face unknowns; (ii) a stabilization penalizing the difference between the trace of the cell unknown and the face unknowns. The discrete problem is then assembled elementwise.

\paragraph{Local gradient reconstruction.}

For all $T \in \T_h$, we define 
$\dgrad{k}{T} : \dU{T}{k} \to \PP^{k}(T)^d$ such that, for all $\u v_T \in \dU{T}{k}$,
\bae\label{eq:G}
  (\dgrad{k}{T} \u v_T,\b\tau)_T = (\GRAD v_T,\b\tau)_T + (v_{\partial T}-v_T,\b\tau \cdot \b n_{\partial T})_{\partial T}\qquad \forall \b\tau \in  \PP^{k}(T)^d.
\eae
Relation \eqref{eq:G} mimics an integration by parts, the face unknowns playing the role of the trace on $\partial T$. As in~\cite{Di-Pietro.Droniou:17}, the gradient is sought in the full space $\PP^k(T)^d$, which yields the commutation property \eqref{eq:G:proj} below.

\paragraph{Local potential reconstruction.}

We define 
$\drec{k+1}{T} : \dU{T}{k} \to \PP^{k+1}(T)$ such that, for all $\u v_T \in \dU{T}{k}$,
\bse \label{eq:rT}
\ba
(\GRAD \drec{k+1}{T} \u v_T,\GRAD w)_T &= (\GRAD v_T,\GRAD w)_T + (v_{\partial T}-v_T,\GRAD w \cdot \b n_{\partial T})_{\partial T} \qquad \forall w \in  \PP^{k+1}(T),\\
  (\drec{k+1}{T} \u v_T-v_T,1)_T &= 0.
\ea
\ese
Recalling the elliptic projector (see \cite[(2.14)]{Di-Pietro.Droniou:20}), we have for all $v \in W^{1,1}(T)$,
\bae
\drec{k+1}{T}\Interp{T}{k}v = \proj{T}{1,k+1} v.
\eae
The potential reconstruction is only used in the stabilization, to which it provides higher-order consistency.

\paragraph{Boundary residual operator.}

We define 
$\dfbres{k}{\partial T} : \dU{T}{k} \rightarrow \PP^{k}(\F_T)$ such that, for all $\u v_T \in \dU{T}{k}$,
   \bae\label{eq:dfbres}
     \dfbres{k}{\partial T} \u v_T \eq
     \frac{1}{h_T}\proj{\partial T}{k}\left[
       v_T-v_{\partial T}+(I-\proj{T}{\widetilde k})(\drec{k+1}{T} \u v_T)
       \right].
   \eae
For every $F\in\F_T$, $\dfbres{k}{F}\u v_T$ denotes the restriction to $F$ of $\dfbres{k}{\partial T}\u v_T$.
\br[Unified treatment of the equal- and mixed-order formulations]
 In the mixed-order formulation, i.e., $\widetilde k = k+1$, for all $T \in \T_h$ and all $\u v_T \in \dU{T}{k}$, observing that $\drec{k+1}{T} \u v_T \in \PP^{k+1}(T)$, the term $(I-\proj{T}{\widetilde k})(\drec{k+1}{T} \u v_T)$ in \eqref{eq:dfbres} vanishes and we get the usual mixed-order boundary residual operator. 
\er

The operator $\dfbres{k}{\partial T}$ measures the mismatch between the trace of the cell unknown and the face unknowns, scaled as a gradient, and vanishes on interpolates of polynomials of degree $k+1$ (see \eqref{eq:sT:polynomial.consistency}).

We define the global versions $\dgrad{k}{h}$ and $\dfbres{k}{\partial h}$ of the local operators $\dgrad{k}{T}$ and $\dfbres{k}{\partial T}$. 

\subsection{Stabilization function}

\paragraph{Global skeleton scalar product.}
We use the broken skeleton $\partial\T_h^{\mathrm{br}}\coloneq\bigsqcup_{T\in\T_h}\partial T$, where each interior face is represented separately for its two adjacent elements. We define for all $q \in \Plogpm(\Om)$ and all $(w,v) \in L^{q'(\cdot)}(\partial\T_h^{\mathrm{br}}) \times L^{q(\cdot)}(\partial\T_h^{\mathrm{br}})$,
\bae
(w,v)_{\partial h} \eq \sum_{T \in \T_h}h_T(w,v)_{\partial T}.
\eae

\paragraph{Stabilization function.}
The stabilization is built on $\dfbres{k}{\partial T}$ through a flux function $\stabfl_F$, which plays for the face residuals the role of $\stress$ for the gradient.
We define the stabilization function $\cst{s}_h : \dU{h}{k} \times \dU{h}{k} \to \RR$ such that, for all $\u v_h,\u w_h \in \dU{h}{k}$,
\bse
\bae\label{eq:sh}
  \cst{s}_h(\u w_h, \u v_h) 
  \eq (\stabfl_{\partial h}(\cdot,\dfbres{k}{\partial h}\u w_h),\dfbres{k}{\partial h}\u v_h)_{\partial h} = \sum_{T\in\T_h}\cst{s}_T(\u w_T,\u v_T),
\eae
with for all $T \in \T_h$, 
\bae\label{eq:sT}
\cst{s}_T(\u w_T,\u v_T) \eq h_T(\stabfl_{\partial T}(\cdot,\dfbres{k}{\partial T}\u w_T),\dfbres{k}{\partial T}\u v_T)_{\partial T},
\eae
\ese
where for all $F \in \F_T$, $(\stabfl_{\partial h})\res{F} = (\stabfl_{\partial T})\res{F} = \stabfl_{F}$, and we consider for all $F \in \F_h$, the stabilization flux function $\stabfl_F : F \times \RR \to \RR$ satisfying the following assumption. 

 \bas[Stabilization flux function]\label{ass:sT}
For all $F \in \F_h$, the function $\stabfl_F : F \times \RR \to \RR$ is measurable and there exists $\delta_F \in L^{p(\cdot)}(F;[0,+\infty))$ such that for all $v,w \in \RR$ with $v \ne w$ and almost every $\b x \in F$,
   \bse\label{eq:S:holder.continuity.strong.monotonicity}
     \ba 
       |\stabfl_F(\b x,w)-\stabfl_F(\b x,v)| &\lesssim \sigma_\cst{c}\left(\delta_F(\b x)^{\sob(\b x)}+|w|^{\sob(\b x)}+|v|^{\sob(\b x)}\right)^\frac{\sob(\b x)-2}{\sob(\b x)}| w-v |,\label{eq:S:holder.continuity}  \\
       \left(\stabfl_F(\b x,w)-\stabfl_F(\b x,v)\right) \left(w-v\right) &\gtrsim \sigma_\cst{m}\left(\delta_F(\b x)^{\sob(\b x)}+|w|^{\sob(\b x)}+|v|^{\sob(\b x)}\right)^\frac{\sob(\b x)-2}{\sob(\b x)}|w-v|^{2},\label{eq:S:strong.monotonicity}\\
       \stabfl_F(\b x,0) &= 0. \label{eq:S:zero}
  \ea
\ese
\eas

 \begin{example}[Stabilization flux function]
   Local stabilization flux functions that match Assumption \ref{ass:sT} can be obtained by setting, for all $F \in \F_h$, all $w \in \RR$, and all $\b x \in F$,
   \bae\label{eq:sF}
     \stabfl_F(\b x,w) \eq \mu_F(\b x)\left(\delta_F(\b x)^{a(\b x)}+|w|^{a(\b x)}\right)^\frac{\sob(\b x)-2}{a(\b x)}w,
   \eae
   with a measurable weight $\mu_F:F\to(0,\infty)$ satisfying $\sigma_\cst{m}\lesssim\mu_F(\b x)\lesssim\sigma_\cst{c}$ for almost every $\b x\in F$. Here $a:F\to[a_{\min},a_{\max}]\subset(0,\infty)$ is a bounded measurable function, uniformly bounded away from zero. The numerical tests use the particular choice $a\equiv2$.
   It can be checked that $\stabfl_F$ is a $p(\cdot)$-power-framed function satisfying \eqref{eq:S:holder.continuity.strong.monotonicity}; see \cite[Appendix A]{Botti.Castanon-Quiroz.ea:20} for a proof.
 \end{example}

 \br[Weight and degeneracy of the stabilization flux function]
 For all $F \in \F_h$, it is possible to adapt the weight of $\stabfl_F$ by choosing, e.g., $\mu_F$ equal to a convex combination $\theta \sigma_\cst{m}+ (1-\theta)\sigma_\cst{c}$ with $\theta \in [0,1]$. The weight may also include a factor depending on the polynomial degree $k$ to improve numerical convergence; in our experiments, we choose $\mu_F=2(k+1)(k+2)$. The implicit comparison constants may then depend on $k$, which is fixed in the analysis.
 Concerning the degeneracy parameter, it is possible to always choose $\delta_F = 1$, even if $\delta = 0$. This flexibility comes from the fact that the stabilization is a purely numerical device: $\stabfl_F$ may differ substantially from the flux function $\stress$.
 \er

 \subsection{Discrete problem}

 The discrete problem reads: Find $\u u_h \in \dU{h,0}{k}$ such that
 \bae
   \label{eq:lr.discrete} \cst{a}_h(\u u_h,\u v_h) = (f,v_h)_\Om \qquad \forall \u v_h \in \dU{h,0}{k},
 \eae
 where the \emph{discrete diffusion function} $\cst{a}_h : \dU{h}{k} \times \dU{h}{k} \to \RR$ is defined such that, for all $\u v_h,\u w_h \in \dU{h}{k}$,
\bse\label{eq:ah}
\bae 
  \cst{a}_h(\u w_h, \u v_h) 
  \eq \sum_{T\in\T_h}\cst{a}_T(\u w_T,\u v_T),
\eae
with for all $T \in \T_h$,
\bae\label{eq:aT}
\cst{a}_T(\u w_T,\u v_T) \eq (\stress(\cdot,\dgrad{k}{T} \u w_T),\dgrad{k}{T} \u v_T)_T+\cst{s}_T(\u w_T,\u v_T).
\eae
\ese

\section{Variable-exponent discrete norms and inequalities}\label{sec:discrete.norms}

We define the variable-exponent discrete seminorms and develop the inverse, trace, and embedding inequalities used in the convergence analysis.

\paragraph{Variable-exponent discrete seminorm.}

For all $q \in \Plogp(\Om)$, we define on $\dU{h}{k}$ the $W^{1,q(\cdot)}(\Omega)$-like modular $\rho_{q(\cdot),h}$ such that, for all $\u v_h \in \dU{h}{k}$,
\bse\label{eq:dseminorm:T:mod}
 \bae
    \rho_{q(\cdot),h}(\u v_h) \eq \sum_{T \in \T_h}\rho_{q(\cdot),T}(\u v_T),
 \eae
  with for all $T \in \T_h$,
  \bae
      \rho_{q(\cdot),T}(\u v_T) \eq \rho_{\tL^{q(\cdot)}(T)}(\GRAD v_T)+ h_T\rho_{L^{q(\cdot)}(\partial T)}(h_T^{-1}(v_{\partial T} - v_T)),
\eae
\ese
and we define via \eqref{eq:luxembourg} the associated Luxemburg seminorms $\| {\cdot} \|_{q(\cdot),h}$ and $\| {\cdot} \|_{q(\cdot),T}$ for all $T\in\T_h$.

\paragraph{Reconstructed-gradient discrete seminorm.}

For all $q\in\Plogp(\Om)$, we also define the HHO reconstructed modular
\bse\label{eq:dseminorm:G:mod}
\bae
 \rho_{q(\cdot),\cst{G},h}(\u v_h)
 \eq \sum_{T\in\T_h}\rho_{q(\cdot),\cst{G},T}(\u v_T),
\eae
with, for all $T\in\T_h$,
\bae
 \rho_{q(\cdot),\cst{G},T}(\u v_T)
 \eq
 \rho_{\tL^{q(\cdot)}(T)}(\dgrad{k}{T}\u v_T)
 +h_T\rho_{L^{q(\cdot)}(\partial T)}
 \big(\dfbres{k}{\partial T}\u v_T\big),
\eae
\ese
and we define via \eqref{eq:luxembourg} the associated Luxemburg seminorms $\|{\cdot}\|_{q(\cdot),\cst{G},h}$ and $\|{\cdot}\|_{q(\cdot),\cst{G},T}$ for all $T\in\T_h$.

\paragraph{Element-mean discrete seminorm.}

For all $q\in\Plogp(\Om)$, we define the element-mean discrete modular
\bse
 \bae
    \rho_{q_h,h}(\u v_h) \eq \sum_{T\in\T_h}\|\u v_T\|_{q_T,T}^{q_T},
 \eae
  with for all $T \in \T_h$, $q_T \eq |T|^{-1}\int_T q(\b x)\,\cst d\b x$, and 
  \bae
      \|\u v_T\|_{q_T,T}^{q_T} &\eq \|\GRAD v_T\|_{\tL^{q_T}(T)}^{q_T}+ h_T\|h_T^{-1}(v_{\partial T} - v_T)\|_{L^{q_T}(\partial T)}^{q_T},
\eae
\ese
and we define via \eqref{eq:luxembourg} the associated Luxemburg seminorms $\| {\cdot} \|_{q_h,h}$ and $\| {\cdot} \|_{q_T,T}$ for all $T \in \T_h$.

\subsection{Equivalence between discrete seminorms}

\bpr[Local and global equivalence of the discrete seminorms]\label{pr:discrete_norms_equiv}
For all $q\in\Plogp(\Om)$, all $T \in \T_h$ and all $\u v_T \in \dU{T}{k}$,
\bse
\ba
     \| \u v_T \|_{q(\cdot),T}
     \eqsim \| \u v_T \|_{q(\cdot),\cst{G},T}
     \eqsim \| \u v_T \|_{q_T,T},
\ea
\ese
and, for all $\u v_h \in \dU{h}{k}$,
\bse\label{eq:disc:mod-norm:equiv}
\ba
    \| \u v_h \|_{q(\cdot),h}
    \eqsim \| \u v_h \|_{q(\cdot),\cst{G},h}
    \eqsim \| \u v_h \|_{q_h,h}.
\ea
\ese
The hidden constants depend only on $\Om$, $k$, $d$, $\kappa$, $q_\cst{log}$, $q_\petm$, and $q_\petp$.
\epr

\bp

Let $T\in\T_h$. If $\|\u v_T\|_{q_T,T}=0$, the local equivalence is immediate. Otherwise, set $\u{\tilde v}_T\eq\|\u v_T\|_{q_T,T}^{-1}\u v_T$. Using the fixed-exponent polynomial Lebesgue embedding \eqref{eq:lebesgue.emb:fixed} as in the proof of Theorem~\ref{thm:mean_exp_equiv}, we get
\bae\label{eq:disc:mean:inverse}
|\GRAD\tilde v(\b x)|
&\lesssim \max(1,h_T^{-d/q_T}) \qquad \forall \b x\in T,\\
|h_T^{-1}(\tilde v_F-\tilde v_T)(\b x)|
&\lesssim \max(1,h_T^{-d/q_T}) \qquad \forall F\in\F_T,\ \forall \b x\in F.
\eae
Consider the disjoint union of $T$ and its faces, endowed with the measure $\cst d\b x$ on $T$ and $h_T\cst ds$ on every $F\in\F_T$, and the piecewise polynomial quantity given by $\GRAD\tilde v$ in $T$ and $h_T^{-1}(\tilde v_F-\tilde v_T)$ on $F$. Its frozen $q_T$-modular is exactly $\|\u{\tilde v}_T\|_{q_T,T}^{q_T}=1$. Moreover, the total measure is equivalent to $h_T^d$, and the inverse estimates \eqref{eq:disc:mean:inverse} and log-H\"older continuity provide the same upper bound as in \eqref{eq:mean.exp.6}. Therefore, applying the $U_\ge$--$U_<$ splitting \eqref{eq:rhoUplus}--\eqref{eq:mean.exp.upper} simultaneously to the cell and face contributions yields
\bae\label{eq:norm:equiv:3}
\rho_{q(\cdot),T}(\u{\tilde v}_T) \eqsim 1.
\eae
By \eqref{eq:mod_norm_equiv}, this proves $\|\u v_T\|_{q(\cdot),T}\eqsim\|\u v_T\|_{q_T,T}$. The global equivalence follows with the same normalization $\u{\tilde v}_h\eq\|\u v_h\|_{q_h,h}^{-1}\u v_h$, applying the splitting \eqref{eq:rhoUplus}--\eqref{eq:mean.exp.upper} on the disjoint union of all cell and weighted face contributions; the constants remain uniform since $\sum_{T\in\T_h}(|T|+h_T|\partial T|)\lesssim|\Om|$.

It remains to compare the reconstructed seminorm with the element-mean one. For $T\in\T_h$, normalize again by $\u{\tilde v}_T\eq\|\u v_T\|_{q_T,T}^{-1}\u v_T$. The fixed-exponent HHO stability (see, e.g., \cite{Di-Pietro.Droniou:17*1,Di-Pietro.Droniou:20}) gives
\bae\label{eq:HHO:fixed:stability}
 \|\dgrad{k}{T}\u{\tilde v}_T\|_{\tL^{q_T}(T)}^{q_T}
 +h_T\|\dfbres{k}{\partial T}\u{\tilde v}_T\|_{L^{q_T}(\partial T)}^{q_T}
 \eqsim1.
\eae
Since $\dgrad{k}{T}\u{\tilde v}_T$ and $\dfbres{k}{\partial T}\u{\tilde v}_T$ are polynomial, the splitting argument \eqref{eq:rhoUplus}--\eqref{eq:mean.exp.upper} applied to \eqref{eq:HHO:fixed:stability}, together with \eqref{eq:mean.exp.diam.equiv}, yields $\rho_{q(\cdot),\cst{G},T}(\u{\tilde v}_T)\eqsim1$. Hence,
$\|\u v_T\|_{q(\cdot),\cst{G},T}\eqsim\|\u v_T\|_{q_T,T}$ by \eqref{eq:mod_norm_equiv}. The global equivalence follows by the same normalization and the splitting \eqref{eq:rhoUplus}--\eqref{eq:mean.exp.upper} over all elements.

\ep

\subsection{Variable-exponent discrete inverse and trace inequalities} 

\bpr[Variable-exponent discrete inverse inequality]\label{pr:dis:inverse:ineq}
Let $q\in\Plogp(\Om)$ and $\ell \in \NN$. We have for all $T\in\T_h$ and all $v \in \PP^\ell(T)$,
\bae\label{eq:dis:inverse:ineq}
\|\GRAD v\|_{\tL^{q(\cdot)}(T)} \lesssim h_T^{-1}\|v\|_{L^{q(\cdot)}(T)}.
\eae
\epr

\bp
Let $T\in \T_h$ and set $q_T \eq |T|^{-1}\int_T q(\b x)\,\cst d\b x$. For $v\in \PP^\ell(T)$, the mean-exponent norm equivalence \eqref{eq:mean.exp.norm} gives
$ \|\GRAD v\|_{\tL^{q(\cdot)}(T)}\eqsim
 \|\GRAD v\|_{\tL^{q_T}(T)}$, and 
 $ \|v\|_{L^{q(\cdot)}(T)}\eqsim\|v\|_{L^{q_T}(T)}$.
Hence, the fixed-exponent inverse inequality \cite[Lemma~1.28]{Di-Pietro.Droniou:20}, applied with exponent $q_T$, yields the result.
\ep

\bpr[Continuous local trace inequality]\label{pr:cont:trace:ineq}
Let $q\in\Plogpm(\Om)$. We have for all $T\in\T_h$ and all  $v \in W^{1,q(\cdot)}(T)$,
\bae\label{eq:cont:trace:ineq}
\|v\|_{L^{q(\cdot)}(\partial T)} \lesssim \|h_T^{-1/q(\cdot)}v\|_{L^{q(\cdot)}(T)}+h_T\|h_T^{-1/q(\cdot)}\GRAD v\|_{\tL^{q(\cdot)}(T)}.
\eae
\epr

\bp
We apply the trace argument on the simplices of the regular matching submesh of $T$ (see Section~\ref{sec:discrete.setting}). Let $K$ be one of these simplices, with $h_K\eqsim h_T$, and denote by $q_K\eq |K|^{-1}\int_K q(\b x)\,\cst{d}\b x$ its mean exponent. Choose a uniformly shape-regular affine map $\b\Phi_K:\widehat K\to K$ from the reference simplex, and set
\bae\label{eq:trace:pullback:def}
\widetilde v\eq v\circ\b\Phi_K,\qquad
\widetilde q\eq q\circ\b\Phi_K,\qquad
|\det \GRAD\b\Phi_K|\eqsim h_K^d,\qquad
J_{\partial\b\Phi_K}\eqsim h_K^{d-1}.
\eae
The log-H\"older condition yields, uniformly in $K$ and $\b x\in\overline K$,
\bae\label{eq:log-holder-cont-h}
h_K^{1/q(\b x)-1/q_K}\eqsim1.
\eae
Changing variables in the volume and boundary modulars and using \eqref{eq:log-holder-cont-h} give the two norm equivalences
\bse\label{eq:trace:scaling}
\ba
\|v\|_{L^{q(\cdot)}(\partial K)}
&\eqsim h_K^{(d-1)/q_K}
\|\widetilde v\|_{L^{\widetilde q(\cdot)}(\partial\widehat K)},\\
\|h_K^{-1/q(\cdot)}v\|_{L^{q(\cdot)}(K)}
&\eqsim h_K^{(d-1)/q_K}
\|\widetilde v\|_{L^{\widetilde q(\cdot)}(\widehat K)}.
\ea
\ese
Moreover, $\widehat\GRAD\widetilde v=(\GRAD\b\Phi_K)^\top(\GRAD v)\circ\b\Phi_K$, and uniform shape regularity yields
\bae\label{eq:trace:scaling:grad}
h_K\|h_K^{-1/q(\cdot)}\GRAD v\|_{\tL^{q(\cdot)}(K)}
\eqsim h_K^{(d-1)/q_K}
\|\widehat\GRAD\widetilde v\|_{\tL^{\widetilde q(\cdot)}(\widehat K)}.
\eae
The variable-exponent trace inequality on $\widehat K$ (see \cite[Corollary~12.1.5]{diening2011lebesgue}) applies with uniform constants: the maps $\b\Phi_K$ have uniformly controlled shape parameters, while the pulled-back exponents $\widetilde q$ have the same uniform endpoint bounds and a uniformly controlled local log-H\"older modulus. It follows from \eqref{eq:trace:scaling}--\eqref{eq:trace:scaling:grad} that
\bae\label{eq:trace:local:simplex}
\|v\|_{L^{q(\cdot)}(\partial K)}
\lesssim\|h_K^{-1/q(\cdot)}v\|_{L^{q(\cdot)}(K)}
+h_K\|h_K^{-1/q(\cdot)}\GRAD v\|_{\tL^{q(\cdot)}(K)}.
\eae
Each face of $T$ is a union, up to sets of surface measure zero, of boundary facets of the submesh. Since the number of simplices is uniformly bounded and their diameters are comparable to $h_T$, summing \eqref{eq:trace:local:simplex} over the boundary simplices and using the triangle inequality yields \eqref{eq:cont:trace:ineq}.
\ep

\bpr[Variable-exponent discrete trace inequality]\label{pr:dis:trace:ineq}
Let $q\in\Plogpm(\Om)$ and $\ell \in \NN$. We have for all $T\in\T_h$ and all $v \in \PP^\ell(T)$,
\bae\label{eq:dis:trace:ineq}
\|v\|_{L^{q(\cdot)}(\partial T)} \lesssim
\|h_T^{-1/q(\cdot)}v\|_{L^{q(\cdot)}(T)}.
\eae
\epr

\bp
Let $v \in \PP^\ell(T)$ and set $q_T \eq |T|^{-1}\int_T q(\b x)\,\cst d\b x$. As above, \eqref{eq:mean.exp.diam.equiv} implies
\bae
\|h_T^{-1/q(\cdot)}v\|_{L^{q(\cdot)}(T)}\eqsim
h_T^{-1/q_T}\|v\|_{L^{q(\cdot)}(T)},
\qquad
\|h_T^{-1/q(\cdot)}\GRAD v\|_{\tL^{q(\cdot)}(T)}\eqsim
h_T^{-1/q_T}\|\GRAD v\|_{\tL^{q(\cdot)}(T)}.
\eae
Applying \eqref{eq:dis:inverse:ineq} to the polynomial $v$ therefore gives
\bae
 h_T\|h_T^{-1/q(\cdot)}\GRAD v\|_{\tL^{q(\cdot)}(T)}
 \lesssim \|h_T^{-1/q(\cdot)}v\|_{L^{q(\cdot)}(T)}.
\eae
Hence, the continuous trace inequality \eqref{eq:cont:trace:ineq} yields the result.
\ep

\subsection{Variable-exponent discrete Poincar\'e--Wirtinger inequality}

\bpr[Discrete Poincar\'e--Wirtinger inequality]
Let $q \in \Plogp(\Om)$ and $\ell \in \NN$. We have for all $T\in\T_h$ and $v  \in \PP^\ell(T)$ such that $\int_T v\,\cst{d}\b x = 0$,
  \bae\label{eq:T:discrete.Sobolev.Poincare.Wirtinger}
    \|  v \|_{L^{q(\cdot)}(T)} &\lesssim h_T\|  \GRAD v \|_{\tL^{q(\cdot)}(T)}, 
\eae
  where the hidden constant depends only on $\ell$, $d$, $\kappa$, $q_\cst{log}$, and $q_\petp$.
  \epr

  \bp 
  Let $T \in \T_h$ and $v  \in \PP^\ell(T)$ such that $\int_T v\,\cst{d}\b x = 0$. Using the mean-exponent equivalence \eqref{eq:mean.exp.norm} and the fixed-exponent Poincar\'e--Wirtinger inequality, we have
 \bae 
  \|  v \|_{L^{q(\cdot)}(T)} = \|  v -\proj{T}{0}v\|_{L^{q(\cdot)}(T)} \lesssim h_T\|  \GRAD v \|_{\tL^{q(\cdot)}(T)}.
  \eae
  \ep 

\subsection{Variable-exponent discrete Sobolev--Poincar\'e inequality}

For all $q \in \Plog(\Om)$, we define the Sobolev exponent function $q^*  : \Om \to [1,\infty]$ such that for all $\b x \in \Om$,
\bae\label{eq:exp:sobolev} 
       q^*(\b x) \eq \begin{cases} \dfrac{dq(\b x)}{d-q(\b x)} & \text{if } q(\b x) < d, \\[2mm] \infty & \text{if } q(\b x) \ge d. \end{cases}
 \eae

The following Sobolev--Poincar\'e inequality is proved in \cite[Theorem 8.3.1]{diening2011lebesgue} under the hypothesis $q_+<d$, and we generalize it to the case where $q$ is possibly larger than $d$.

\bpr[Sobolev--Poincar\'e inequality]\label{pr:sob-poin-wir}
Let {$q\in\Plogpm(\Om)$ and $r \in \Plogp(\Om)$} such that $r(\cdot) \le q^*(\cdot)$ in $\Om$. For all $v  \in W_0^{1,q(\cdot)}(\Om)$,
  \bae\label{eq:Sobolev.Poincare}
    \|  v \|_{L^{r(\cdot)}(\Omega)} \lesssim \|   \GRAD v \|_{\tL^{q(\cdot)}(\Om)},
  \eae
  where the hidden constant depends only on $d$, $\Omega$, $q_\petm$, $q_\petp$, $r_\petm$, $r_\petp$, $q_\cst{log}$, and $r_\cst{log}$.
  \epr

\bp
See \cite[Theorem 8.3.1]{diening2011lebesgue} for a proof when $q_+<d$ (together with the H\"older inequality \eqref{eq:holder} to introduce $r(\cdot)$ since $\Om$ is bounded).
Let $v\in W_0^{1,q(\cdot)}(\Omega)$. If $d=1$, the
fundamental theorem of calculus and the H\"older
inequality \eqref{eq:holder} give for a.e. $x \in \Om$,
\bae
|v(x)| \le \int_\Om |v'(y)|\,\cst{d}y \le {2}\|v'\|_{L^{q(\cdot)}(\Omega)}\|1\|_{L^{q'(\cdot)}(\Omega)} \lesssim \|v'\|_{L^{q(\cdot)}(\Omega)}.
\eae
Hence, since $\Om$ is bounded, $\|v\|_{L^{r(\cdot)}(\Omega)} \lesssim \|v\|_{L^{\infty}(\Omega)} \lesssim \|v'\|_{L^{q(\cdot)}(\Omega)}$.
We now assume $d\ge2$ in the reduction below.
In the general case, set
\bae\label{eq:s0:eps:s}
 s_0(\cdot)\eq\frac{dr(\cdot)}{d+r(\cdot)},\qquad
 \varepsilon_0\eq\frac12\min(q_- -1,d-1),\qquad
 s(\cdot)\eq\max(s_0(\cdot),1+\varepsilon_0).
\eae
Since $r\le q^*$, we have $s_0\le q$ in $\Om$; hence $1<s<d$, $s\le q$, and $r=\dfrac{d s_0}{d-s_0}\le\dfrac{d s}{d-s}=s^*$ in $\Om$, since $0<s_0\le s<d$ and $t\mapsto dt/(d-t)$ is increasing on $(0,d)$. Using \eqref{eq:Sobolev.Poincare} in the case where the exponent on the right-hand side is smaller than $d$, we obtain
\bae\label{eq:proof:SPW}
 \|v\|_{L^{r(\cdot)}(\Omega)}
 \lesssim \|v\|_{L^{s^*(\cdot)}(\Omega)}
 \lesssim \|\GRAD v\|_{\tL^{s(\cdot)}(\Om)}
 \lesssim \|\GRAD v\|_{\tL^{q(\cdot)}(\Om)},
\eae
where we concluded with the H\"older inequality \eqref{eq:holder} together with the fact that $\Om$ is bounded.
\ep

We first recall the node-averaging interpolator (also called Oswald interpolator) and state the variable-exponent estimate used below. Let $\mathfrak T_h$ be a matching simplicial submesh of $\mathcal M_h$ in the sense of \cite[Definition~1.8]{Di-Pietro.Droniou:20}, and set $m_{\mathrm{av}}\eq\max(1,\widetilde k)$. For every Lagrange node $V$ of $\mathfrak T_h$, let $\mathfrak T_V$ denote the set of simplices sharing $V$. The node-averaging operator
$\mathcal I_{\cst{av},h}^{m_{\mathrm{av}}}:\PP^{\widetilde k}(\T_h)\to \PP^{m_{\mathrm{av}}}(\mathfrak T_h)\cap W_0^{1,1}(\Omega)$ is defined by
\bae\label{eq:Iav:def}
(\mathcal I_{\cst{av},h}^{m_{\mathrm{av}}}v_h)(V)\eq
\begin{cases}
\displaystyle \frac{1}{\#\mathfrak T_V}\sum_{K\in\mathfrak T_V}v_h\restri{K}(V),&V\in\Omega,\\[2mm]
0,&V\in\partial\Omega.
\end{cases}
\eae
For $F\in\F_h^{\rm i}$ shared by $T_1,T_2\in\T_h$, we set
$\llbracket v_h\rrbracket_F\eq v_{T_1}\restri F-v_{T_2}\restri F$, while
$\llbracket v_h\rrbracket_F\eq v_T\restri F$ for $F\in\F_h^{\rm b}$.

\bl[Variable-exponent boundedness of the node-averaging operator]\label{pr:Iav:bound}
Let $q\in\Plogpm(\Omega)$. For every $v_h\in\PP^{\widetilde k}(\T_h)$,
\bae\label{eq:Iav:bound}
\|\nabla_h(v_h-\mathcal I_{\cst{av},h}^{m_{\mathrm{av}}}v_h)\|_{\tL^{q(\cdot)}(\T_h)}
\lesssim
\big\|h_{\F}^{-1/q'(\cdot)}\llbracket v_h\rrbracket\big\|_{L^{q(\cdot)}(\F_h)}.
\eae
If, in addition, $q_\petp<d$, then
\bae\label{eq:Iav:bound:sob}
\|v_h-\mathcal I_{\cst{av},h}^{m_{\mathrm{av}}}v_h\|_{L^{q^*(\cdot)}(\Omega)}
\lesssim
\big\|h_{\F}^{-1/q'(\cdot)}\llbracket v_h\rrbracket\big\|_{L^{q(\cdot)}(\F_h)}.
\eae
The hidden constants depend only on $\Om$, $d$, $k$, $\kappa$, $q_\petm$, $q_\petp$, and $q_\cst{log}$.
\el

\bp
For a constant exponent, \eqref{eq:Iav:bound} is exactly \cite[Proposition~14]{Botti.Castanon-Quiroz.ea:20}. More precisely, its proof also gives the local $L^q$ estimate \cite[Eq.~(35)]{Botti.Castanon-Quiroz.ea:20}. Setting $e_h\eq v_h-\mathcal I_{\cst{av},h}^{m_{\mathrm{av}}}v_h$, let $\mathcal N(K)$ be the set of Lagrange nodes of degree $m_{\mathrm{av}}$ on a simplex $K\in\mathfrak T_h$, and define its associated face patch by
\bae
\mathfrak F_K\eq\{F\in\F_h:\overline F\cap\mathcal N(K)\ne\varnothing\}.
\eae
At every $V\in\mathcal N(K)$, the difference between a local nodal value and its Oswald average is controlled by telescoping the jumps across mesh faces incident to $V$; for $V\in\partial\Omega$, the boundary-face jump is used. All these faces belong to $\mathfrak F_K$. The fixed-exponent local estimate and the polynomial inverse inequality \eqref{eq:dis:inverse:ineq} give, for the frozen exponent $q_K\eq |K|^{-1}\int_K q(\b x)\,\cst{d}\b x$,
\bse\label{eq:Iav:frozen}
\ba
\|e_h\|_{L^{q_K}(K)}^{q_K}
&\lesssim \sum_{F\in\mathfrak F_K}h_F\|\llbracket v_h\rrbracket\|_{L^{q_K}(F)}^{q_K},\\
\|\GRAD e_h\|_{\tL^{q_K}(K)}^{q_K}
&\lesssim h_K^{-q_K}\|e_h\|_{L^{q_K}(K)}^{q_K}
\lesssim \sum_{F\in\mathfrak F_K}h_F^{1-q_K}\|\llbracket v_h\rrbracket\|_{L^{q_K}(F)}^{q_K}.
\ea
\ese
To pass to $q(\cdot)$, one repeats that proof on each simplex $K\in\mathfrak T_h$ with the frozen exponent $q_K$. The constants are uniform since $q_K\in[q_\petm,q_\petp]$ and $\mathfrak T_h$ is uniformly shape regular. Indeed, the matching, uniformly shape-regular simplicial submesh has uniformly bounded nodal stars. Each $\mathfrak F_K$ contains a uniformly bounded number of mesh faces; for $F\in\mathfrak F_K$, one has $h_F\eqsim h_K$ and $\operatorname{dist}(K,F)\lesssim h_K$. Conversely, each mesh face belongs to at most a uniformly bounded number of such patches. Thus the local constants and the overlap constant are independent of $h$. The mean-exponent equivalence of Theorem~\ref{thm:mean_exp_equiv}, in particular \eqref{eq:mean.exp.diam.equiv}, applied on the simplices {and, in intrinsic dimension $d-1$ with surface measure, on their facets}, gives
$h_F^{1/q'(\b x)-1/q_K'}\eqsim1$ on every patch. More explicitly, if $(a_F)_{F\in\F_h}$ is non-negative, the uniform finite overlap gives
\bae
 \sum_{K\in\mathfrak T_h}\sum_{F\in\mathfrak F_K}a_F
 \lesssim \sum_{F\in\F_h}a_F.
\eae
Writing $J_h\eq\|h_{\F}^{-1/q'(\cdot)}\llbracket v_h\rrbracket\|_{L^{q(\cdot)}(\F_h)}$ and assuming $J_h>0$, the $U_\ge$--$U_<$ splitting \eqref{eq:rhoUplus}--\eqref{eq:mean.exp.upper}, applied on the disjoint union of the weighted patch faces, gives from \eqref{eq:Iav:frozen}
\bae
\sum_{K\in\mathfrak T_h}\|\GRAD e_h/J_h\|_{\tL^{q_K}(K)}^{q_K}
&\lesssim
\sum_{K\in\mathfrak T_h}\sum_{F\in\mathfrak F_K}
h_F^{1-q_K}\|\llbracket v_h\rrbracket/J_h\|_{L^{q_K}(F)}^{q_K}
\lesssim1.
\eae
Applying \eqref{eq:rhoUplus}--\eqref{eq:mean.exp.upper} to the left-hand side yields $\rho_{\tL^{q(\cdot)}(\T_h)}(\GRAD e_h/J_h)\lesssim1$. (The case $J_h=0$ is immediate.) The uniformly finite overlap of the nodal patches and \eqref{eq:mod_norm_equiv} then yield \eqref{eq:Iav:bound}. 

If $q_\petp<d$, we bound $e_h$ itself rather than its gradient: starting from the first estimate in \eqref{eq:Iav:frozen} and using the discrete Lebesgue embedding \eqref{eq:disc:leb:emb}, we have for every $K\in\mathfrak T_h$,
\bae
\|e_h\|_{L^{q^*(\cdot)}(K)}
&\lesssim h_K^{-1}\|e_h\|_{L^{q(\cdot)}(K)}
\lesssim
\big\|h_{\F}^{-1/q'(\cdot)}\llbracket v_h\rrbracket\big\|_{L^{q(\cdot)}(\mathfrak F_K)}.
\eae
Summing the corresponding modular estimates over $K$ and using again the uniformly finite overlap of the nodal patches gives \eqref{eq:Iav:bound:sob}.
\ep

\bt[Discrete Sobolev--Poincar\'e inequality]\label{pr:discrete-sob-poin}
Let {$q\in\Plogpm(\Om)$ and $r \in \Plogp(\Om)$} such that $r(\cdot) \le q^*(\cdot)$ in $\Om$. For all $\u v_h  \in \dU{h,0}{k}$,
  \bae\label{eq:discrete.Sobolev.Poincare}
    \|  v_h \|_{L^{r(\cdot)}(\Omega)} \lesssim \|   \u v_h \|_{q(\cdot),h},
  \eae
  where the hidden constant depends only on $d$, $\Omega$, $k$, $\kappa$, $q_\petm$, $q_\petp$, $r_\petm$, $r_\petp$, $q_\cst{log}$, and $r_\cst{log}$.
  By virtue of \eqref{eq:discrete.Sobolev.Poincare}, for all {$q\in \Plogpm(\Om)$}, $\| {\cdot} \|_{q(\cdot),h}$ is a norm on $\dU{h,0}{k}$, as well as $\| {\cdot} \|_{q_h,h}$ and $\|{\cdot}\|_{q(\cdot),\cst{G},h}$, thanks to \eqref{eq:disc:mod-norm:equiv}.
  \et

\bp
The case of constant exponents is \cite[Theorem~5.3]{Di-Pietro.Droniou.ea:18}.

When $d=1$, we argue directly. Each $T$ is an interval, and the fundamental theorem of calculus, used successively from a boundary endpoint (where the face unknown vanishes), gives
\bae
 \|v_h\|_{L^\infty(\Omega)}
 \le \sum_{T\in\T_h}\int_T |v_T'(x)|\,\cst{d}x
 +\sum_{T\in\T_h}\sum_{F\in\F_T}|v_T(F)-v_F|.
\eae
Let $M\coloneq\|\u v_h\|_{q(\cdot),h}$. The case $M=0$ is immediate. For $M>0$, the definition of this Luxemburg seminorm yields
\bae
 \sum_{T\in\T_h}\int_T\left|\frac{v_T'(x)}{M}\right|^{q(x)}\,\cst{d}x
 +\sum_{T\in\T_h}\sum_{F\in\F_T}h_T
 \left|\frac{v_T(F)-v_F}{h_TM}\right|^{q(F)}\le1.
\eae
Apply variable-exponent H\"older inequality to the volume terms and to the weighted counting measure
$\sum_{T\in\T_h}h_T\mathcal H^0|_{\partial T}$, whose total mass is $2|\Omega|$. Since $1<q_-\le q_+<\infty$ and $\Omega$ is bounded, we obtain $\|v_h\|_{L^\infty(\Omega)}\lesssim M$, hence \eqref{eq:discrete.Sobolev.Poincare} for every admissible $r$.

Assume now $d\ge2$. As in the proof of Proposition~\ref{pr:sob-poin-wir}, it is enough to establish the result when $q_\petp<d$ and $r=q^*$; the general case follows by the same reduction to an exponent below $d$ and by the embedding $L^{q^*(\cdot)}(\Omega)\hookrightarrow L^{r(\cdot)}(\Omega)$ obtained from \eqref{eq:holder}.

Let $\u v_h\in\dU{h,0}{k}$ and recall that $v_h\in\PP^{\widetilde k}(\T_h)$ denotes its broken cell polynomial. Since $\mathcal I_{\cst{av},h}^{m_{\mathrm{av}}}v_h\in W_0^{1,q(\cdot)}(\Omega)$, \eqref{eq:Sobolev.Poincare}, the triangle inequality, and \eqref{eq:Iav:bound}--\eqref{eq:Iav:bound:sob} give
\bse\label{eq:SPW:Iav}
\ba
\|v_h\|_{L^{q^*(\cdot)}(\Omega)}
&\le \|v_h-\mathcal I_{\cst{av},h}^{m_{\mathrm{av}}}v_h\|_{L^{q^*(\cdot)}(\Omega)}
 +\|\mathcal I_{\cst{av},h}^{m_{\mathrm{av}}}v_h\|_{L^{q^*(\cdot)}(\Omega)}\\
&\lesssim
\|\nabla_hv_h\|_{\tL^{q(\cdot)}(\Omega)}
+\big\|h_{\F}^{-1/q'(\cdot)}\llbracket v_h\rrbracket\big\|_{L^{q(\cdot)}(\F_h)}.
\ea
\ese
It remains to control the jump term by the HHO norm. For an interior face $F=T_1\cap T_2$, using the single-valued face unknown $v_F$,
\bae\label{eq:jump:split}
|\llbracket v_h\rrbracket_F|
\le |v_{T_1}-v_F|+|v_{T_2}-v_F|,
\eae
while, for $F\in\F_h^{\rm b}$, $v_F=0$ and $|\llbracket v_h\rrbracket_F|=|v_T-v_F|$. Hence, using \eqref{eq:jump:split}, $h_F\eqsim h_T$ for $F\in\F_T$, and \eqref{eq:sum-power},
\bae\label{eq:Iav:HHOjump}
\rho_{L^{q(\cdot)}(\F_h)}\!\left(h_{\F}^{-1/q'(\cdot)}\llbracket v_h\rrbracket\right)
\lesssim
\sum_{T\in\T_h}h_T\rho_{L^{q(\cdot)}(\partial T)}\!\left(h_T^{-1}(v_{\partial T}-v_T)\right).
\eae
By \eqref{eq:luxembourg}, \eqref{eq:dseminorm:T:mod}, and \eqref{eq:mod_norm_equiv}, \eqref{eq:Iav:HHOjump} implies
\bae\label{eq:Iav:HHOcontrol}
\|\nabla_hv_h\|_{\tL^{q(\cdot)}(\Omega)}
+\big\|h_{\F}^{-1/q'(\cdot)}\llbracket v_h\rrbracket\big\|_{L^{q(\cdot)}(\F_h)}
\lesssim \|\u v_h\|_{q(\cdot),h}.
\eae
Combining \eqref{eq:SPW:Iav} and \eqref{eq:Iav:HHOcontrol} yields \eqref{eq:discrete.Sobolev.Poincare} for $r=q^*$, and therefore for every $r\le q^*$.
\ep

From now on, we will equip $\dU{h,0}{k}$ with the norm $\|{\cdot}\|_{p(\cdot),h}$.

\section{Projectors, discrete operators, and monotonicity}\label{sec:properties}

In the following, to avoid the proliferation of generic constants, for all $a,b\in[0,\infty)$, we write henceforth $a\lesssim b$ (resp. $a\gtrsim b$) for the inequality $a\le Cb$ (resp. $a\ge Cb$) where $C \in (0,\infty)$ is independent of $h$, the functions $u$, $\stress$, $f$ in \eqref{eq:intro.problem}, the parameters $\delta,\sigma_\cst{c},\sigma_\cst{m}$ in Assumption \ref{ass:stress}, and, for local inequalities, of any mesh element $T \in \T_h$ or face $F \in \F_h$ on which the inequality holds, and $\delta_F$ in Assumption \ref{ass:sT}.
We also write $a\eqsim b$ to mean $C^{-1}b \le a \le C b$, which is equivalent to $b \lesssim a \lesssim b$.
The dependencies of the hidden constants are further specified when needed.

\subsection{Variable-exponent approximation properties of the projectors}

The three projector results used below are established and proved in detail in \cite{harnist-romain-a-priori}, following the abstract projector strategy of \cite[Lemmas~2.1--2.2]{Di-Pietro.Droniou:17*1} and the $L^2$-projector approximation result \cite[Theorem~1.45]{Di-Pietro.Droniou:20}. We recall their statements here for convenience. The bounded-projector approximation lemma is included because it underlies the variable-exponent approximation proposition established in the companion a priori analysis; only the resulting approximation and stability properties are invoked below.

\bl[Approximation for bounded projectors]\label{lem:bounded-projector-approx}
Suppose that $U$ is either star-shaped with respect to a ball of radius $\vartheta h_U$, or is the union of a uniformly bounded number of shape-regular conforming simplices of diameter comparable to $h_U$, with a connected full-facet adjacency graph. Let $p(\cdot)\in\Plogpm(U)$, and let $\Pi_\ell:L^{p(\cdot)}(U)\to\PP^\ell(U)$ be a projector reproducing $\PP^\ell(U)$ and such that
\bae\label{eq:abstract-projector-Lp-stability}
 \|\Pi_\ell w\|_{L^{p(\cdot)}(U)}\le C_\Pi\|w\|_{L^{p(\cdot)}(U)}
 \qquad\forall w\in L^{p(\cdot)}(U).
\eae
Assume moreover that the polynomial inverse inequality $|w|_{W^{m,p(\cdot)}(U)}\le C_{\rm I}h_U^{-m}\|w\|_{L^{p(\cdot)}(U)}$ holds for all $w\in\PP^\ell(U)$ and all integers $m\ge1$. Then, for every integer $s\in\{1,\ldots,\ell+1\}$, every $m\in\{0,\ldots,s\}$, and every $v\in W^{s,p(\cdot)}(U)$,
\bae\label{eq:abstract-projector-approx}
 |v-\Pi_\ell v|_{W^{m,p(\cdot)}(U)}
 \lesssim h_U^{s-m}|v|_{W^{s,p(\cdot)}(U)}.
\eae
The hidden constant depends on $d$, $\ell$, $s$, $p_\petm$, $p_\petp$, $p_\cst{log}$, $C_\Pi$, $C_{\rm I}$, and the relevant uniform geometric parameter (either $\vartheta$ or the simplicial submesh regularity bounds).
\el

\bl[Boundedness of orthogonal projectors on polynomial subspaces]\label{bl:boundedness:proj}
Let $U$ satisfy \eqref{eq:mean.exp.regularity} and $h_U\le h_0$ for some fixed $h_0<\infty$, let $p(\cdot)\in\Plogpm(U)$, and let $\mathcal P\subset\PP^\ell(U)$. The $L^2$-orthogonal projector $\pi_{\mathcal P}:L^1(U)\to\mathcal P$ satisfies
\bae\label{eq:L2proj-variable-bounded}
 \|\pi_{\mathcal P}\varphi\|_{L^{p(\cdot)}(U)}
 \lesssim \|\varphi\|_{L^{p(\cdot)}(U)}.
\eae
The hidden constant is independent of the particular subspace $\mathcal P$ and depends only on $d$, $\ell$, the regularity constant in \eqref{eq:mean.exp.regularity}, $p_\petm$, $p_\petp$, $p_\cst{log}$, and $h_0$.
\el

\bpr[Variable-exponent approximation properties]\label{pr:variable-projector-approx}
Let $q\in\Plogpm(\Om)$ (in particular, $q=p$ or $q=p'$), $T\in\T_h$, $1\le n\le\ell+1$, and $0\le m\le n$. For every $v\in W^{n,q(\cdot)}(T)$,
\bae\label{eq:proj:app:norm:T}
 |v-\proj{T}{\ell}v|_{W^{m,q(\cdot)}(T)}
 \lesssim h_T^{n-m}|v|_{W^{n,q(\cdot)}(T)}.
\eae
Moreover, the scaled trace estimate
\bae\label{eq:proj:app:norm:F}
 \big\|h_T^{1/q(\cdot)}(v-\proj{T}{\ell}v)\big\|_{L^{q(\cdot)}(\partial T)}
 \lesssim h_T^{n}|v|_{W^{n,q(\cdot)}(T)}
\eae
holds. The hidden constants are uniform along the mesh family and depend only on $d$, $\kappa$, $\ell$, $n$, $q_\petm$, $q_\petp$, and $q_\cst{log}$.
\epr

\subsection{Properties of the discrete operators}

We have the following properties for the discrete local operators.

\paragraph{Boundedness of the interpolator.}
The following boundedness property for $\Interp{T}{k}$ is proved for fixed exponents in \cite[Proposition 6.24]{Di-Pietro.Droniou:20}.

\bpr[Boundedness of the interpolator]\label{pr:interpolator-bounded}
For all $q \in \Plogpm(\Om)$, all $T \in \T_h$, and all $v \in W^{1,q(\cdot)}(T)$,
\bae\label{eq:I:boundedness:norm}
  \|\Interp{T}{k} v\|_{q(\cdot),T} \lesssim |v|_{W^{1,q(\cdot)}(T)}.
\eae
The hidden constant depends only on $k$, $d$, $\kappa$, $q_-$, $q_+$, and $q_\cst{log}$.
\epr

\bp
As in \cite[Proposition 6.24]{Di-Pietro.Droniou:20}, the gradient part follows from the $W^{1,q(\cdot)}$ stability of the local $L^2$ projector, obtained from \eqref{eq:L2proj-variable-bounded} and \eqref{eq:abstract-projector-approx}. For the cell--face part, use the definition \eqref{eq:I:def} of $\Interp{T}{k}$, the projector boundedness \eqref{eq:L2proj-variable-bounded}, and the scaled trace inequality \eqref{eq:cont:trace:ineq}. Combining the two contributions gives \eqref{eq:I:boundedness:norm}.
\ep

\paragraph{Commutation of the gradient reconstruction.}
By design, the following relation holds (see \cite[Section 7.2.5]{Di-Pietro.Droniou:20}):
for all $T \in \T_h$ and all $v\in W^{1,1}(T)$,
\bae\label{eq:G:proj}
  \dgrad{k}{T} \Interp{T}{k} v = \PROJ{T}{k}(\GRAD v).
\eae

\paragraph{Consistency of the boundary residual operator.}

The boundary residual operator is polynomially consistent (cf. \cite{Di-Pietro.Droniou:20}), i.e., for all $T \in \T_h$,
 \bae\label{eq:sT:polynomial.consistency}
   \dfbres{k}{\partial T}(\Interp{T}{k} w) = 0\qquad\forall w\in\PP^{k+1}(T).
 \eae
The polynomial consistency of this residual is the counterpart of the difference-operator consistency in \cite[Proposition~2.6]{Di-Pietro.Droniou:20}.

\subsection{Strong monotonicity of the discrete diffusion function}

For every $T\in\T_h$, let $\delta_{\partial T}$ denote the function on $\partial T$ defined by
$(\delta_{\partial T})\res{F}=\delta_F$ for all $F\in\F_T$, and set
\bae\label{eq:degeneracy:Dh}
\mathcal D_h
\eq \rho_{L^{\sob(\cdot)}(\Omega)}(\delta)
+\sum_{T\in\T_h}h_T\rho_{L^{\sob(\cdot)}(\partial T)}(\delta_{\partial T}).
\eae

\bl[Strong monotonicity of $\cst{a}_h$]\label{lem:ah:holder.continuity.strong.monotonicity}
Under Assumptions \ref{ass:stress} and \ref{ass:sT}, for all $\u v_h,\u w_h\in\dU{h}{k}$, denoting $\u e_h=\u v_h-\u w_h$,
\bae\label{eq:ah:strong.monotonicity}
\cst{a}_h(\u v_h,\u e_h)-\cst{a}_h(\u w_h,\u e_h)
\gtrsim \sigma_\cst{m}
\Big(
\mathcal D_h
+\rho_{\sob(\cdot),\cst{G},h}(\u v_h)
+\rho_{\sob(\cdot),\cst{G},h}(\u w_h)
\Big)^{1-\frac{\sobmax_\petm}{\sob_\petm}}
\rho_{\sob(\cdot),\cst{G},h}(\u e_h)^{\frac{\sobmax_\petm}{\sob_\petm}}.
\eae
\el

\bp
If $\rho_{\sob(\cdot),\cst{G},h}(\u e_h)=0$, the assertion is immediate. Assume
$\rho_{\sob(\cdot),\cst{G},h}(\u e_h)>0$ and let $\varepsilon \in (0,1]$.
For every $T\in\T_h$, applying the same reasoning as in \eqref{eq:a:sm:1} and \eqref{eq:eps:pm} gives
\bae\label{eq:ah:sm:bulk}
\rho_{\tL^{\sob(\cdot)}(T)}(\dgrad{k}{T}\u e_T)
&\lesssim\varepsilon\Big[
\rho_{L^{\sob(\cdot)}(T)}(\delta)
+\rho_{\tL^{\sob(\cdot)}(T)}(\dgrad{k}{T}\u v_T)
+\rho_{\tL^{\sob(\cdot)}(T)}(\dgrad{k}{T}\u w_T)\Big]\\
&\quad+\varepsilon^{1-\frac{\sobmax_\petm}{\sob_\petm}}{\sigma_\cst{m}^{-1}}
\big(\stress(\cdot,\dgrad{k}{T}\u v_T)-\stress(\cdot,\dgrad{k}{T}\u w_T),
\dgrad{k}{T}\u e_T\big)_T.
\eae
For a face $F\in\F_T$, the same rewriting of \eqref{eq:S:strong.monotonicity} gives, pointwise on $F$,
\bae
|\dfbres{k}{F}\u e_T|^{\sob}
\lesssim
\big(\delta_F^{\sob}+|\dfbres{k}{F}\u v_T|^{\sob}+|\dfbres{k}{F}\u w_T|^{\sob}\big)^{1-\sob/\sobmax}
\Big[\sigma_\cst{m}^{-1}
(\stabfl_F(\cdot,\dfbres{k}{F}\u v_T)-\stabfl_F(\cdot,\dfbres{k}{F}\u w_T))
\dfbres{k}{F}\u e_T\Big]^{\sob/\sobmax}.
\eae
Thus, applying the same reasoning as in \eqref{eq:a:sm:1} and \eqref{eq:eps:pm}, multiplying by $h_T$, and summing over the faces of $T$ yields
\bae\label{eq:ah:sm:face}
h_T\rho_{L^{\sob(\cdot)}(\partial T)}(\dfbres{k}{\partial T}\u e_T)&\lesssim\varepsilon h_T\Big[
\rho_{L^{\sob(\cdot)}(\partial T)}(\delta_{\partial T})
+\rho_{L^{\sob(\cdot)}(\partial T)}(\dfbres{k}{\partial T}\u v_T)
+\rho_{L^{\sob(\cdot)}(\partial T)}(\dfbres{k}{\partial T}\u w_T)\Big]\\
&\qquad+\varepsilon^{1-\frac{\sobmax_\petm}{\sob_\petm}}\sigma_\cst{m}^{-1}
\big[\cst{s}_T(\u v_T,\u e_T)-\cst{s}_T(\u w_T,\u e_T)\big].
\eae
Summing \eqref{eq:ah:sm:bulk} and \eqref{eq:ah:sm:face} over $T \in \T_h$ gives
\bae\label{eq:ah:sm:rho}
\rho_{\sob(\cdot),\cst{G},h}(\u e_h)
&\lesssim \varepsilon\Big(
\mathcal D_h
+\rho_{\sob(\cdot),\cst{G},h}(\u v_h)
+\rho_{\sob(\cdot),\cst{G},h}(\u w_h)\Big)+\varepsilon^{1-\frac{\sobmax_\petm}{\sob_\petm}}\sigma_\cst{m}^{-1}
\big[\cst{a}_h(\u v_h,\u e_h)-\cst{a}_h(\u w_h,\u e_h)\big].
\eae
As in \eqref{eq:eps:choice}, thanks to the linearity of $\dgrad{k}{T}$ and $\dfbres{k}{\partial T}$ and \eqref{eq:sum-power}, we set
\bae\label{eq:ah:eps:choice}
\varepsilon={\theta}
\frac{\rho_{\sob(\cdot),\cst{G},h}(\u e_h)}
{\mathcal D_h
+\rho_{\sob(\cdot),\cst{G},h}(\u v_h)
+\rho_{\sob(\cdot),\cst{G},h}(\u w_h)}{\in(0,1]},
\eae
{where $\theta>0$ is chosen sufficiently small, depending only on the exponent bounds and fixed mesh parameters, to ensure $\varepsilon\le1$ and to absorb the first term in \eqref{eq:ah:sm:rho}.}
Substituting \eqref{eq:ah:eps:choice} in \eqref{eq:ah:sm:rho}, the first term is absorbed and \eqref{eq:ah:strong.monotonicity} follows.
\ep

\section{Well-posedness of the discrete problem}\label{sec:well-posedness}

\bt[{ Discrete well-posedness and variable-exponent stability bound}]\label{thm:well-posedness}
Under Assumptions \ref{ass:stress} and \ref{ass:sT}, there exists a unique solution $\u u_h\in\dU{h,0}{k}$ to \eqref{eq:lr.discrete}. Moreover,
\bae\label{eq:discrete.solution:bounds}
&\rho_{\sob(\cdot),\cst{G},h}(\u u_h)\lesssim\max(1,\sigma_\cst{m}^{-p'_\petm},\sigma_\cst{m}^{-p'_\petp})
\max\!\left(
\rho_{L^{\sob'(\cdot)}(\Om)}(f)^{p'_\petp/p'_\petm},
\rho_{L^{\sob'(\cdot)}(\Om)}(f)^{p'_\petm/p'_\petp},
\mathcal D_h
\right).
\eae
{The hidden constant is independent of $h$, $\sigma_\cst{m}$, $\sigma_\cst{c}$, and the degeneracy parameters; it may depend on $p_\petm$, $p_\petp$, $p_\cst{log}$, $\Omega$, $k$, and the mesh regularity parameter $\kappa$.}
\et

\bp 
(i) \emph{Existence.} Existence of a solution to this problem for a fixed $h$ can be proved by adapting the arguments of \cite[Theorem 4.5]{Di-Pietro.Droniou:17}.
  Specifically, equip $\dU{h,0}{k}$ with an inner product $(\cdot,\cdot)_{W,h}$, denote by $\|{\cdot}\|_{W,h}$ the induced norm, and let $\Phi_h : \dU{h,0}{k} \to \dU{h,0}{k}$ be such that, for all $\u w_h\in\dU{h,0}{k}$, $(\Phi_h(\u w_h),\u v_h)_{W,h} = \cst{a}_h(\u w_h,\u v_h)$ for all $\u v_h\in\dU{h,0}{k}$.
  For all $\u v_h\in\dU{h,0}{k}$ such that $\rho_{\sob(\cdot),\cst{G},h}(\u v_h) \ge \max(1,\mathcal D_h)$ (which happens for $\|\u v_h\|_{\sob(\cdot),\cst{G},h}$ large enough), 
\bae
  (\Phi_h(\u v_h),\u v_h)_{W,h} = \cst{a}_h(\u v_h,\u v_h) ~&\refover{\eqref{eq:ah:strong.monotonicity}}\gtrsim {\sigma_\cst{m}} (\mathcal D_h+\rho_{\sob(\cdot),\cst{G},h}(\u v_h))^{1-\frac{\sobmax_\petm}{\sob_\petm}}\rho_{\sob(\cdot),\cst{G},h}(\u v_h)^{\frac{\sobmax_\petm}{\sob_\petm}} \\
  &\gtrsim {\sigma_\cst{m}}\rho_{\sob(\cdot),\cst{G},h}(\u v_h)\refovern{\eqref{eq:modular-norm-equiv}}\ge {\sigma_\cst{m}}\|\u v_h\|_{\sob(\cdot),\cst{G},h}^{\sob_\petm} \ge {\sigma_\cst{m}} C_h^{\sob_\petm} \| \u v_h\|_{W,h}^{\sob_\petm},
\eae
   where $C_h \in (0,\infty)$ comes from the equivalence of the norms $\|{\cdot}\|_{\sob(\cdot),\cst{G},h}$ and $\|{\cdot}\|_{W,h}$ on the finite-dimensional $\dU{h,0}{k}$.
  Hence, we can conclude as in the proof of Theorem~\ref{prop:a-priori} and there exists a solution to \eqref{eq:lr.discrete}.
  \par\medskip

  (ii) \emph{A priori estimates.}
Suppose $\rho_{\sob(\cdot),\cst{G},h}(\u u_h)\ge\mathcal D_h$. As in the continuous case, define $\#\eq\operatorname{sign}(1-\|\u u_h\|_{\sob(\cdot),\cst{G},h})$, with $\operatorname{sign}(0)=+1$ for the choice of index. By \eqref{eq:ah:strong.monotonicity}, the discrete problem \eqref{eq:lr.discrete}, Hölder's inequality \eqref{eq:holder}, the discrete Sobolev--Poincaré bound \eqref{eq:discrete.Sobolev.Poincare}, and the modular--norm relations \eqref{eq:modular-norm-equiv} and \eqref{eq:disc:mod-norm:equiv},
{
\bae\label{eq:well-posedness:0}
\sigma_\cst{m}\rho_{\sob(\cdot),\cst{G},h}(\u u_h)
&\lesssim\cst{a}_h(\u u_h,\u u_h)=(f,u_h)_\Om\\
&\lesssim\|f\|_{L^{\sob'(\cdot)}(\Om)}\|\u u_h\|_{\sob(\cdot),\cst{G},h}\\
&\lesssim\|f\|_{L^{\sob'(\cdot)}(\Om)}\rho_{\sob(\cdot),\cst{G},h}(\u u_h)^{1/p_\#}.
\eae
}
For a nonzero discrete modular, rearranging gives
{
\bae\label{eq:proof:apriori:d}
\rho_{\sob(\cdot),\cst{G},h}(\u u_h)
\lesssim\sigma_\cst{m}^{-p'_\#}\|f\|_{L^{\sob'(\cdot)}(\Om)}^{p'_\#}.
\eae
}
Using \eqref{eq:modular-norm-equiv}, and treating the case $\rho_{\sob(\cdot),\cst{G},h}(\u u_h)<\mathcal D_h$ separately, yields \eqref{eq:discrete.solution:bounds}. The hidden constants depend only on the fixed exponent bounds, log-Hölder constant, domain, polynomial degree, and mesh regularity, and are independent of $\sigma_\cst{m}$.
  \par\medskip
  (iii) \emph{Uniqueness.}
  Let $\u u_h,\u u'_h \in \dU{h,0}{k}$ be two solutions of \eqref{eq:lr.discrete} and define $\u e_h = \u u_h-\u u'_h$. Then, 
  \bae
   0 &= (f,e_h)_\Om-(f,e_h)_\Om\\ &\refover{\eqref{eq:lr.discrete}}= \cst{a}_h(\u u_h,\u e_h)-\cst{a}_h(\u u'_h,\u e_h) \\
   &\refover{\eqref{eq:ah:strong.monotonicity}}\gtrsim \sigma_\cst{m}(\mathcal D_h+\rho_{\sob(\cdot),\cst{G},h}(\u u_h)+\rho_{\sob(\cdot),\cst{G},h}(\u u'_h))^{1-\frac{\sobmax_\petm}{\sob_\petm}}\rho_{\sob(\cdot),\cst{G},h}(\u e_h)^{\frac{\sobmax_\petm}{\sob_\petm}}.
  \eae
  Thus, $\rho_{\sob(\cdot),\cst{G},h}(\u e_h) = 0 = \| \u e_h\|_{W,h}$, which implies $\u u_h = \u u'_h$.
\ep

\section{Convergence}
\label{sec:convergence}

{The aim of this section is to prove the convergence result of Theorem~\ref{thm:convergence} under minimal regularity.}

\subsection{Variable-exponent discrete compactness}

\bl[Discrete compactness]\label{lem:discrete.comp}
Let $(\u v_h)_{h\in\mathfrak H} \in (\dU{h,0}{k})_{h\in\mathfrak H}$ be a sequence for which there exists $C \in (0,\infty)$ such that for all $h \in \mathfrak H$,
\bae\label{eq:discretecomp:assump}
\|\u v_h\|_{\sob(\cdot),h} \le C.
\eae 
Then, there exists $v \in W^{1,\sob(\cdot)}_0(\Om)$ such that, up to a subsequence, 
\begin{enumerate}[label=(\roman*)]
\item \label{thm:compac:velocity} $v_h \CV{h}{0} v$ strongly in $L^{q(\cdot)}(\Omega)$ for all {$q \in \Plogp(\Om)$} such that {there exists $\varepsilon_q>0$ with $q(\cdot)+\varepsilon_q\le\sob^*(\cdot)$ a.e. in $\Om$};
\item \label{thm:compac:grad} $\dgrad{k}{h} \u v_h \CV{h}{0} \GRAD  v$ weakly in $\tL^{\sob(\cdot)}(\Omega)$.
\end{enumerate}
\el
\br[Uniform subcriticality]
The pointwise inequality $q(\b x)<p^*(\b x)$ for almost every $\b x\in\Omega$ does not, in general, imply the existence of $\varepsilon_q>0$ satisfying $q(\cdot)+\varepsilon_q\le p^*(\cdot)$ almost everywhere. The uniform gap in Lemma~\ref{lem:discrete.comp} is therefore stated explicitly.
\er

\bp
Using the same reasoning as in the proof of  \cite[Theorem 5.6]{Di-Pietro.Droniou.ea:18}, for all $\bxi \in \RR^d$, we have
\bae
\|v_h(\cdot+\bxi)-v_h\|_{L^1(\RR^d)}
 &\lesssim |\bxi|\left(\sum_{T\in\T_h}\|\GRAD v_h\|_{\tL^1(T)}+\sum_{F\in\F_h}\|\llbracket v_h\rrbracket\|_{L^1(F)}\right)\\
 &\lesssim |\bxi|\|\u v_h\|_{1,h}
 \lesssim |\bxi|\|\u v_h\|_{\sob(\cdot),h}
 \lesssim |\bxi|.
\eae
{Here we used the triangle inequality and H\"older inequality~\eqref{eq:holder}.}
Hence, by the Kolmogorov--Riesz compactness theorem, $(v_h)_{h\in\mathfrak H}$ is relatively compact in $L^1(\Om)$. After extraction, there exists $v\in L^1(\Om)$ such that $v_h\to v$ strongly in $L^1(\Om)$ and almost everywhere in $\Om$.

{Let $q\in\Plogp(\Om)$ and fix $\varepsilon_q>0$ such that $q(\cdot)+\varepsilon_q\le\sob^*(\cdot)$ a.e. in $\Om$. Set $r\coloneq q+\varepsilon_q$. Then $r\in\Plogp(\Om)$ and, for a.e. $\b x\in\Om$,}
\bae\label{eq:exponent:gap}
 r(\b x)=q(\b x)+\varepsilon_q\le\sob^*(\b x).
\eae
The discrete Sobolev inequality \eqref{eq:discrete.Sobolev.Poincare} and \eqref{eq:discretecomp:assump} show that $(v_h)_{h\in\mathfrak H}$ is uniformly bounded in $L^{r(\cdot)}(\Om)$. The almost-everywhere convergence and Fatou's lemma imply $v\in L^{r(\cdot)}(\Om)$, so $(v_h-v)_{h\in\mathfrak H}$ is also bounded in $L^{r(\cdot)}(\Om)$. 
{Set $e_h\coloneq v_h-v$. Then
\bae\label{eq:compactness:interpolation:bounds}
\|e_h\|_{L^1(\Om)}\to0,
\qquad \sup_{h\in\mathfrak H}\rho_{L^{r(\cdot)}(\Om)}(e_h)<\infty.
\eae
Since $r_+<\infty$ and $r(\b x)=q(\b x)+\varepsilon_q$, let
\bae
\theta\coloneq\frac{\varepsilon_q}{2(r_+-1)}\in(0,1).
\eae
For almost every $\b x\in\Om$, we have $q(\b x)\le\theta+(1-\theta)r(\b x)$. Splitting the integral into $\{|e_h|\le1\}$ and $\{|e_h|>1\}$ and then applying the classical H\"older inequality with the constant conjugate exponents $1/\theta$ and $1/(1-\theta)$ gives
\bae\label{eq:compactness:interpolation}
\rho_{L^{q(\cdot)}(\Om)}(e_h)
&=\int_{\{|e_h|\le1\}}|e_h|^{q(\b x)}\,\cst d\b x+
  \int_{\{|e_h|>1\}}|e_h|^{q(\b x)}\,\cst d\b x\\
&\le\|e_h\|_{L^1(\Om)}+
 \int_{\{|e_h|>1\}}|e_h|^{\theta}|e_h|^{(1-\theta)r(\b x)}\,\cst d\b x\\
&\lesssim\|e_h\|_{L^1(\Om)}+
 \|e_h\|_{L^1(\Om)}^\theta
 \rho_{L^{r(\cdot)}(\Om)}(e_h)^{1-\theta}
\CV{h}{0}0.
\eae
Consequently, the modular--norm relation \eqref{eq:modular-norm-equiv} implies $\|v_h-v\|_{L^{q(\cdot)}(\Om)}\to0$.}
This proves \emph{(i)}.
It remains to show that $v \in W^{1,\sob(\cdot)}_0(\Om)$ and \emph{(ii)}. Now, we use the same reasoning as in the proof of \cite[Proposition 5.6]{Di-Pietro.Droniou:17} to write for all $\btau \in \b{C}^\infty(\RR^d)$,
{
\bae
(\dgrad{k}{h}\u v_h,\btau)_\Om
&=\I_1+\I_2-(v_h,\DIV\btau)_{\Om},
\eae
where
\ban
\I_1 &\eq (\dgrad{k}{h}\u v_h-\nabla_h v_h,
\btau-\PROJ{h}{k}\btau)_\Om,\\
\I_2 &\eq \sum_{T\in\T_h}
(v_{\partial T}-v_T,(\PROJ{T}{k}\btau-\btau)\cdot\b n_{\partial T})_{\partial T}.
\ean}
By definition of the reconstructed modular \eqref{eq:dseminorm:G:mod} and \eqref{eq:disc:mod-norm:equiv},
\bae\label{eq:dpgrho1}
 \|\dgrad{k}{h}\u v_h\|_{\tL^{p(\cdot)}(\Om)}
 \le \|\u v_h\|_{p(\cdot),\cst{G},h}
 \lesssim \|\u v_h\|_{p(\cdot),h}
 \lesssim 1.
\eae
Thus, $(\dgrad{k}{h}\u v_h)_{h\in\mathfrak H}$ is bounded in $\tL^{p(\cdot)}(\Om)$ and, by reflexivity, converges weakly in $\tL^{p(\cdot)}(\Om)$, up to a subsequence, to some $G\in\tL^{p(\cdot)}(\Om)$.
Therefore, H\"older's inequality \eqref{eq:holder}, the triangle inequality, \eqref{eq:discretecomp:assump}, and \eqref{eq:dpgrho1} yield
\bae
|\I_1| \le (\|\dgrad{k}{h}\u v_h\|_{\tL^{p(\cdot)}(\Om)}+\|\nabla_h v_h\|_{\tL^{p(\cdot)}(\Om)})\|\btau-\PROJ{h}{k}\btau\|_{\tL^{p'(\cdot)}(\Om)}
\refovern{\eqref{eq:discretecomp:assump}, \eqref{eq:proj:app:norm:T}}
\lesssim h|\btau|_{\tW^{1,\sob'(\cdot)}(\Om)}.
\eae
Moving to $\I_2$, for all $T\in\T_h$, apply the variable-exponent H\"older inequality \eqref{eq:holder} on each term $(v_{\partial T}-v_T,(\PROJ{T}{k}\btau-\btau)\cdot\b n_{\partial T})_{\partial T} = h_T(h_T^{-1}(v_{\partial T}-v_T),(\PROJ{T}{k}\btau-\btau)\cdot\b n_{\partial T})_{\partial T}$. The first factor is exactly controlled by the cell--face part of $\|\u v_h\|_{\sob(\cdot),h}$ defined by \eqref{eq:dseminorm:T:mod}. For the second factor, observe that
\bae
 h_T\rho_{L^{\sob'(\cdot)}(\partial T)}(\PROJ{T}{k}\btau-\btau)
 =\rho_{L^{\sob'(\cdot)}(\partial T)}
 \bigl(h_T^{1/\sob'(\cdot)}(\PROJ{T}{k}\btau-\btau)\bigr).
\eae
Since $\btau$ is smooth, the pointwise polynomial approximation estimate $|\PROJ{T}{k}\btau-\btau|\lesssim h_T\|\GRAD\btau\|_{L^\infty(T)}$ yields the required scaled trace modular bound directly. Summing these bounds over the mesh and using \eqref{eq:modular-norm-equiv} gives
\bae
 |\I_2|\lesssim h\|\u v_h\|_{\sob(\cdot),h}\,
 \|\GRAD\btau\|_{L^\infty(\Om)}.
\eae
Hence we can conclude as in \cite[Proposition 5.6]{Di-Pietro.Droniou:17}: letting $h$ tend to $0$ gives $(G,\btau)_\Om = -(v,\DIV\btau)_{\Om}$ and taking $\btau$ compactly supported in $\Om$ yields $G = \GRAD v$. Finally, since $(\GRAD v,\btau)_\Om = -(v,\DIV\btau)_{\Om}$ for all $\btau \in \b{C}^\infty(\RR^d)$, an integration by parts shows that $v\restri{\partial\Om} = 0$, hence $v \in W^{1,\sob(\cdot)}_0(\Om)$.
\ep

\subsection{Convergence under minimal regularity}

\bt[Convergence under minimal regularity]\label{thm:convergence}
  Let $(\u u_h)_{h\in\mathfrak H}$ be a sequence such that, for every $h\in\mathfrak H$, $\u u_h\in\dU{h,0}{k}$ solves \eqref{eq:lr.discrete}.
  Then, under Assumptions \ref{ass:stress} and \ref{ass:sT}, and assuming in addition that
  $\sup_{h\in\mathfrak H}\mathcal D_h<+\infty$, there exists a unique $u \in U$ solving \eqref{eq:lr.weak} such that: \begin{enumerate}[label=(\roman*)]
\item $u_h \CV{h}{0} u$ strongly in $L^{q(\cdot)}(\Omega)$ for every $q\in\Plogp(\Om)$ for which there exists $\varepsilon_q>0$ such that $q(\cdot)+\varepsilon_q\le\sob^*(\cdot)$ almost everywhere in $\Om$;
\item $\dgrad{k}{h}\u u_h \CV{h}{0}\GRAD u$ strongly in $\tL^{\sob(\cdot)}(\Omega)$;
\item $\cst{s}_h(\u u_h,\u u_h)\CV{h}{0}0$.
\end{enumerate}
\et

\bp[Proof of Theorem \ref{thm:convergence}]\label{proof:thm:convergence}
  \textbf{Step 1.} \textit{Existence of a limit}. Since, for all $h \in \mathfrak H$, $\u u_h \in \dU{h,0}{k}$ solves \eqref{eq:lr.discrete}, the modular bound \eqref{eq:discrete.solution:bounds}, the modular--norm relations \eqref{eq:modular-norm-equiv}, and \eqref{eq:disc:mod-norm:equiv} imply that the sequence $(\|\u u_h\|_{\sob(\cdot),h})_{h\in\mathfrak H}$ is bounded uniformly in $h$.
  Thus, Lemma~\ref{lem:discrete.comp} yields, up to a subsequence, the strong convergence (i) and only the weak convergence of $(\dgrad{k}{h}\u u_h)_{h\in\mathfrak H}$ to $\nabla u$ in $\boldsymbol L^{p(\cdot)}(\Omega)$. 

  \textbf{Step 2.} \textit{Identification of the limit.} 
  Since the sequence $(\|\u u_h\|_{\sob(\cdot),h})_{h \in \mathfrak H}$ is bounded uniformly in $h$, recalling \eqref{eq:dpgrho1}, the sequence $(\|\dgrad{k}{h}\u u_h\|_{\tL^{\sob(\cdot)}(\Omega)})_{h \in \mathfrak H}$ is also bounded uniformly in $h$.
  Combined with \eqref{eq:dpgrho1} and the H\"older continuity \eqref{eq:ass:stress:cont} of $\stress$, this implies that the sequence $(\|\stress(\cdot,\dgrad{k}{h}\u u_h)\|_{\tL^{\sob'(\cdot)}(\Omega)})_{h \in \mathfrak H}$ is also bounded uniformly in $h$.
  Therefore, $(\stress(\cdot,\dgrad{k}{h}\u u_h))_{h \in \mathfrak H}$ weakly converges to some $\stress_{u} \in \tL^{\sob'(\cdot)}(\Omega)$ up to a subsequence.
  Furthermore, using \eqref{eq:lr.discrete} together with the definition \eqref{eq:ah} of $\cst{a}_h$, for all $\varphi \in C_\cst{c}^\infty(\Omega)$, letting $\u{\hat\varphi}_h \eq \Interp{h}{k} \varphi$, we get
  \bae\label{eq:convergence:G}
   (\stress(\cdot,\dgrad{k}{h}\u u_h), \dgrad{k}{h}\u{\hat\varphi}_h)_\Om = (f, \proj{h}{\widetilde k}\varphi)_\Om - \cst{s}_h(\u u_h,\u{\hat\varphi}_h).
  \eae
  We next show that the stabilization term vanishes. From \eqref{eq:S:holder.continuity}--\eqref{eq:S:zero}, with $v=0$, and the elementary bound
  $|z|\le(\delta_F^{p(\b x)}+|z|^{p(\b x)})^{1/p(\b x)}$, we obtain
  \bae
   |\stabfl_F(\b x,z)|^{p'(\b x)}
   \lesssim C_{\sigma_\cst{c}}\big(\delta_F(\b x)^{p(\b x)}+|z|^{p(\b x)}\big).
  \eae
  Here $C_{\sigma_\cst{c}}\eq \max\{1,\sigma_\cst{c}^{p'_\petm},\sigma_\cst{c}^{p'_\petp}\}$; hence, by \eqref{eq:degeneracy:Dh}, \eqref{eq:discrete.solution:bounds}, and the uniform bound on $\mathcal D_h$,
  \bae
   \sum_{T\in\T_h}h_T\rho_{L^{p'(\cdot)}(\partial T)}
   \!\left(\stabfl_{\partial T}(\cdot,\dfbres{k}{\partial T}\u u_T)\right)
   \lesssim C_{\sigma_\cst{c}}\big(\mathcal D_h+\rho_{p(\cdot),\cst{G},h}(\u u_h)\big)\lesssim C_{\sigma_\cst{c}}.
  \eae
  On the other hand, polynomial consistency \eqref{eq:sT:polynomial.consistency} gives, element by element,
  \bae
   \dfbres{k}{\partial T}\Interp{T}{k}\varphi
   =\dfbres{k}{\partial T}\Interp{T}{k}
   (\varphi-\proj{T}{k+1}\varphi).
  \eae
  \eqref{eq:I:boundedness:norm}, \eqref{eq:disc:mod-norm:equiv}, and \eqref{eq:proj:app:norm:T} therefore yield
  \bae
   \sum_{T\in\T_h}h_T\rho_{L^{p(\cdot)}(\partial T)}
   \!\left(\dfbres{k}{\partial T}\Interp{T}{k}\varphi\right)
   \longrightarrow0.
  \eae
  Applying the variable-exponent H\"older inequality \eqref{eq:holder} on the weighted disjoint union of the element boundaries consequently gives
  \bae
     |\cst{s}_h(\u u_h,\u{\hat\varphi}_h)|\longrightarrow0.
  \eae
  Thus the stabilization term vanishes for interpolants of smooth test functions.
Also, \eqref{eq:proj:app:norm:T} implies that $\proj{h}{\widetilde k}\varphi \CV{h}{0}{}\varphi$ strongly in $L^{\sob(\cdot)}(\Omega)$ and $\dgrad{k}{h}\u{\hat\varphi}_h \CV{h}{0}{} \GRAD \varphi$ strongly in $\tL^{\sob(\cdot)}(\Omega)$ (thanks to  \eqref{eq:G:proj}).
Thus, passing to the limit along the subsequence in \eqref{eq:convergence:G} yields
  \bae\label{eq:convergence:sigma.u}
    (\stress_{u},\GRAD  \varphi)_\Om = (f ,\varphi)_\Om.
  \eae
  By density, \eqref{eq:convergence:sigma.u} is valid for $\varphi\in U$.
  On the other hand, using the definition \eqref{eq:ah} of $\cst{a}_h$, the fact that $\cst{s}_h(\u u_h,\u u_h) \ge 0$ from \eqref{eq:sh}, \eqref{eq:S:strong.monotonicity}, and \eqref{eq:S:zero}, and \eqref{eq:lr.discrete}, we infer
  \bae\label{eq:convergence:bound:1}
      (\stress(\cdot,\dgrad{k}{h}\u u_h), \dgrad{k}{h}\u u_h)_\Om = \cst{a}_h(\u u_h,\u u_h) - \cst{s}_h(\u u_h,\u u_h) \leq \cst{a}_h(\u u_h,\u u_h) = (f,u_h)_\Om.
  \eae
  Here $(f,u_h)_\Omega\to(f,u)_\Omega$ follows from (i) with $q=p$: indeed $p^*(\b x)-p(\b x)\ge p_\petm^2/d$ when $p(\b x)<d$ (and $p^*(\b x)=\infty$ otherwise), so one can choose $\varepsilon_p=p_\petm^2/(2d)>0$.
  Hence, using the strong monotonicity \eqref{eq:ass:stress:mono} of $\stress$,  we obtain, for all $\boldsymbol{\Lambda} \in \tL^{p(\cdot)}(\Omega)$,
  \bae\label{eq:convergence:bound:2}
      0 &\le (\stress(\cdot,\dgrad{k}{h}\u u_h)-\stress(\cdot,\boldsymbol{\Lambda}),\dgrad{k}{h}\u u_h-\boldsymbol{\Lambda})_\Om \\
      &\refover{\eqref{eq:convergence:bound:1}}\le (f,u_h)_\Om - (\stress(\cdot,\dgrad{k}{h}\u u_h),\boldsymbol{\Lambda})_\Om -(\stress(\cdot,\boldsymbol{\Lambda}),\dgrad{k}{h}\u u_h-\boldsymbol{\Lambda})_\Om \\
      &\CV{h}{0} (f,u)_\Om - (\stress_{u}, \boldsymbol{\Lambda})_\Om -(\stress(\cdot,\boldsymbol{\Lambda}),\GRAD u-\boldsymbol{\Lambda})_\Om,
    \eae
  where the limit of the third term is justified by the fact that $\stress(\cdot,\boldsymbol{\Lambda}) \in \tL^{p'(\cdot)}(\Omega)$ thanks to the H\"older continuity \eqref{eq:ass:stress:cont} of $\stress$ and the weak convergence $\dgrad{k}{h}\u u_h \rightharpoonup \GRAD u$ in $\tL^{p(\cdot)}(\Omega)$.
  It follows from the classical Minty trick \cite{Minty:63,Leray.Lions:65} (see \cite[Theorem 4.6]{Di-Pietro.Droniou:17} and \cite[Theorem 17]{harnist3} concerning its application to HHO methods) that $u$ satisfies the continuous weak problem \eqref{eq:lr.weak}.
   Indeed, using \eqref{eq:convergence:sigma.u} with $\varphi=u$ in the limit form of \eqref{eq:convergence:bound:2} gives, for every $\boldsymbol\Lambda\in\boldsymbol L^{p(\cdot)}(\Omega)$,
   \bae\label{eq:Minty.trick}
   \bigl(\stress_u-\stress(\cdot,\boldsymbol\Lambda),\nabla u-\boldsymbol\Lambda\bigr)_\Omega\ge0.
   \eae
   Taking $\boldsymbol\Lambda=\nabla u\pm t\boldsymbol\Psi$ with arbitrary $\boldsymbol\Psi\in\boldsymbol L^{p(\cdot)}(\Omega)$, dividing by $t>0$, and letting $t\to0$ yields
   \bae
   (\stress_u,\boldsymbol\Psi)_\Omega=(\stress(\cdot,\nabla u),\boldsymbol\Psi)_\Omega,
   \eae
   by continuity of the associated Nemytskii map under the growth assumptions \eqref{eq:ass:stress:cont} and \eqref{eq:ass:stress:zero}. Hence $\stress_u=\stress(\cdot,\nabla u)$ in $\boldsymbol L^{p'(\cdot)}(\Omega)$, and \eqref{eq:convergence:sigma.u} becomes \eqref{eq:lr.weak}.
Hence, $u \in U$ solves the weak formulation \eqref{eq:lr.weak}, and since it is unique by Theorem~\ref{prop:a-priori}, the convergences extend to the whole sequence.
  \par\medskip
  \textbf{Step 3.} \textit{Strong convergence of the reconstructed gradient and convergence of the stabilization function.}
  Passing to the upper limit inequality \eqref{eq:convergence:bound:1}, and using \eqref{eq:convergence:sigma.u} with $\varphi$ replaced by $u$ (by density), we obtain
  \bae \label{eq:convergence:bound.Fatou}
      \limsup\limits_{h \to 0}( \stress(\cdot,\dgrad{k}{h}\u u_h) , \dgrad{k}{h}\u u_h)_\Om \leq (f,u)_\Om = (\stress_{ u}, \GRAD u)_\Om.
    \eae
  Setting $M_h\coloneq \rho_{L^{\sob(\cdot)}(\Om)}(\delta)
   +\rho_{\tL^{\sob(\cdot)}(\Om)}(\dgrad{k}{h}\u u_h)
   +\rho_{\tL^{\sob(\cdot)}(\Om)}(\GRAD u)$, and repeating the argument leading to \eqref{eq:a:strong.monotonicity} for the pair
  $(\dgrad{k}{h}\u u_h,\GRAD u)$ gives
  \bae\label{eq:convergence:limit.0}
   \sigma_\cst{m}M_h^{1-\frac{\sobmax_\petm}{\sob_\petm}}
   \rho_{\tL^{\sob(\cdot)}(\Om)}(\dgrad{k}{h}\u u_h-\GRAD u)^{\frac{\sobmax_\petm}{\sob_\petm}}
   \lesssim
   (\stress(\cdot,\dgrad{k}{h}\u u_h)-\stress(\cdot,\GRAD u),
   \dgrad{k}{h}\u u_h-\GRAD u)_\Om.
  \eae
  The right-hand side has nonpositive superior limit by \eqref{eq:convergence:bound.Fatou}, the weak convergences established above, and the identity $\stress_u=\stress(\cdot,\nabla u)$. On the other hand, since $(M_h)_{h\in\mathfrak H}$ is uniformly bounded and $1-\frac{\sobmax_\petm}{\sob_\petm}\le0$, \eqref{eq:convergence:limit.0} therefore forces
  \bae
   \rho_{\tL^{\sob(\cdot)}(\Om)}(\dgrad{k}{h}\u u_h-\GRAD u)\CV{h}{0}0,
  \eae
  hence $\dgrad{k}{h}\u u_h\to\GRAD u$ strongly in $\boldsymbol L^{p(\cdot)}(\Omega)$.
  It remains to identify the convergence of the flux. From \eqref{eq:ass:stress:cont} with $\b\eta=0$ and \eqref{eq:ass:stress:zero}, one obtains the standard $p(\cdot)$-growth bound for the Nemytskii operator generated by $\stress$. Since $1<p_-\le p_+<\infty$ and $\dgrad{k}{h}\u u_h\to\nabla u$ strongly in $\boldsymbol L^{p(\cdot)}(\Omega)$, the continuity of $\stress(\b x,\cdot)$ and uniform integrability of the corresponding $p'(\cdot)$-powers imply that $\stress(\cdot,\dgrad{k}{h}\u u_h) \to \stress(\cdot,\nabla u)$ strongly in $\b L^{p'(\cdot)}(\Omega)$.
  In particular, we have
  \bae\label{eq:convergence:sh.uh.0}
    \cst{s}_h(\u u_h,\u u_h) \refovern{\eqref{eq:ah}}= (f , u_h)_\Om  - (\stress(\cdot,\dgrad{k}{h}\u u_h),\dgrad{k}{h}\u u_h)_\Om \CV{h}{0} 0.
  \eae
\ep

\section{Numerical results}
\label{sec:numerics-convergence}

This section illustrates Theorem~\ref{thm:convergence}: convergence of the HHO approximation
$\u u_h$ to the unique weak solution $u$ holds under \emph{no regularity
assumption on $u$ beyond} $u\in U$; the manufactured solutions below have nonzero boundary traces, which are handled as described in Remark~\ref{rem:nonhomogeneous}.

\subsection{General setting}   
\label{sec:numerics-convergence-setting}   

\paragraph{Computational domain and meshes.} All tests are performed on the unit
hypercube $\Omega=(0,1)^d$, for $d\in\{2,3,4\}$, with the Carreau--Yasuda
flux of Example~\ref{ex:Carreau--Yasuda}, taken with $\mu\equiv1$,
$\alpha\equiv2$ and $\delta\in\{0,1\}$. The stabilization uses a flux of the
same form, $\stabfl_F(\b x,w)=\mu_F\big(\delta^2+|w|^2\big)^{(p(\b x)-2)/2}w$
with $\mu_F=2(k+1)(k+2)$ and $\delta_F=\delta$. Regardless of $d$, the mesh is a Cartesian grid, uniform
and axis-aligned (quadrilateral, hexahedral, or tesseract cells depending
on $d$; see Section~\ref{sec:numerics-dimension}). Throughout this
section, $h_T$ denotes the side length of the cell $T$, that is
$\mathrm{diam}(T)/\sqrt{d}$. This differs from the geometric diameter by the fixed factor $1/\sqrt{d}$ and hence only rescales the stabilization term $\Delta_{\partial T}^k$ by a uniform constant. The quantity reported
as $h$ on the horizontal axis of every figure in this section is the true
geometric mesh size, $h=\max_{T\in\mathcal{T}_h}\mathrm{diam}(T)$. The HHO scheme is the
one described in Section~\ref{sec:scheme}, in its equal-order form: cell and face unknowns are polynomials of the same degree $k\in\{0,1,2,3\}$, and the potential reconstruction has degree $k+1$.

\paragraph{Solver.} The discrete nonlinear system is solved by a linearized
iteration, wrapped in a continuation in the exponent field. Two
linearizations are used, and the choice between them is dictated by the
range of $p(\cdot)$. We write the flux as
$\stress(\b x,\b\tau)=\nu(\b x,\b\tau)\,\b\tau$, with the nonlinear weight
$\nu(\b x,\b\tau)=(\delta^2+|\b\tau|^2)^{(p(\b x)-2)/2}$. The Picard
iteration freezes $\nu$ at the previous iterate, and thus never involves
the derivative of the flux; this makes it the natural choice near the
lower bound $p\to1$, where, at $\delta=0$, that derivative becomes
singular. Its limitation lies at the other end of the range: by freezing
$\nu$, it overlooks how $\nu$ varies with $\dgrad{k}{h}\u u_h$, an effect
which grows with $p(\b x)-2$ and prevents convergence once $p_+$ exceeds
about~$3$. The Newton iteration, which
differentiates the flux exactly, retains this effect and converges well
beyond that threshold, but for the same reason inherits the singularity of
the derivative, and fails at $\delta=0$ as $p\to1$. The two are therefore
complementary: the Picard iteration is used throughout, and the Newton
iteration takes over on the large-exponent range of
Section~\ref{sec:numerics-extreme-p}, where the former no longer converges.
Wherever both converge, they reach the same discrete solution. The two
iterations read as follows:
\begin{enumerate} 
  \item {Under-relaxed Picard iteration:} the raw update
  $\u u_{h,\star}^{(n+1)}$ solves~\eqref{eq:lr.discrete} with the flux frozen
  at the previous iterate,
  \bae
    \stress^{(n)}(\cdot,\b\tau) \eq \nu\big(\cdot,\dgrad{k}{h}\u u_h^{(n)}\big)\,\b\tau,
  \eae
  and the stabilization frozen in the same way. A damping is needed because the
  weight $\nu(\cdot,\dgrad{k}{h}\u u_h)$ can vary sharply between iterates
  once $p(\b x)$ departs far from $2$:
  \bae
    \u u_h^{(n+1)} = \omega\,\u u_{h,\star}^{(n+1)} + (1-\omega)\,\u u_h^{(n)},
    \qquad
    \omega = \min\left(p_\petm-1,\ \frac{1}{p_\petp-1}\right)^{\!\beta},
  \eae
  where $n\geq0$ is the iteration index and $\beta=\frac{1}{2}$. The
  iteration stops once the relative increment
  $\|\u u_h^{(n+1)}-\u u_h^{(n)}\|/\|\u u_h^{(n)}\|$ drops below $10^{-8}$, where $\|{\cdot}\|$ is the Euclidean norm of the vector of face unknowns.

  \item {Damped Newton iteration:} the flux is differentiated exactly, and
  each step reads
  \bae
    \u u_h^{(n+1)} = \u u_h^{(n)} + \lambda_n\,\u v_h^{(n)},
  \eae
  where $\u v_h^{(n)}$ is the Newton increment, solution of the problem
  linearized at $\u u_h^{(n)}$, in which the frozen flux
  $\stress^{(n)}(\cdot,\b\tau)$ of the Picard iteration, with
  $\b\xi=\dgrad{k}{h}\u u_h^{(n)}$ and $\b\tau=\dgrad{k}{h}\u w_h$, is
  replaced by the exact derivative of the flux
  \bae
    \GRAD\stress(\cdot,\b\xi)\,\b\tau =
    \nu(\cdot,\b\xi)\left[\b\tau + \big(p(\b x)-2\big)
    \frac{\b\xi\cdot\b\tau}{\delta^{2}+|\b\xi|^{2}}\,\b\xi\right],
  \eae
  the stabilization term being linearized in the same way. The damping
  $\lambda_n\in(0,1]$ is halved from $1$ until the residual decreases. The
  iteration stops once the relative increment drops below $10^{-8}$, as for
  the Picard iteration, and the residual has been reduced by $10^{-9}$
  relative to the load.
\end{enumerate} 

\paragraph{$\theta$-continuation.} Either linearization is run inside the same
continuation. Rather than assembling directly at $p(\cdot)$, it is applied
to the auxiliary field
\bae
  p_\theta(\b x) = 2+\theta\,(p(\b x)-2), \qquad \theta:0\to1,
\eae
with $p_0\equiv2$ the linear problem, solved in a single shot, and
$p_1=p(\cdot)$ the actual target, $p_\theta$ interpolating $p(\cdot)$'s
spatial profile linearly in between. $\theta$ is increased from $0$ to $1$,
each step warm-started from the previous solution and refined by adaptive
bisection whenever the iteration stagnates, thereby making the severe
exponent contrasts and singular solutions below tractable, without altering
the problem actually solved once for $\theta=1$.

\paragraph{Exact solutions.} A smooth solution and a family of singular solutions are contrasted below; in
each case, the source term $f$ and the Dirichlet boundary data are computed
from the exact solution.
\begin{enumerate}
  \item A smooth test case, with a nowhere-vanishing gradient:
    \bae
      u(x,y)=e^{x+y}. \label{eq:sol-smooth}
    \eae 
  \item Singular test cases built from the same tunable radial singularity solution,
  \bae
      u(\b x) = |\b x-\b x_0|^\gamma, \qquad \b x_0=(0.5,0.5), \label{eq:sol-singular}
  \eae
  centered at the domain center so that its reduced regularity affects
  the global error. For $\gamma\in(1,2)$, $\nabla u$ is continuous and
  vanishes at $\b x_0$, while the Hessian of $u$ blows up there. Two
  regularity settings are considered: a moderate deficit,
  $\gamma=\frac{3}{2}$, and a severe deficit, $\gamma=1.1$, close to
  the limiting value $\gamma=1$ below which $\nabla u$ is unbounded at
  $\b x_0$.
\end{enumerate}

\paragraph{Error measure.} For $\u v_h\in\dU{h}{k}$, we consider the modular
\bae
  \rho_{\mathrm{e},h}(\u v_h) = \sum_{T\in\T_h}\Big(
    \rho_{\tL^{p(\cdot)}(T)}(\dgrad{k}{T}\u v_T)
    + \mu_F\,h_T\,\rho_{L^{p(\cdot)}(\partial T)}(\dfbres{k}{\partial T}\u v_T)
  \Big),
\eae
and the associated Luxemburg norm
$\|\u v_h\|_{\mathrm{e},h}=\inf\{\lambda>0:\rho_{\mathrm{e},h}(\u v_h/\lambda)\leq1\}$.
The first term of $\rho_{\mathrm{e},h}$ is the modular of the reconstructed
gradient, and the second is the stabilization energy at $\delta=0$. The
quantity reported in every figure below, denoted $\|u-u_h\|_{\mathrm e,h}$
for short, is the discrete error $\|\u u_h-\Interp{h}{k}u\|_{\mathrm{e},h}$.
The integrals are evaluated by quadrature,
with the local value of $p$ at each quadrature node, and the infimum, which has no
closed form for a variable exponent, is computed by bisection.

\br[General meshes]\label{rem:numerics-general-meshes}
All the tests of this section use Cartesian meshes. In the companion
paper~\cite{harnist-romain-a-priori}, which establishes error estimates for
the same scheme, the convergence rates are also verified on general
polygonal meshes, namely Voronoi meshes in $d=2$.
\er

\subsection{Robustness with respect to the regularity of the solution}
\label{sec:numerics-regularity}   

We consider here the effect of the regularity of the solution $u$ on the convergence behavior, using the three exact solutions introduced above (\eqref{eq:sol-smooth}, and \eqref{eq:sol-singular} with $\gamma=3/2$ and $\gamma=1.1$). The two spatially varying exponent fields used in this section are
\bse
\ba
  p(x,y) &= 1.5 + x, \label{eq:law-linear}\\  
  p(x,y) &= 1.5 + \sin(\pi x)\sin(\pi y), \label{eq:law-sinusoidal}
\ea
\ese
which we call the linear and sinusoidal exponent fields, respectively.
Both have the range $[1.5,2.5]$. 

\begin{figure}[htbp]
  \centering
  \includegraphics[width=0.48\textwidth]{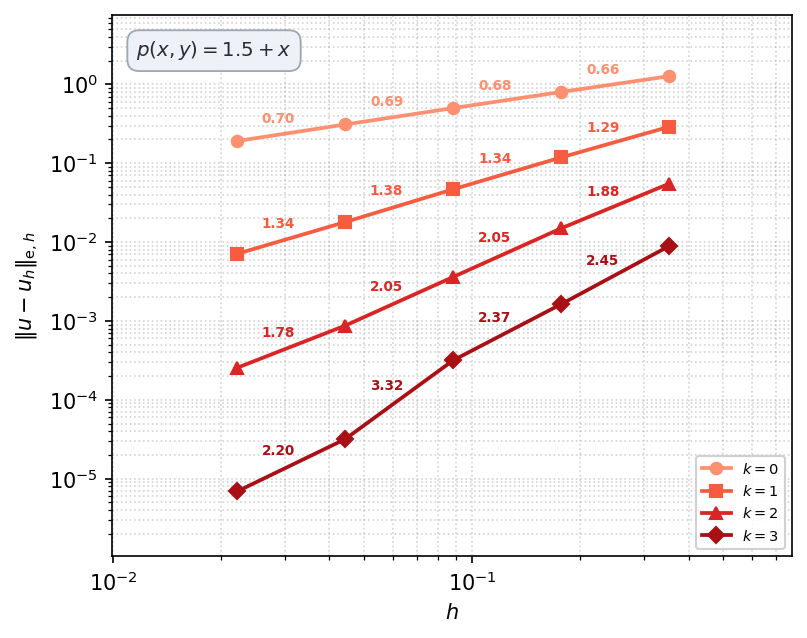}
  \hfill
  \includegraphics[width=0.48\textwidth]{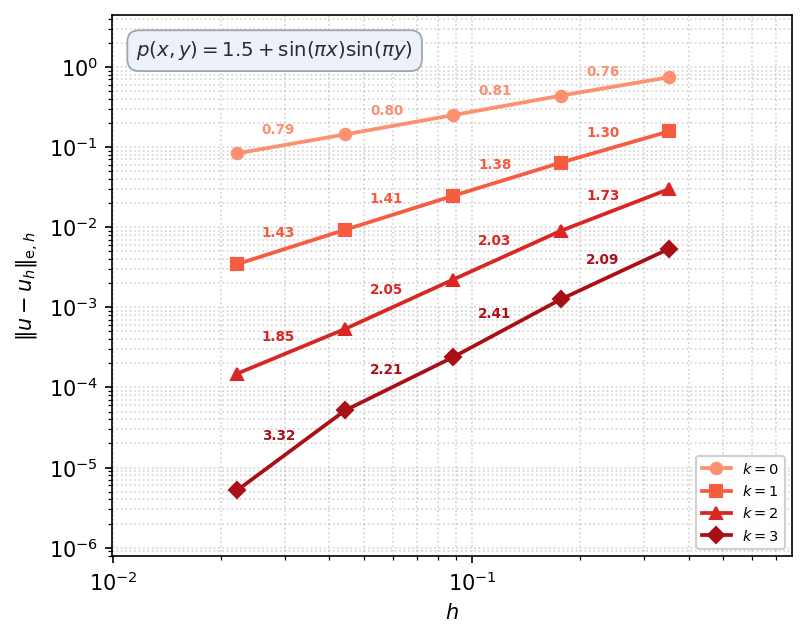}
  \caption{Errors $\|u-u_h\|_{\mathrm e,h}$ as a function of the mesh size
  for the smooth solution~\eqref{eq:sol-smooth}, $\delta=0$.
  \textbf{Left}: $p=1.5+x$. \textbf{Right}: $p=1.5+\sin(\pi x)\sin(\pi
  y)$.}
  \label{fig:optimal}
\end{figure}

We first consider the smooth solution~\eqref{eq:sol-smooth}.
Figure~\ref{fig:optimal} shows the HHO norm error
$\|u-u_h\|_{\mathrm e,h}$ as a function of the mesh size $h$, for both exponent fields and $k\in\{0,1,2,3\}$. As expected, since $u$'s regularity places no ceiling on $k$, the rate keeps improving with $k$ throughout.

\begin{figure}[htbp]
  \centering
  \includegraphics[width=0.48\textwidth]{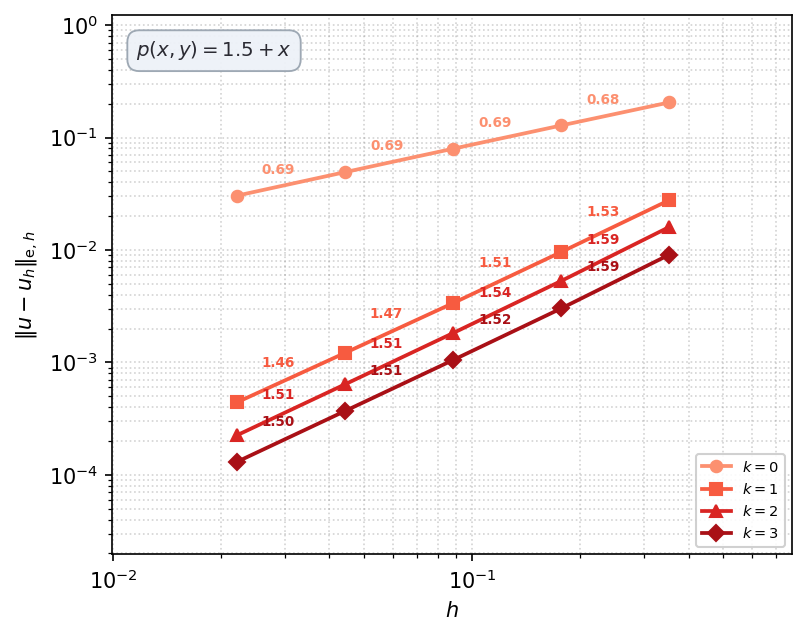}
  \hfill
  \includegraphics[width=0.48\textwidth]{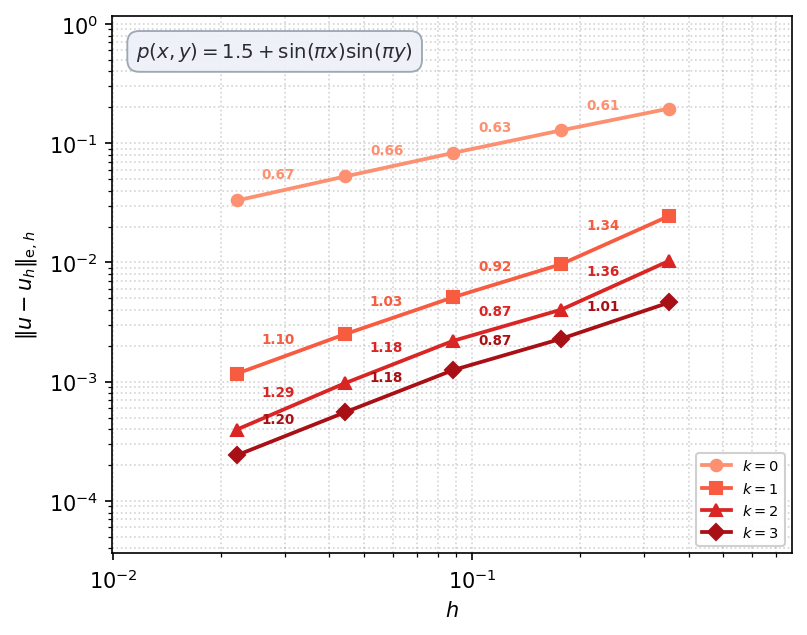}
  \caption{Errors $\|u-u_h\|_{\mathrm e,h}$ as a function of the mesh size
  for the moderate deficit~\eqref{eq:sol-singular} ($\gamma=3/2$), 
  $\delta=0$. \textbf{Left}: $p=1.5+x$. \textbf{Right}:
  $p=1.5+\sin(\pi x)\sin(\pi y)$.}
  \label{fig:mild}
\end{figure}

We consider now the moderate deficit of
regularity~\eqref{eq:sol-singular} ($\gamma=3/2$). Unlike the previous
case, the regularity of this solution limits the attainable rate, whatever $k$. Figure~\ref{fig:mild} reports the same error as a function of $h$. For
the linear exponent field, the rate plateaus starting at $k=1$, at a value close
to $\gamma=3/2$. For the sinusoidal exponent field, the rate also stops improving
from $k=1$ on, at a lower level: on the finest meshes, it lies between $1.1$
and $1.3$ and has not fully settled yet. In both cases,
convergence never fails: only the rate stops improving, once $k$
exceeds what $\gamma=3/2$ can support. 

\begin{figure}[htbp]
  \centering
  \includegraphics[width=0.48\textwidth]{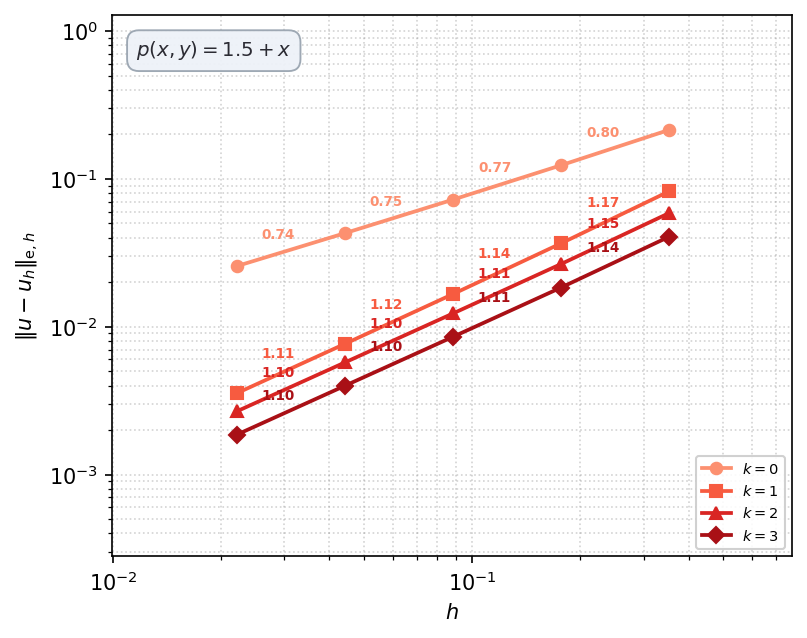}
  \hfill
  \includegraphics[width=0.48\textwidth]{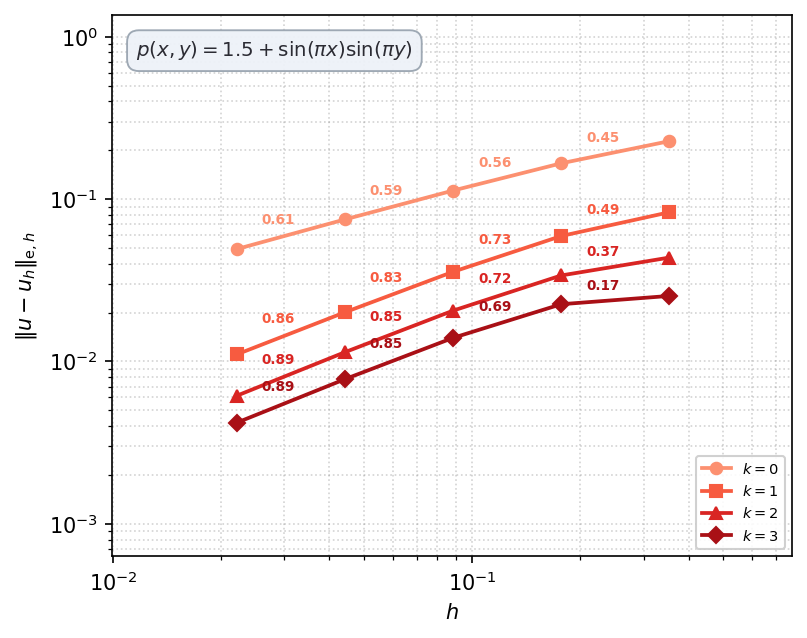}
  \caption{Errors $\|u-u_h\|_{\mathrm e,h}$ as a function of the mesh size
  for the severe deficit~\eqref{eq:sol-singular} ($\gamma=1.1$),
  $\delta=0$. \textbf{Left}: $p=1.5+x$. \textbf{Right}:
  $p=1.5+\sin(\pi x)\sin(\pi y)$.}
  \label{fig:sharp}
\end{figure}

Finally, we evaluate the severe deficit of
regularity~\eqref{eq:sol-singular} ($\gamma=1.1$). Figure~\ref{fig:sharp}
shows the same effect, more pronounced: the rate again plateaus from $k=1$
on, for both exponent fields.
This plateau also sits markedly lower than for the moderate deficit. On the
linear exponent field, the error decreases like $h^{1.1}$ for every $k\geq1$, against
$h^{1.5}$ for $\gamma=3/2$; for the sinusoidal exponent field, it decreases like
$h^{0.9}$, against $h^{1.2}$. A severe deficit of regularity thus degrades
the convergence rates attainable by the scheme as soon as $k\geq1$, while
convergence itself is again observed in every configuration.

\subsection{Robustness with respect to extreme values of $p(\cdot)$}
\label{sec:numerics-extreme-p}

We study here the robustness of Theorem~\ref{thm:convergence} with respect to extreme values
of $p(\cdot)$, using the smooth solution~\eqref{eq:sol-smooth} throughout.
The moderate range $p\in[1.5,2.5]$ of Section~\ref{sec:numerics-regularity}
is replaced by two much wider ranges, each stressing one end of the
admissible interval: $p\in[1.05,2.95]$, which comes within $0.05$ of the
lower bound $p>1$, and $p\in[2,6]$, which leaves the lower end at the
Laplacian and drives the exponent to large values. Both ranges are run at
$\delta=0$, the pure $p(\cdot)$-Laplacian flux, and at $\delta=1$, its
Carreau--Yasuda regularization.

\begin{figure}[htbp]
  \centering
  \includegraphics[width=0.48\textwidth]{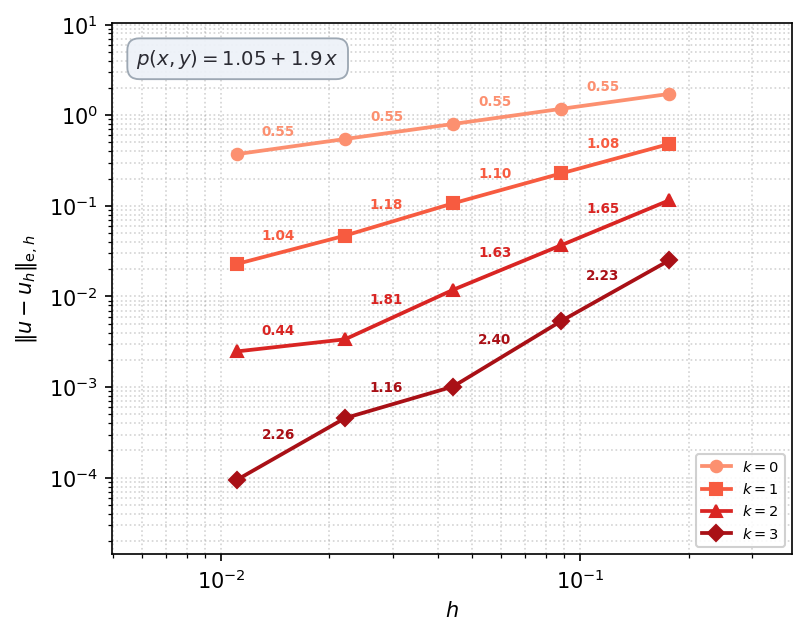}
  \hfill
  \includegraphics[width=0.48\textwidth]{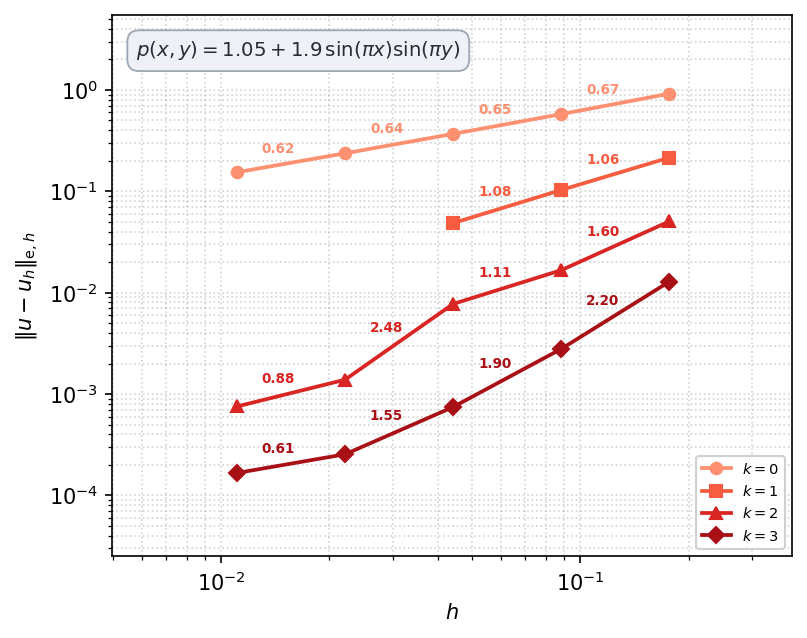}

  \vspace{0.4cm}
  \includegraphics[width=0.48\textwidth]{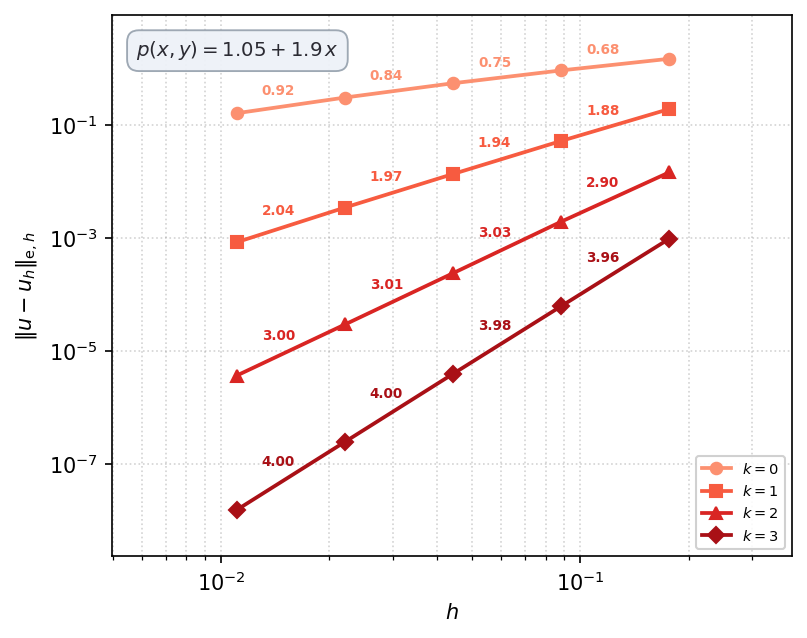}
  \hfill
  \includegraphics[width=0.48\textwidth]{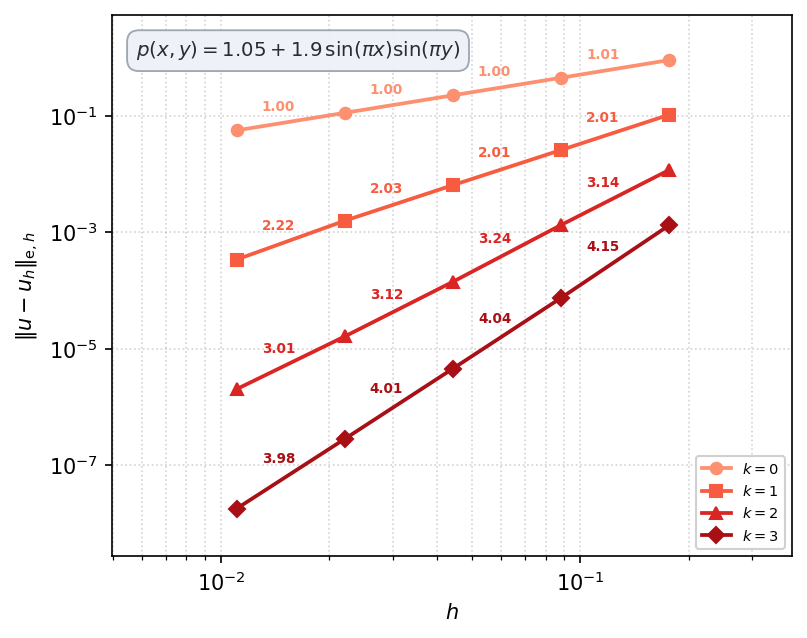}
  \caption{Errors $\|u-u_h\|_{\mathrm e,h}$ as a function of the mesh size for 
  the exponent fields
  \eqref{eq:law-extreme-linear}--\eqref{eq:law-extreme-sinusoidal},
  $p\in[1.05,2.95]$. \textbf{Top}: $\delta=0$. \textbf{Bottom}: $\delta=1$.
  \textbf{Left}: $p=1.05+1.9x$. \textbf{Right}:
  $p=1.05+1.9\sin(\pi x)\sin(\pi y)$. }
  \label{fig:extreme}
\end{figure}

\paragraph{Near the lower admissibility bound.} On this range, the nonlinear
systems are solved with the Picard iteration. We consider the fields
\bse
\ba
  p(x,y) &= 1.05 + 1.9x, \label{eq:law-extreme-linear}\\
  p(x,y) &= 1.05 + 1.9\sin(\pi x)\sin(\pi y), \label{eq:law-extreme-sinusoidal}
\ea
\ese
both ranging over $p\in[1.05,2.95]$. The two fields differ in where the
lower value $p_-=1.05$ is attained: on the single edge $x=0$ for the linear
field, on the whole of $\partial\Omega$ for the sinusoidal field, which makes
the latter the more demanding of the two. Figure~\ref{fig:extreme} reports
the error $\|u-u_h\|_{\mathrm e,h}$ against the mesh size $h$, for both
fields, $k\in\{0,1,2,3\}$ and both values of $\delta$.

The error decreases under refinement at every degree and on both fields, at
$\delta=0$ as well as at $\delta=1$. The only gap concerns the sinusoidal
law at $k=1$ and $\delta=0$, for which the top-right panel reports no error
on the two finest meshes: there, the $\theta$-continuation of
Section~\ref{sec:numerics-convergence-setting} stalled at
$\theta\approx0.9999$, so that the target exponent field was never reached.
These missing points therefore reflect a limitation of the nonlinear
solver, not of the scheme. At $\delta=0$, the slopes of the curves
$k=2$ and $k=3$ also become irregular on the finest meshes. These last points
are sensitive to the stopping tolerance of the Picard iteration: lowering it
from $10^{-8}$ to $10^{-9}$ shifts them by $10$ to $34\%$. Like the missing
points, these irregularities thus stem from the nonlinear solver.

\begin{figure}[htbp]
  \centering
  \includegraphics[width=0.48\textwidth]{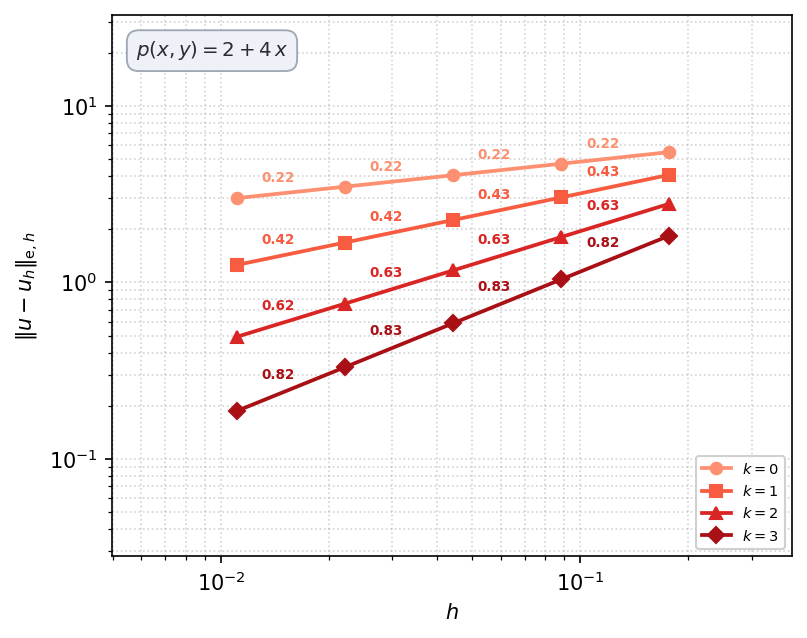}
  \hfill
  \includegraphics[width=0.48\textwidth]{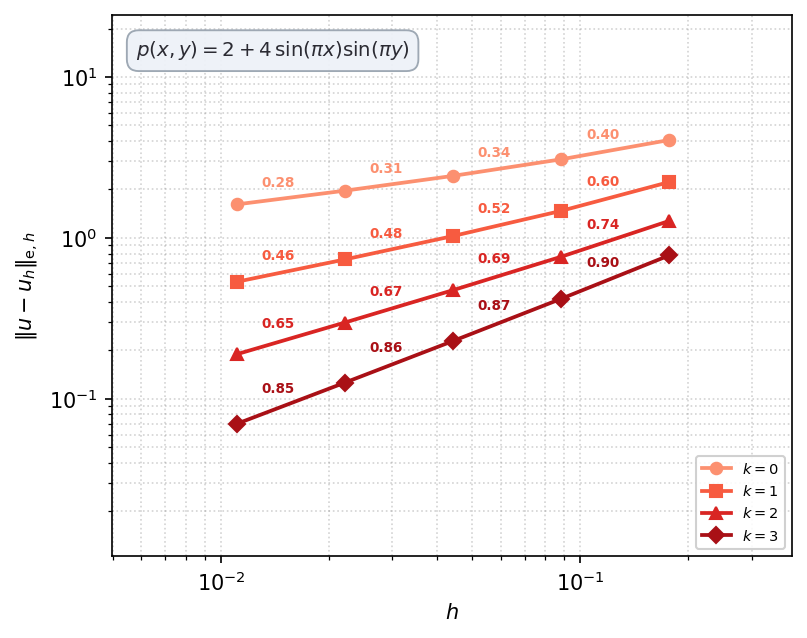}

  \vspace{0.4cm}
  \includegraphics[width=0.48\textwidth]{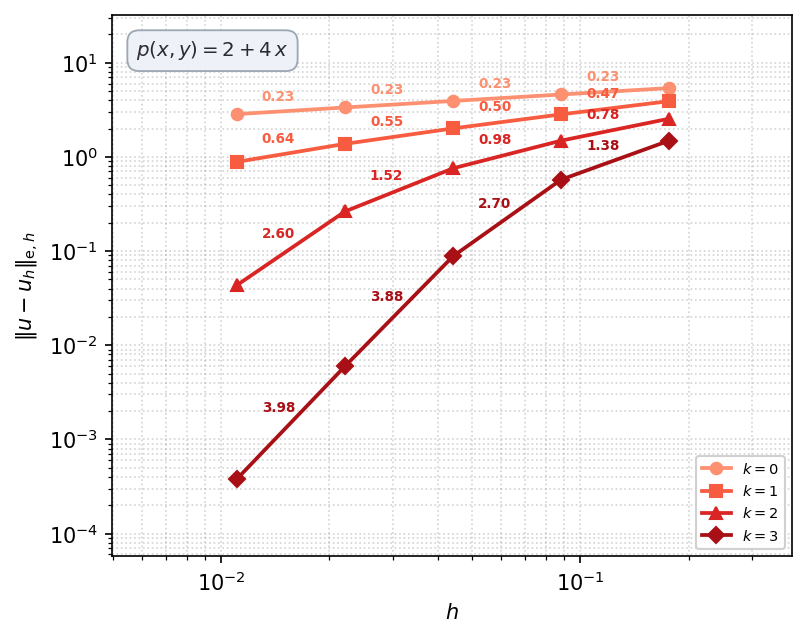}
  \hfill
  \includegraphics[width=0.48\textwidth]{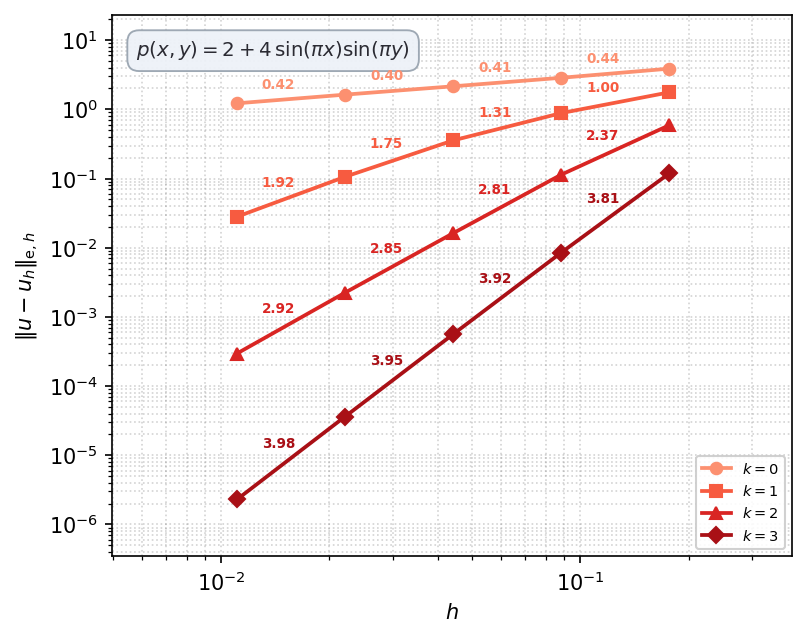}
  \caption{Errors $\|u-u_h\|_{\mathrm e,h}$ as a function of the mesh size for
  the large-exponent fields
  \eqref{eq:law-high-linear}--\eqref{eq:law-high-sinusoidal},
  $p\in[2,6]$. \textbf{Top}:  $\delta=0$. \textbf{Bottom}: $\delta=1$.   
  \textbf{Left}: $p=2+4x$. \textbf{Right}: $p=2+4\sin(\pi x)\sin(\pi y)$.}
  \label{fig:high-p}  
\end{figure}

\paragraph{Large exponents.} On this range, where the Picard iteration no
longer converges, the nonlinear systems are solved with the Newton
iteration. We now push $p_+$ upward with
\bse
\ba
  p(x,y) &= 2 + 4x, \label{eq:law-high-linear}\\
  p(x,y) &= 2 + 4\sin(\pi x)\sin(\pi y), \label{eq:law-high-sinusoidal}
\ea
\ese
both ranging over $p\in[2,6]$. Setting $p_-=2$ excludes the singular regime
$p<2$, so that the large-$p$ end alone is being stressed, and $p_+=6$ lies far
beyond the exponents encountered in applications. Figure~\ref{fig:high-p} 
reports the errors.
Convergence is again observed in all sixteen configurations, at $\delta=0$ as
well as at $\delta=1$. The two values of $\delta$ give errors of the same
size at $k=0$, but they separate as $k$ increases. Convergence itself, which is
what Theorem~\ref{thm:convergence} asserts, is unaffected. 

Taken together, Figures~\ref{fig:extreme} and~\ref{fig:high-p} exhibit
convergence over exponent fields spanning $[1.05,2.95]$ and $[2,6]$, at both
$\delta=0$ and $\delta=1$, from values close to the admissibility bound $p=1$ up to $p_+=6$, in agreement with
Theorem~\ref{thm:convergence}, which only requires $p(\cdot)$ to be
log-H\"older continuous with $1<p_-\leq p_+<\infty$.

\subsection{Extension to higher dimension}
\label{sec:numerics-dimension}

Theorem~\ref{thm:convergence} makes no distinction between $d=2$, $3$ or
$4$. The detailed study of Sections~\ref{sec:numerics-regularity}
and~\ref{sec:numerics-extreme-p} is carried out in $d=2$; in higher
dimension, the scheme and the nonlinear solver apply without any
modification, and we restrict ourselves to a few representative cases in
$d=3$ and $d=4$.

\paragraph{Test case.} The smooth solution~\eqref{eq:sol-smooth} and the
sinusoidal exponent fields are extended to higher dimension by letting them
act on every coordinate. In $d=3$, with coordinates $(x,y,z)$,
\bse
\ba
  u(x,y,z) &= e^{x+y+z}, \label{eq:sol-smooth-3d}\\
  p(x,y,z) &= 1.5 + \sin(\pi x)\sin(\pi y)\sin(\pi z),
    \label{eq:law-sinusoidal-3d}\\
  p(x,y,z) &= 1.05 + 1.9\sin(\pi x)\sin(\pi y)\sin(\pi z),
    \label{eq:law-extreme-sinusoidal-3d}
\ea
\ese
and in $d=4$, with coordinates $(x,y,z,w)$,
\bse
\ba
  u(x,y,z,w) &= e^{x+y+z+w}, \label{eq:sol-smooth-4d}\\
  p(x,y,z,w) &= 1.5 + \sin(\pi x)\sin(\pi y)\sin(\pi z)\sin(\pi w),
    \label{eq:law-sinusoidal-4d}\\
  p(x,y,z,w) &= 1.05 + 1.9\sin(\pi x)\sin(\pi y)\sin(\pi z)\sin(\pi w).
    \label{eq:law-extreme-sinusoidal-4d}
\ea
\ese
In both dimensions, the two exponent fields range over $p\in[1.5,2.5]$ and
$p\in[1.05,2.95]$ respectively, as their two-dimensional
counterparts~\eqref{eq:law-sinusoidal} and~\eqref{eq:law-extreme-sinusoidal}.

\begin{figure}[htbp]
  \centering
  \resizebox{!}{4.25cm}{\begin{tikzpicture}
  \foreach \i in {0,...,4}{
    \draw[gray!65,line width=0.4pt] (\i,0) -- (\i,4);
    \draw[gray!65,line width=0.4pt] (0,\i) -- (4,\i);
  }
  \draw[black,line width=1.1pt] (0,0) rectangle (4,4);
\end{tikzpicture}}
  \hspace{0.5cm}
  \resizebox{!}{4.25cm}{\begin{tikzpicture}
  \foreach \y in {0,...,3}{
    \foreach \z in {0,...,3}{
      \draw[gray!65,line width=0.4pt]
        ({0+0.5*\y},{\z+0.35*\y}) -- ({3+0.5*\y},{\z+0.35*\y});
    }
  }
  \foreach \x in {0,...,3}{
    \foreach \y in {0,...,3}{
      \draw[gray!65,line width=0.4pt]
        ({\x+0.5*\y},{0+0.35*\y}) -- ({\x+0.5*\y},{3+0.35*\y});
    }
  }
  \foreach \x in {0,...,3}{
    \foreach \z in {0,...,3}{
      \draw[gray!65,line width=0.4pt]
        ({\x+0.5*0},{\z+0.35*0}) -- ({\x+0.5*3},{\z+0.35*3});
    }
  }
  \foreach \y in {0,3}{
    \foreach \z in {0,3}{
      \draw[black,line width=1.1pt]
        ({0+0.5*\y},{\z+0.35*\y}) -- ({3+0.5*\y},{\z+0.35*\y});
    }
  }
  \foreach \x in {0,3}{
    \foreach \y in {0,3}{
      \draw[black,line width=1.1pt]
        ({\x+0.5*\y},{0+0.35*\y}) -- ({\x+0.5*\y},{3+0.35*\y});
    }
  }
  \foreach \x in {0,3}{
    \foreach \z in {0,3}{
      \draw[black,line width=1.1pt]
        ({\x+0.5*0},{\z+0.35*0}) -- ({\x+0.5*3},{\z+0.35*3});
    }
  }
\end{tikzpicture}}
  \hspace{0.5cm}
  \resizebox{!}{4.25cm}{\begin{tikzpicture}
  \newcommand{\drawcubegrid}[6]{%
    \pgfmathsetmacro{\Nn}{#1}
    \pgfmathsetmacro{\ss}{#2}
    \pgfmathsetmacro{\ox}{#3}
    \pgfmathsetmacro{\oy}{#4}
    \foreach \y in {0,...,\Nn}{
      \foreach \z in {0,...,\Nn}{
        \draw[#5,line width=0.3pt]
          ({\ox+\ss*(0+0.5*\y)},{\oy+\ss*(\z+0.35*\y)}) --
          ({\ox+\ss*(\Nn+0.5*\y)},{\oy+\ss*(\z+0.35*\y)});
      }
    }
    \foreach \x in {0,...,\Nn}{
      \foreach \y in {0,...,\Nn}{
        \draw[#5,line width=0.3pt]
          ({\ox+\ss*(\x+0.5*\y)},{\oy+\ss*(0+0.35*\y)}) --
          ({\ox+\ss*(\x+0.5*\y)},{\oy+\ss*(\Nn+0.35*\y)});
      }
    }
    \foreach \x in {0,...,\Nn}{
      \foreach \z in {0,...,\Nn}{
        \draw[#5,line width=0.3pt]
          ({\ox+\ss*(\x+0.5*0)},{\oy+\ss*(\z+0.35*0)}) --
          ({\ox+\ss*(\x+0.5*\Nn)},{\oy+\ss*(\z+0.35*\Nn)});
      }
    }
    \foreach \y in {0,\Nn}{
      \foreach \z in {0,\Nn}{
        \draw[#6,line width=1.0pt]
          ({\ox+\ss*(0+0.5*\y)},{\oy+\ss*(\z+0.35*\y)}) --
          ({\ox+\ss*(\Nn+0.5*\y)},{\oy+\ss*(\z+0.35*\y)});
      }
    }
    \foreach \x in {0,\Nn}{
      \foreach \y in {0,\Nn}{
        \draw[#6,line width=1.0pt]
          ({\ox+\ss*(\x+0.5*\y)},{\oy+\ss*(0+0.35*\y)}) --
          ({\ox+\ss*(\x+0.5*\y)},{\oy+\ss*(\Nn+0.35*\y)});
      }
    }
    \foreach \x in {0,\Nn}{
      \foreach \z in {0,\Nn}{
        \draw[#6,line width=1.0pt]
          ({\ox+\ss*(\x+0.5*0)},{\oy+\ss*(\z+0.35*0)}) --
          ({\ox+\ss*(\x+0.5*\Nn)},{\oy+\ss*(\z+0.35*\Nn)});
      }
    }
  }

  \drawcubegrid{2}{1}{0}{0}{gray!60}{black}
  \drawcubegrid{2}{0.5}{0.75}{0.675}{gray!45}{gray!35}

  \foreach \x in {0,2}{
    \foreach \y in {0,2}{
      \foreach \z in {0,2}{
        \draw[black!55,densely dashed,line width=0.45pt]
          ({0+1*(\x+0.5*\y)},{0+1*(\z+0.35*\y)}) --
          ({0.75+0.5*(\x+0.5*\y)},{0.675+0.5*(\z+0.35*\y)});
      }
    }
  }
\end{tikzpicture}}
  \caption{The three conforming mesh families: quadrilaterals ($d=2$, left), hexahedra
  ($d=3$, middle), and a tesseract grid ($d=4$, right); the tesseract is
  drawn as two nested cube grids (a slice at $w=0$ and one at $w=1$)
  connected corner to corner along the fourth axis.}
  \label{fig:mesh-families}
\end{figure}
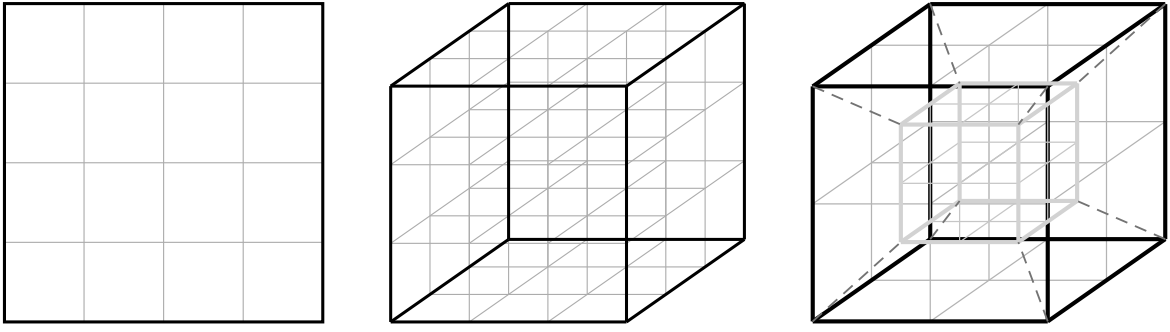

\paragraph{Meshes.} The mesh cells match the dimension being tested:
quadrilaterals for $d=2$, hexahedra for $d=3$, and tesseracts (the
4-dimensional hypercube cell) for $d=4$. Figure~\ref{fig:mesh-families}
illustrates the three families at a coarse level, for legibility, rather
than the actual meshes tested.

\begin{figure}[htbp]
  \centering
  \includegraphics[width=0.48\textwidth]{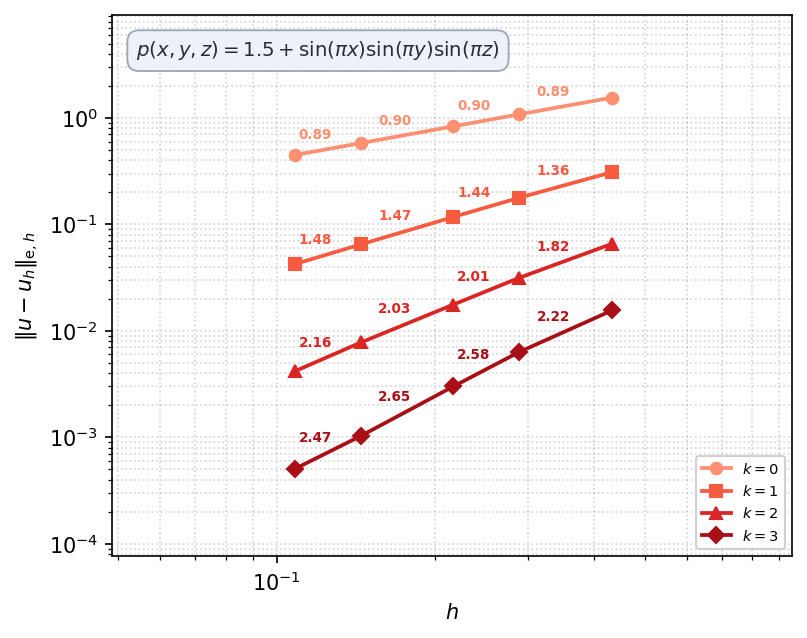}
  \hfill
  \includegraphics[width=0.48\textwidth]{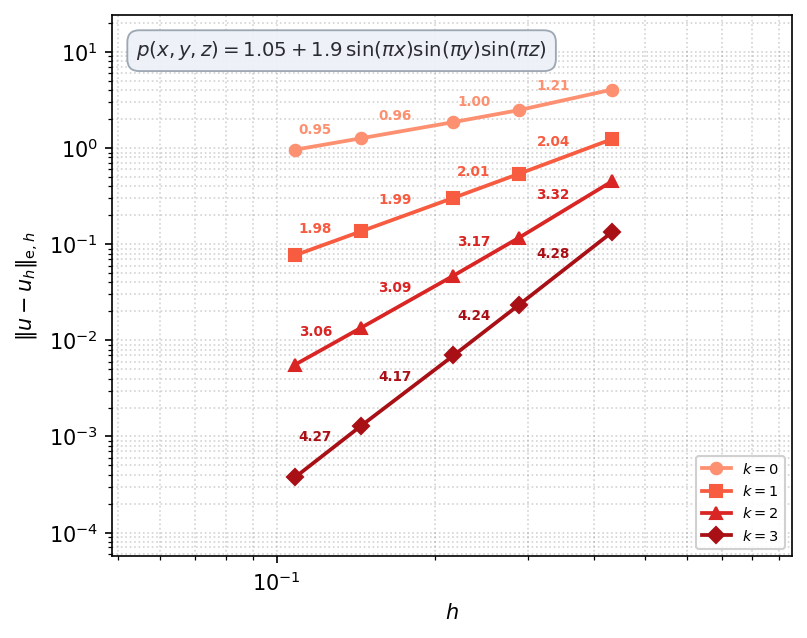}
  \caption{Errors $\|u-u_h\|_{\mathrm e,h}$ as a function of the mesh size
  for the smooth solution~\eqref{eq:sol-smooth-3d}, $d=3$.
  \textbf{Left}: $p$ given by~\eqref{eq:law-sinusoidal-3d}, $\delta=0$.
  \textbf{Right}: $p$ given by~\eqref{eq:law-extreme-sinusoidal-3d},
  $\delta=1$.}
  \label{fig:dim3}
\end{figure}

\begin{figure}[htbp]
  \centering
  \includegraphics[width=0.48\textwidth]{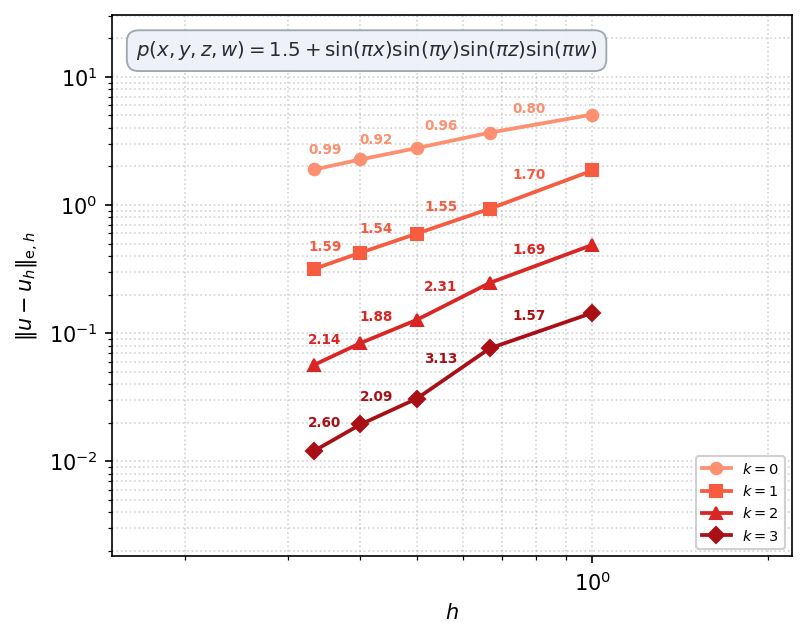}
  \hfill
  \includegraphics[width=0.48\textwidth]{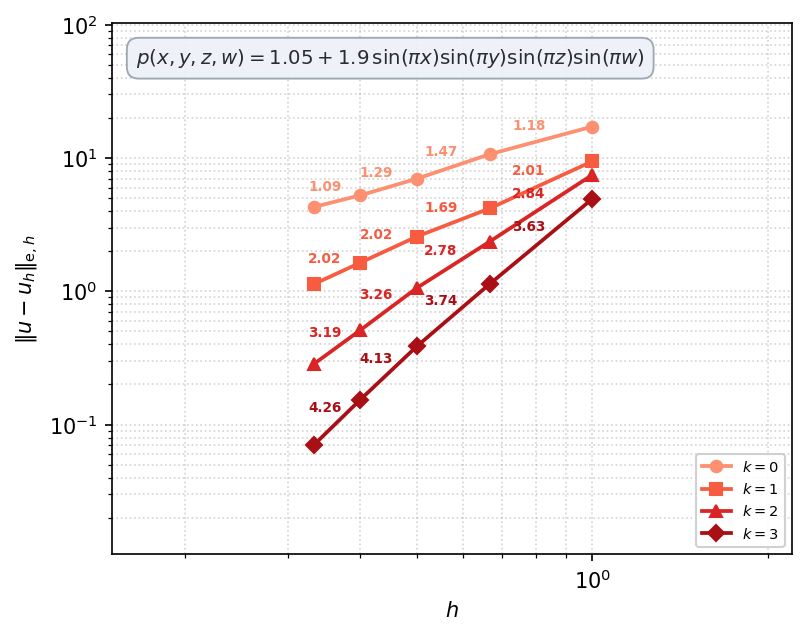}
  \caption{Errors $\|u-u_h\|_{\mathrm e,h}$ as a function of the mesh size
  for the smooth solution~\eqref{eq:sol-smooth-4d}, $d=4$.
  \textbf{Left}: $p$ given by~\eqref{eq:law-sinusoidal-4d}, $\delta=0$.
  \textbf{Right}: $p$ given by~\eqref{eq:law-extreme-sinusoidal-4d},
  $\delta=1$.}
  \label{fig:dim4}
\end{figure}

\paragraph{Results.} The nonlinear systems are solved with the Picard
iteration, at $\delta=0$ on the moderate range and at $\delta=1$ on the
range close to the lower bound, with the same settings as in
Sections~\ref{sec:numerics-regularity} and~\ref{sec:numerics-extreme-p}.
Figures~\ref{fig:dim3} and~\ref{fig:dim4} report the error
$\|u-u_h\|_{\mathrm e,h}$ against the mesh size $h$, for $d=3$ and $d=4$
respectively and $k\in\{0,1,2,3\}$.
The error decreases under refinement for every degree, on both exponent
fields and in both dimensions, and the curves are ordered by $k$ as in
$d=2$. In $d=4$, the observed rates still vary from one mesh to the next
on both ranges: the tesseract meshes being necessarily
coarse, the asymptotic regime is not yet fully reached, although the error
decreases throughout.

\section{Conclusion}
\label{sec:conclusion}

We have proposed and analyzed a Hybrid High-Order discretization of
Leray--Lions problems with a variable exponent $p(\cdot)$ on general
polytopal meshes. A central point of the discrete analysis is the
mean-exponent equivalence of Theorem~\ref{thm:mean_exp_equiv}, which turns
variable-exponent norms of local polynomials into fixed-exponent norms at
the mean exponent and yields the discrete Lebesgue embedding
\eqref{eq:disc:leb:emb} with the scale factor
$h_U^{d/p_U-d/q_U}$. Together with the HHO reconstructions and stabilization,
this provides the discrete norm equivalences, Sobolev--Poincar\'e inequality,
and compactness mechanism needed for the variable-growth setting.

The method covers general Leray--Lions fluxes satisfying the structural
continuity and monotonicity assumptions, without requiring a variational
structure, and includes both equal-order and mixed-order formulations for
all polynomial degrees $k\ge0$. We have proved well-posedness of the discrete
problem and convergence under the sole natural energy regularity of the
exact solution. More precisely, the discrete potentials converge strongly in
every {uniformly subcritical} variable-exponent Lebesgue space, the reconstructed
gradients converge strongly in $L^{p(\cdot)}$, and the stabilization energy
vanishes; see Theorem~\ref{thm:convergence}.

{ The numerical results confirm error decay over several polynomial degrees and space dimensions, while reduced-regularity tests exhibit the expected saturation in observed rates. In the extreme-exponent experiments, the missing $k=1$ data and the irregular slopes observed at higher order are associated with nonlinear solver limitations rather than with a failure of the discrete convergence theorem. The tests thus illustrate the robustness of the method and the computational limits of the nonlinear solvers.}

Theorem~\ref{thm:convergence} is qualitative. Quantitative estimates under
additional approximation regularity, including localized variable-exponent
rates and a regime-dependent non-degenerate refinement, are developed in the
companion a priori error paper~\cite{harnist-romain-a-priori}.

\section*{Acknowledgments}
A. Harnist acknowledges financial support from UTC’s Research Directorate through the Early-Career Researchers Program for the project \emph{``ANACONDAFUN – ANAlyses CONtinues et Discrétisations Avancées de FlUides Non-newtoniens''}.

\appendix
\section{Variable-exponent Lebesgue and Sobolev background}\label{app:ve-background}
This appendix collects the standard functional-analytic material used throughout the paper.

\subsection{Lebesgue and Sobolev spaces with variable exponents}

Let $U\subset\RR^d$ be a bounded measurable set. For volume subsets, we use the exponent inherited from $\Om$. On mesh faces, including boundary faces, any exponent used in surface integrals is the restriction of the unique finite-valued continuous extension to $\overline\Om$ of an exponent in $\Plogp(\Om)$; intrinsic surface measure is used there. We write $\Plog(U)$ for the class of measurable exponents $p:U\to[1,\infty)$ satisfying the local log-H\"older condition: there exists $p_\cst{log}\in(0,\infty)$ such that, for all distinct $\b x,\b y\in U$,
\bae
 \left|\frac1{p(\b x)}-\frac1{p(\b y)}\right|
 \le \frac{p_\cst{log}}{\log(e+1/|\b x-\b y|)}.\label{eq:Plog:1}
\eae
\br[Unbounded domains]
For exponents defined on an unbounded set, one generally supplements local log-H\"older continuity by a decay condition at infinity, involving a limit $p_\infty$ and the factor $\log(e+|\b x|)$; see \cite[Definition~4.1.1]{diening2011lebesgue}. This additional condition is not needed here, since $\Omega$ and all its subsets used below are bounded. The exponent class can be defined on measurable sets $U$, without requiring them to be open; Sobolev and trace results are applied only on domains with the requisite geometry.
\er
For $p\in\Plog(U)$, we define for every measurable function $f : U \to \RR$, the $L^{p(\cdot)}(U)$-modular and $L^{p(\cdot)}(U)$-norm
\bae
\rho_{L^{p(\cdot)}(U)}(f) &\coloneq \int_U |f(\b x)|^{p(\b x)}\,\cst{d}\b x, \\
\|f\|_{L^{p(\cdot)}(U)} &\coloneq \inf\{\lambda > 0 : \rho_{L^{p(\cdot)}(U)}(f/\lambda) \le 1 \}.
\eae 
\br[Luxemburg norm]
Notice that when $p(\cdot)$ is constant, we have $\|f\|_{L^{p}(U)} =  \rho_{L^{p}(U)}(f)^{1/p}$. When $p(\cdot)$ is non-constant with respect to the integration variable, we need to use the Luxemburg norm as above. In the following, for every modular $\rho_\star(\cdot)$, we will define its corresponding Luxemburg norm 
\bae\label{eq:luxembourg}
\|{\cdot}\|_\star \eq \inf\{\lambda > 0 : \rho_\star({\cdot}/\lambda) \le 1 \}.
\eae  
\er
For a multi-index $\alpha=(\alpha_1,\ldots,\alpha_d)\in\mathbb N^d$, we set $|\alpha|\eq \sum_{j=1}^d\alpha_j$ and $\partial_\alpha\eq \partial_{x_1}^{\alpha_1}\cdots\partial_{x_d}^{\alpha_d}$.
For every integer $m \in \NN$, we define the $W^{m,p(\cdot)}(U)$-modular and semi-modular such that for every measurable function $f : U \to \RR$,
\bae
\rho_{W^{m,p(\cdot)}(U)}(f) &\coloneq \sum_{|\alpha| \le m} \rho_{L^{p(\cdot)}(U)}(\partial_\alpha f) \quad \text{ and } \quad \rho_{\dot W^{m,p(\cdot)}(U)}(f) &\coloneq \sum_{|\alpha| = m} \rho_{L^{p(\cdot)}(U)}(\partial_\alpha f).
\eae
The Lebesgue spaces are defined for measurable $U$, whereas the Sobolev spaces below are understood on open $U$, with distributional derivatives.
We define the variable-exponent Lebesgue space $L^{p(\cdot)}(U)$ made of the measurable functions $f : U \to \RR$ such that $\|f\|_{L^{p(\cdot)}(U)} < \infty$; we analogously define the variable-exponent Sobolev space $W^{m,p(\cdot)}(U)$.
These definitions extend to vector-valued functions $\b f : U \to \RR^d$ component by component, and the corresponding spaces are denoted in bold, e.g., $\tL^{p(\cdot)}(U)$.
The scalar inequality and the proposition below are standard; see, e.g., \cite{diening2011lebesgue}.

We recall the generalized Young inequality applied on scalars: for all $x,y \in [0,\infty)$, all $s\in[1,\infty)$ and $p,q\in(1,\infty)$ with $\frac{1}{s}=\frac{1}{p}+\frac{1}{q}$, and all $\varepsilon \in (0,\infty)$,
\bae\label{eq:gener:young:scal}
(xy)^s \le \varepsilon x^p + \varepsilon^{-q/p}y^q. 
\eae

\bpr
For all $p(\cdot),q(\cdot),s(\cdot) \in \Plog(U)$ such that $\frac{1}{s(\cdot)} = \frac{1}{p(\cdot)}+\frac{1}{q(\cdot)}$ in $U$, and all $f \in L^{p(\cdot)}(U)$, $g \in L^{q(\cdot)}(U)$, we have the H\"older and Young inequalities
\bse
\ba
\|fg\|_{L^{s(\cdot)}(U)} &\le 2\|f\|_{L^{p(\cdot)}(U)}\|g\|_{L^{q(\cdot)}(U)},\label{eq:holder}\\
\rho_{L^{s(\cdot)}(U)}(fg) &\le \rho_{L^{p(\cdot)}(U)}(f)+\rho_{L^{q(\cdot)}(U)}(g).\label{eq:young}
\ea
\ese
Moreover, we can deduce the generalized Young inequality writing for all $\varepsilon \in (0,\infty)$,
\bae\label{eq:young:gener}
\rho_{L^{s(\cdot)}(U)}(fg) \le \varepsilon\rho_{L^{p(\cdot)}(U)}(f)+\rho_{L^{q(\cdot)}(U)}(\varepsilon^{-1/p(\cdot)}g).
\eae 
\epr

\subsection{Bounded variable-exponent class}

In the following, we will denote for all $p(\cdot) \in \Plog(U)$ the constants
\bae
\sob_\petm = \essinf_U p \in [1,\infty) \quad \text{ and } \quad \sob_\petp = \esssup_U p \in [1,\infty].
\eae
We define $\Plogm(U)$ (resp. $\Plogp(U)$) as the set of all $p(\cdot) \in \Plog(U)$ for which $\sob_\petm > 1$ (resp. $\sob_\petp < \infty$), and we define $\Plogpm(U) \eq \Plogm(U) \cap \Plogp(U)$.
For $p\in\Plogpm(U)$, the conjugate exponent and the auxiliary exponents are defined by
\ba\label{eq:exp:conj.lower.uper}
 p' \eq \frac{p}{p-1},\qquad \sobmax\coloneq\max(p,2),\qquad \sobmin\coloneq\min(p,2).
\ea
In this case $p',\sobmin,\sobmax\in\Plogpm(U)$.
We denote the conjugates of the endpoint values by $p'_\petm\coloneq(p_\petm)'$ and $p'_\petp\coloneq(p_\petp)'$ (and similarly, locally, by $p'_{T,-}\coloneq(p_{T,-})'$ and $p'_{T,+}\coloneq(p_{T,+})'$). Since taking the conjugate reverses the order, the essential extrema of $p'(\cdot)$ are instead $(p')_\petm=p'_\petp$ and $(p')_\petp=p'_\petm$, whenever $p_\petm>1$.
The following basic modular--norm relations hold for bounded exponents.
\bpr
For all $p(\cdot) \in \Plogp(U)$ and {$0\neq f \in L^{p(\cdot)}(U)$}, we have 
\bse
\ba
&\hspace{-15mm}\rho_{L^{p(\cdot)}(U)}\left(\frac{f}{\|f\|_{L^{p(\cdot)}(U)}}\right) = 1,\label{eq:mod_over_norm_1}\\
\|f\|_{L^{p(\cdot)}(U)} \le 1 \quad &\Longleftrightarrow \quad \rho_{L^{p(\cdot)}(U)}(f) \le 1.\label{eq:equiv_norm_mod_1}
\ea
\ese
\epr
The following result is shown in \cite{diening2011lebesgue}.
\bpr 
For all $p(\cdot) \in \Plogp(U)$ and $f \in L^{p(\cdot)}(U)$, we have 
\bse\label{eq:modular-norm-equiv}
\ba
\min(\|f\|_{L^{p(\cdot)}(U)}^{\sob_\petm},\|f\|_{L^{p(\cdot)}(U)}^{\sob_\petp}) &\le \rho_{L^{p(\cdot)}(U)}(f) \le \max(\|f\|_{L^{p(\cdot)}(U)}^{\sob_\petm},\|f\|_{L^{p(\cdot)}(U)}^{\sob_\petp}),\\
\min(\rho_{L^{p(\cdot)}(U)}(f)^{1/\sob_\petm},\rho_{L^{p(\cdot)}(U)}(f)^{1/\sob_\petp}) &\le \|f\|_{L^{p(\cdot)}(U)} \le \max(\rho_{L^{p(\cdot)}(U)}(f)^{1/\sob_\petm},\rho_{L^{p(\cdot)}(U)}(f)^{1/\sob_\petp}).
\ea
\ese
\epr
We need the following up-to-a-constant version of \eqref{eq:equiv_norm_mod_1}, stated so that the dependence of the constants is explicit.
\bl\label{lm:mod_imp_norm} Let $p(\cdot) \in \Plogp(U)$. Let $v \in L^{p(\cdot)}(U)$ and $\lambda,C \in (0,\infty)$. We have
\bse\label{eq:mod_norm_equiv}
\ba 
\rho_{L^{p(\cdot)}(U)}(v/\lambda) \le C  \ &\Longrightarrow \   \|v\|_{L^{p(\cdot)}(U)}/\lambda \le \max(C^{1/\sob_\petm},C^{1/\sob_\petp}),\\
\rho_{L^{p(\cdot)}(U)}(v/\lambda) \ge C  \ &\Longrightarrow \   \|v\|_{L^{p(\cdot)}(U)} /\lambda\ge \min(C^{1/\sob_\petm},C^{1/\sob_\petp}).
\ea 
Equivalently, 
\ba 
\|v\|_{L^{p(\cdot)}(U)}/\lambda \le C  \ &\Longrightarrow \     \rho_{L^{p(\cdot)}(U)}(v/\lambda) \le \max(C^{\sob_\petm},C^{\sob_\petp}),\\
\|v\|_{L^{p(\cdot)}(U)}/\lambda \ge C  \ &\Longrightarrow \   \rho_{L^{p(\cdot)}(U)}(v/\lambda) \ge \min(C^{\sob_\petm},C^{\sob_\petp}).
\ea 
\ese 
\el

\bp
The case $v=0$ is immediate, so we assume $v\neq0$. Let $\mu \eq \|v\|_{L^{p(\cdot)}(U)}$ and $\theta = \lambda/\mu$. By \eqref{eq:mod_over_norm_1}, $\rho_{L^{p(\cdot)}(U)}(v/\mu)=1$.
If $\rho_{L^{p(\cdot)}(U)}(v/\lambda) \le C$, then
\bae
C \ge \rho_{L^{p(\cdot)}(U)}(\theta^{-1}v/\mu) \ge \min(\theta^{-\sob_\petm},\theta^{-\sob_\petp})\rho_{L^{p(\cdot)}(U)}(v/\mu) = \min(\theta^{-\sob_\petm},\theta^{-\sob_\petp}).
\eae
Hence, $\|v\|_{L^{p(\cdot)}(U)}/\lambda = \theta^{-1} \le \max(C^{1/\sob_\petm},C^{1/\sob_\petp})$.
Likewise, if $\rho_{L^{p(\cdot)}(U)}(v/\lambda) \ge C$, then
\bae
C \le \rho_{L^{p(\cdot)}(U)}(\theta^{-1}v/\mu) \le \max(\theta^{-\sob_\petm},\theta^{-\sob_\petp})\rho_{L^{p(\cdot)}(U)}(v/\mu) = \max(\theta^{-\sob_\petm},\theta^{-\sob_\petp}).
\eae
Hence, $\|v\|_{L^{p(\cdot)}(U)}/\lambda = \theta^{-1} \ge \min(C^{1/\sob_\petm},C^{1/\sob_\petp})$.
Now, if $\|v\|_{L^{p(\cdot)}(U)} \le C\lambda$, then $\theta^{-1} \le C$. So
\bae
\rho_{L^{p(\cdot)}(U)}(v/\lambda) = \rho_{L^{p(\cdot)}(U)}(\theta^{-1}v/\mu) \le \max(C^{\sob_\petm},C^{\sob_\petp})\rho_{L^{p(\cdot)}(U)}(v/\mu) =\max(C^{\sob_\petm},C^{\sob_\petp}).
\eae
Likewise, if $\|v\|_{L^{p(\cdot)}(U)} \ge C\lambda$, then $\theta^{-1} \ge C$ and
\bae
\rho_{L^{p(\cdot)}(U)}(v/\lambda) = \rho_{L^{p(\cdot)}(U)}(\theta^{-1}v/\mu) \ge \min(C^{\sob_\petm},C^{\sob_\petp})\rho_{L^{p(\cdot)}(U)}(v/\mu) =\min(C^{\sob_\petm},C^{\sob_\petp}).
\eae
\ep

\section{Variable-exponent polynomial tools}\label{sec:ve-tools}
This appendix gathers the polynomial variable-exponent estimates used in the discrete analysis. 
Whenever an exponent $q\in\Plogp(\Omega)$ is evaluated on a mesh face, we use its unique finite-valued continuous extension to $\overline\Omega$. The log-H\"older condition alone guarantees a continuous extension of $1/q$, but not necessarily a finite-valued extension of $q$.
Polynomial Lebesgue embeddings on faces are understood in their affine hulls with the intrinsic $(d-1)$-dimensional measure; they follow by the same fixed-exponent argument in dimension $d-1$, under the induced uniform face-regularity assumptions, rather than by treating faces as sets of positive $d$-dimensional measure.

\subsection{Mean exponent equivalence and discrete Lebesgue embedding}

The aim of this section is to generalize the following Lebesgue embedding to variable exponents, cf. \cite[Lemma 5.1]{Di-Pietro.Droniou:17} and \cite[Remark 1.51]{Di-Pietro.Droniou:20}: for all measurable subsets $U$ of $\RR^d$ with inradius $r_U \in (0,\infty)$ and diameter $h_U \in (0,\infty)$, such that
\bae\label{eq:mean.exp.regularity}
\frac{r_U}{h_U} \ge \kappa > 0,
\eae
where $\kappa$ is a lower bound for the inradius-to-diameter ratio of $U$ (for mesh cells, it can be chosen in terms of the mesh regularity parameter of Section~\ref{sec:discrete.setting}), we have for all $p,q \in [1,\infty]$, all $k \in \NN$, and all $v \in \PP^k(U)$,
\bae\label{eq:lebesgue.emb:fixed}
\|v\|_{L^{p}(U)} \eqsim h_U^{\frac{d}{p}-\frac{d}{q}}\|v\|_{L^{q}(U)}
\eae
where the hidden constant, denoted $C_\cst{le}$, depends only on $k$, $d$, $p$, $q$, and $\kappa$.

\br[Dependency of the exponents]\label{rm:dependency_pq}
Let $C_\cst{le}$ denote the constant in~\eqref{eq:lebesgue.emb:fixed}. For finite $p$ and $q$, it depends on the exponents only through powers $p$, $1/p$, $q$, and $1/q$ of fixed positive numbers. Hence a single constant serves for all exponents in a given bounded range, which is how $C_\cst{le}$ is used in the proof of Theorem~\ref{thm:mean_exp_equiv}. For unions of simplices with suitable regularity, the constant does not even depend on $p$ and $q$, cf. \cite[Remark 1.51]{Di-Pietro.Droniou:20}.
\er

\br[Alternative assumptions]
\label{rm:assumption:lebesgue}
Alternatively, the fixed-exponent polynomial embedding~\eqref{eq:lebesgue.emb:fixed} holds uniformly for cells admitting a bounded number of shape-regular simplices of comparable diameter, with constants controlled by the submesh regularity; see \cite[Lemma~1.25]{Di-Pietro.Droniou:20}.
\er

We start with the following theorem.

\bt[Mean-exponent equivalence]\label{thm:mean_exp_equiv}
Let $U\subset\Om$ be a measurable set satisfying \eqref{eq:mean.exp.regularity} (or any assumption implying \eqref{eq:lebesgue.emb:fixed}, cf. Remark \ref{rm:assumption:lebesgue}).
Let $p \in \Plogp(U)$ and denote 
\bae\label{eq:mean.exp.pU}
p_U \eq \frac{1}{|U|}\int_U p(\b x)\,\cst{d}\b x.
\eae
Then, for all $\b x,\b y \in U$,
\bae\label{eq:mean.exp.diam.equiv}
h_U^{\frac{1}{p(\b x)}-\frac{1}{p_U}} \eqsim h_U^{\frac{1}{p(\b x)}-\frac{1}{p(\b y)}} \eqsim h_U^{p(\b x)-p_U} \eqsim h_U^{p(\b x)-p(\b y)} \eqsim 1. 
\eae
Moreover, for all $k \in \NN$, and all $v \in \PP^k(U)$,
\bae\label{eq:mean.exp.norm}
\|v\|_{L^{p(\cdot)}(U)} \eqsim \|v\|_{L^{p_U}(U)},
\eae
where all the hidden constants depend only on $\Om$, $k$, $d$, $p_\cst{log}$, $p_\petp$, and $C_\cst{le}$.
\et 

\br[{Dependence on $\Om$}]
The dependency on $\Om$ cannot be avoided since the log-H\"older continuity controls the oscillation of the exponent but does not provide a uniform control of arbitrarily large diameters.
\er 

\bp
\emph{(i) Diameter equivalence.}
Throughout the proof, the hidden constants in $\lesssim$ and $\gtrsim$ depend only on $\Om$, $p_\petp$, $p_\cst{log}$, the geometric parameter in~\eqref{eq:mean.exp.regularity} or the uniform polynomial-embedding constant of Remark~\ref{rm:assumption:lebesgue}, $d$, and $k$.
First, let $\b x \in U$. Using \eqref{eq:Plog:1}, we have

\bae\label{eq:mean.exp.0}
\left|\frac{1}{p(\b x)}-\frac{1}{p_U}\right|
&=
\frac{1}{p(\b x)p_U|U|}
\left|\int_U\big(p(\b x)-p(\b y)\big)\,\cst{d}\b y\right|\\
&=
\frac{1}{p(\b x)p_U|U|}
\left|\int_U p(\b x)p(\b y)
\left(\frac{1}{p(\b y)}-\frac{1}{p(\b x)}\right)\,\cst{d}\b y\right|
\le
\frac{p_\cst{log}}{\log(e+1/h_U)}.
\eae

Setting $\tau_U\coloneq p_\petp p_\cst{log}/\log(e+1/h_U)$, we obtain
\bae\label{eq:mean.exp.1}
|p(\b x)-p_U|
&=
p(\b x)p_U
\left|
\frac{1}{p(\b x)}-\frac{1}{p_U}
\right|
\le
\frac{p_\petp p_U p_\cst{log}}
{\log(e+1/h_U)}
=
\tau_U p_U.
\eae
Since $\log(1/h_U)\le\log(e+1/h_U)$,
\bae\label{eq:mean.exp.5}
h_U^{-d\tau_U}
&=
\exp\left(
\frac{
d\,p_\petp p_\cst{log}\log(1/h_U)
}{
\log(e+1/h_U)
}
\right)
\le
e^{d p_\petp p_\cst{log}}.
\eae
Moreover, since $U\subset\Om$ with $\Om$ bounded, $h_U\le \operatorname{diam}(\Om)$. Hence, using \eqref{eq:mean.exp.1}, if $h_U\le1$ we have
$|p(\b x)-p_U||\log h_U|\lesssim \log(1/h_U)/\log(e+1/h_U)\lesssim1$, while if $h_U>1$, then $|\log h_U|\le\log(\max(1,\operatorname{diam}(\Om)))$ and $|p(\b x)-p_U|\lesssim1$. Therefore,
\bae\label{eq:hpexp}
\max\left(
h_U^{p(\b x)-p_U},
h_U^{p_U-p(\b x)}
\right)
&=
\exp\left(
|p(\b x)-p_U||\log h_U|
\right)
\lesssim1,
\eae
or, equivalently, $h_U^{p(\b x)-p_U}\eqsim1$, and we infer $h_U^{p(\b x)-p(\b y)} \eqsim 1$ immediately for all $\b y \in U$. The equivalence $h_U^{\frac{1}{p(\b x)}-\frac{1}{p_U}} \eqsim h_U^{\frac{1}{p(\b x)}-\frac{1}{p(\b y)}} \eqsim 1$ uses the same reasoning (with \eqref{eq:mean.exp.0} instead of \eqref{eq:mean.exp.1}).
\par\medskip
\emph{(ii) Norm-equivalence.}
Let $v\in\PP^k(U)$. The result is trivial if $v=0$, so we assume that $v\neq0$ and set
$\tilde v\eq v/\|v\|_{L^{p_U}(U)}$.
Using \eqref{eq:lebesgue.emb:fixed}, there is
$C_\cst{inv}\in[1,\infty)$ depending only on $d$, $k$, the applicable uniform polynomial-embedding constant, and $p_\petp$
(cf. Remark \ref{rm:dependency_pq}) such that
\bae\label{eq:mean.exp.2}
|\tilde v(\b x)|
&\le
\|\tilde v\|_{L^\infty(U)}
\le
C_\cst{inv}h_U^{-d/p_U}
\|\tilde v\|_{L^{p_U}(U)}
=
C_\cst{inv}h_U^{-d/p_U}
\le
C_\cst{inv}\max(1,h_U^{-d/p_U}).
\eae
Therefore, observing that $\tau_Up_U\le p_\petp^2p_\cst{log}$, we obtain
\bae\label{eq:mean.exp.6}
|\tilde v(\b x)|^{|p(\b x)-p_U|}
\refovern{\eqref{eq:mean.exp.2},\eqref{eq:mean.exp.1}}
\le
C_\cst{inv}^{p_\petp^2p_\cst{log}}
\max(1,h_U^{-d\tau_U})
\refovern{\eqref{eq:mean.exp.5}}\lesssim 1.
\eae
We now set $\lambda_U\eq(2|U|)^{1/p_U}$,
$U_\ge\eq\{\b x\in U:|\tilde v(\b x)|\ge\lambda_U^{-1}\}$, and
$U_<\eq\{\b x\in U:|\tilde v(\b x)|<\lambda_U^{-1}\}$.
Since \eqref{eq:mod_over_norm_1} gives $\rho_{L^{p_U}(U)}(\tilde v)=1$ and
$\rho_{L^{p_U}(U_<)}(\tilde v)\le |U|\lambda_U^{-p_U}=1/2$, we have
\bae\label{eq:rhoUplus}
\rho_{L^{p_U}(U_\ge)}(\tilde v)
&=
\rho_{L^{p_U}(U)}(\tilde v)
-
\rho_{L^{p_U}(U_<)}(\tilde v)
\ge
\frac12.
\eae
Applying the uniform fixed-exponent polynomial embedding~\eqref{eq:lebesgue.emb:fixed} to the constant polynomial $1$, with exponents $1$ and $\infty$, gives $|U|\eqsim h_U^d$. Moreover, $|U|\le|\Om|$, so
$|\log(2|U|)|\lesssim\log(e+1/h_U)$. Therefore, using \eqref{eq:mean.exp.1},
\bae\label{eq:mean.exp.lambda}
\max\left(
\lambda_U^{p(\b x)-p_U},
\lambda_U^{p_U-p(\b x)}
\right)
&=
\exp\left(
\frac{
|p(\b x)-p_U|\,|\log(2|U|)|
}{
p_U
}
\right)
\lesssim 1.
\eae
For every $\b x\in U_\ge$, \eqref{eq:mean.exp.6} gives
$|\tilde v(\b x)|^{|p(\b x)-p_U|}\lesssim1$.
Moreover, since $|\tilde v(\b x)|\ge\lambda_U^{-1}$,
$|\tilde v(\b x)|^{-|p(\b x)-p_U|}
\le\max(1,\lambda_U^{|p(\b x)-p_U|})
\lesssim1$ by \eqref{eq:mean.exp.lambda}.
Consequently,
$|\tilde v(\b x)|^{p_U}\lesssim|\tilde v(\b x)|^{p(\b x)}
\lesssim|\tilde v(\b x)|^{p_U}$,
and hence
\bae\label{eq:mean.exp.7}
\rho_{L^{p_U}(U_\ge)}(\tilde v)
&\lesssim
\rho_{L^{p(\cdot)}(U_\ge)}(\tilde v)
\lesssim
\rho_{L^{p_U}(U_\ge)}(\tilde v).
\eae
Combining \eqref{eq:rhoUplus} and \eqref{eq:mean.exp.7} gives
\bae\label{eq:mean.exp.lower}
\rho_{L^{p(\cdot)}(U)}(\tilde v)
&\ge
\rho_{L^{p(\cdot)}(U_\ge)}(\tilde v)
\gtrsim 1.
\eae
On the other hand, for $\b x\in U_<$, we have
$|\tilde v(\b x)|<\lambda_U^{-1}$, and therefore, by
\eqref{eq:mean.exp.lambda},
$|\tilde v(\b x)|^{p(\b x)}
\le
\lambda_U^{-p(\b x)}
=
\lambda_U^{-p_U}
\lambda_U^{p_U-p(\b x)}
\lesssim
\lambda_U^{-p_U}$.
Hence,
\bae\label{eq:mean.exp.minus}
\rho_{L^{p(\cdot)}(U_<)}(\tilde v)
&\lesssim
|U|\lambda_U^{-p_U}
=
\frac12.
\eae
Using \eqref{eq:mean.exp.7}, \eqref{eq:mean.exp.minus}, and
$\rho_{L^{p_U}(U_\ge)}(\tilde v)\le1$, we obtain
\bae\label{eq:mean.exp.upper}
\rho_{L^{p(\cdot)}(U)}(\tilde v)
&=
\rho_{L^{p(\cdot)}(U_\ge)}(\tilde v)
+
\rho_{L^{p(\cdot)}(U_<)}(\tilde v)
\lesssim
\rho_{L^{p_U}(U_\ge)}(\tilde v)
+
1
\lesssim 1.
\eae
Using \eqref{eq:mean.exp.lower}--\eqref{eq:mean.exp.upper}, we conclude that
\bae\label{eq:mean.exp.8}
\rho_{L^{p(\cdot)}(U)}
\left(
v/\|v\|_{L^{p_U}(U)}
\right) \eqsim 1.
\eae
Finally, applying \eqref{eq:mod_norm_equiv} to \eqref{eq:mean.exp.8}, we get \eqref{eq:mean.exp.norm}.
\ep

\bc[Discrete Lebesgue embedding]\label{cor:discrete-lebesgue-emb}
Let $U\subset\Om$ be a measurable set satisfying \eqref{eq:mean.exp.regularity} (or any assumption implying \eqref{eq:lebesgue.emb:fixed}, cf. Remark \ref{rm:assumption:lebesgue}).
Let $p(\cdot),q(\cdot) \in \Plogpm(U)$ and denote 
\bae
p_U \eq \frac{1}{|U|}\int_U p(\b x)\,\cst{d}\b x \quad \text{and} \quad q_U \eq \frac{1}{|U|}\int_U q(\b x)\,\cst{d}\b x.
\eae 
Then, for all $k \in \NN$, and all $v \in \PP^k(U)$,
\bae\label{eq:disc:leb:emb}
\|v\|_{L^{p(\cdot)}(U)} \eqsim h_U^{\frac{d}{p_U}-\frac{d}{q_U}}\|v\|_{L^{q(\cdot)}(U)},
\eae
where the hidden constant depends only on $\Om$, $k$, $d$, $p_\cst{log}$, $q_\cst{log}$, $p_\petp$, $q_\petp$, and $C_\cst{le}$.
\ec 

\bp
We use the mean-exponent equivalence \eqref{eq:mean.exp.norm} together with the Lebesgue embedding on fixed exponents \eqref{eq:lebesgue.emb:fixed} to write for all $v \in \PP^k(U)$,
\bae
\|v\|_{L^{p(\cdot)}(U)} \eqsim \|v\|_{L^{p_U}(U)} \eqsim h_U^{\frac{d}{p_U}-\frac{d}{q_U}}\|v\|_{L^{q_U}(U)}\eqsim h_U^{\frac{d}{p_U}-\frac{d}{q_U}}\|v\|_{L^{q(\cdot)}(U)}.
\eae
\ep

\begingroup\emergencystretch=2em\sloppy

\bibliographystyle{hho_doi_or_url}
\bibliography{hhollve}
\endgroup

\end{document}